\documentclass{amsproc}

\usepackage{amssymb}

\usepackage{amsxtra,amssymb,amsmath,amscd,url,listings,stmaryrd}
\usepackage[alphabetic,initials]{amsrefs}
\usepackage{scalerel,stackengine}
\stackMath
\usepackage[utf8]{inputenc}
\usepackage{eucal}
\usepackage{fullpage}
\usepackage{scrtime}
\usepackage{hyperref}

\usepackage{datetime}
\usepackage[all]{xy}
\usepackage{graphicx}

\shortdate

\renewcommand{\leq}{\leqslant}
\renewcommand{\geq}{\geqslant}
\newcommand\widecheck[1]{%
\savestack{\tmpbox}{\stretchto{%
  \scaleto{%
    \scalerel*[\widthof{\ensuremath{#1}}]{\kern-.4pt\bigwedge\kern-.4pt}%
    {\rule[-\textheight/2]{1ex}{\textheight}}
  }{\textheight}%
}{0.5ex}}%
\stackon[2pt]{#1}{\scalebox{-1}{\tmpbox}}%
}
\newcommand{\fourier}[1]{\widehat{{#1}}}

 \newcommand{\bessel}[1]{\widecheck{{#1}}}

\numberwithin{equation}{section}

\newcommand{\uple}[1]{\text{\boldmath${#1}$}}

\def\stacksum#1#2{{\stackrel{{\scriptstyle #1}}
{{\scriptstyle #2}}}}

\def\map#1#2#3#4{\begin{matrix}#1&\to &#2
\\#3 &\to &#4
\end{matrix}}

\newcommand{\bfb}{\uple{b}}
\newcommand{\bfa}{\uple{a}}

\newcommand{\bfl}{\uple{l}}
 
\newcommand{\bfx}{\uple{x}}
\newcommand{\bfy}{\uple{y}}

\newcommand{\CF}{C(\mcF)}
\newcommand{\bfI}{\uple{I}}
\newcommand{\bfJ}{\uple{J}}
\newcommand{\eq}{\mathrm{e}_q}
\newcommand{\Std}{\mathrm{Std}}
\newcommand{\Sym}{\mathrm{Sym}}

\newcommand{\Irr}{\mathrm{Irr}}
\newcommand{\geom}{\mathrm{geom}}
\newcommand{\Garith}{G^{\mathrm{arith}}}
\newcommand{\Garithn}{G_n^{\mathrm{arith}}}
\newcommand{\Ggeom}{G^{\mathrm{geom}}}
\newcommand{\GarithF}{G_{\mcF,\mathrm{arith}}}
\newcommand{\GgeomF}{G_{\mcF,\mathrm{geom}}}
\newcommand{\Garithd}[1]{G_{{#1},\mathrm{arith}}}
\newcommand{\Ggeomd}[1]{G_{{#1},\mathrm{geom}}}
\newcommand{\Ksep}{K^{\mathrm{sep}}}

\newcommand{\bfR}{\mathbf{R}}

\newcommand{\bfK}{\mathbf{K}}
\newcommand{\bfM}{\mathbf{M}}

\newcommand{\Cc}{\mathbf{C}}
\newcommand{\Ss}{\mathbf{S}}
\newcommand{\Ct}{\mathbf{C}^\times}
\newcommand{\Nn}{\mathbf{N}}
\newcommand{\Aa}{\mathbf{A}}

\newcommand{\GmFq}{\mathbf{G}_{m,\Fq}}

\newcommand{\Zz}{\mathbf{Z}}
\newcommand{\Pp}{\mathbf{P}}
\newcommand{\Rr}{\mathbf{R}}

\newcommand{\qde}{q^{1/2}}
\newcommand{\qmd}{q^{-1/2}}

\newcommand{\Qq}{\mathbf{Q}}
\newcommand{\Ql}{\mathbf{Q}_{\ell}}
\newcommand{\ovQl}{\ov{\mathbf{Q}_{\ell}}}
\newcommand{\ovFq}{\ov{\Fq}}
\newcommand{\Fp}{{\mathbf{F}_p}}

\newcommand{\Fq}{{\mathbf{F}_q}}
\newcommand{\Fqn}{{\mathbf{F}_{q^n}}}
\newcommand{\Fqnt}{{\mathbf{F}^\times_{q^n}}}
\newcommand{\Fqt}{{\mathbf{F}^\times_q}}

\newcommand{\Ff}{\mathbf{F}}

\newcommand{\mcO}{\mathcal{O}}

\newcommand{\KL}{\mathcal{K}\ell}

\newcommand{\mods}[1]{\,(\mathrm{mod}\,{#1})}

\newcommand{\Id}{\mathrm{Id}}

\newcommand{\what}{\widehat}

\DeclareMathOperator{\frob}{Fr}

\newcommand{\tsum}{\mathcal{S}}

\newcommand{\ra}{\rightarrow}

\DeclareMathOperator{\Spec}{Spec}

\DeclareMathOperator{\rk}{rk}

\DeclareMathOperator{\Frob}{\mathrm{Frob}}
\DeclareMathOperator{\Fr}{\mathrm{Frob}}
\DeclareMathOperator{\Kl}{\mathrm{Kl}}

\DeclareMathOperator{\tr}{\mathrm{tr}}
\DeclareMathOperator{\nr}{\mathrm{Nr}}
\DeclareMathOperator{\Gal}{Gal}

\DeclareMathOperator{\Hom}{Hom}

\DeclareMathOperator{\Aut}{Aut}

\DeclareMathOperator{\swan}{Swan}

\newcommand{\eps}{\varepsilon}
\renewcommand{\rho}{\varrho}

\DeclareMathOperator{\SL}{SL}

\DeclareMathOperator{\GL}{GL}
\DeclareMathOperator{\PGL}{PGL}
\DeclareMathOperator{\PGLd}{PGL_2}

\DeclareMathOperator{\Sp}{Sp}

\DeclareMathOperator{\SU}{SU}

\newcommand{\demi}{{\textstyle{\frac{1}{2}}}}

\newcommand{\sheaf}[1]{\mathcal{{#1}}}

\DeclareMathSymbol{\gena}{\mathord}{letters}{"3C}
\DeclareMathSymbol{\genb}{\mathord}{letters}{"3E}

\newtheorem{theorem}{Theorem}[section]
\newtheorem{lemma}[theorem]{Lemma}
\newtheorem{corollary}[theorem]{Corollary}
\newtheorem{proposition}[theorem]{Proposition}
\newtheorem{conjecture}[theorem]{Conjecture}

\theoremstyle{definition}
\newtheorem{definition}[theorem]{Definition}
\newtheorem{example}[theorem]{Example}
\newtheorem{xca}[theorem]{Exercise}

\theoremstyle{remark}
\newtheorem{remark}[theorem]{Remark}

\numberwithin{equation}{section}

\newcommand{\mcM}{\mathcal{M}}
\newcommand{\mcL}{\mathcal{L}}
\newcommand{\mcS}{\mathcal{S}}
\newcommand{\mcC}{\mathcal{C}}
\newcommand{\mcP}{\mathcal{P}}
\newcommand{\mcF}{\mathcal{F}}
\newcommand{\mcQ}{\mathcal{Q}}
\newcommand{\mcK}{\mathcal{K}}

\newcommand{\mcI}{\mathcal{I}}
\newcommand{\mcG}{\mathcal{G}}
\newcommand{\mcB}{\mathcal{B}}
\newcommand{\mcE}{\mathcal{E}}

\newcommand{\mfa}{\mathfrak{a}}
\newcommand{\mfp}{\mathfrak{p}}

\newcommand{\lf}{\lambda_f}
\newcommand{\rf}{\rho_f}
\newcommand{\lamg}{\lambda_g}
\newcommand{\rhg}{\rho_g}

\newcommand{\vphi}{\varphi}

\renewcommand{\geq}{\geqslant}
\renewcommand{\leq}{\leqslant}

\newcommand{\refs}{\eqref}

\newcommand{\ov}[1]{\overline{#1}}

\newcommand{\peter}[1]{\langle{#1}\rangle}
\newcommand\sumsum{\mathop{\sum\sum}\limits}

\begin{document}

\title[]{Lectures on Applied $\ell$-adic Cohomology}

\author{\'Etienne Fouvry}
\address{Laboratoire de Math\'ematiques d'Orsay, Universit\'e Paris-Sud, 
  CNRS, Universit\'e Paris-Saclay, \linebreak[4]  91405 Orsay, France}
\email{etienne.fouvry@u-psud.fr}

\author{Emmanuel Kowalski}
\address{ETH Z\"urich -- D-MATH\\
  R\"amistrasse 101\\
  CH-8092 Z\"urich\\
  Switzerland} \email{kowalski@math.ethz.ch}

\author{Philippe Michel} \address{EPFL/SB/TAN, Station 8, CH-1015
  Lausanne, Switzerland } \email{philippe.michel@epfl.ch}

\author{Will Sawin }\address{Mathematics Department,
Rm 411, MC 4439 2990 Broadway
New York NY 10027, Columbia University, USA }
\email{will.sawin@columbia.edu}

\subjclass[2000]{Primary }

\date{\today}

\begin{abstract} We describe how a systematic use of the deep methods from $\ell$-adic cohomology pioneered by Grothendieck and Deligne and further developed by Katz and Laumon help make progress on various classical questions from analytic number theory. This text is an extended version of a series of lectures given during the 2016 Arizona Winter School.
 \end{abstract}

\maketitle

\setcounter{tocdepth}{1}
\tableofcontents

\tableofcontents

\section{Introduction}
One of the most basic question in number theory is to understand how various  sets of integers behave when restricted to (i.e.~intersected with) {\em congruence classes}, a notion that goes back at least to Euclid and was exposed systematically by Gauss in his 1801 {\em Disquisitiones Arithmeticae} (following works of Fermat, Euler, Wilson, Lagrange, Legendre and their predecessors from the middle ages and antiquity), and which is fundamental to number theory.

Let us recall that given an integer $q\in\Zz-\{0\}$, a {\em congruence class} (a.k.a. an {\em arithmetic progression}) modulo  $q$ is a subset of $\Zz$ of the shape
$$a\mods q=a+q\Zz\subset \Zz$$
for some integer $a$. The set of congruence classes modulo $q$ is denoted $\Zz/q\Zz$; it is a finite ring of cardinality $q$ (with addition  and multiplication induced by that of $\Zz$).

In number theory, especially analytic number theory, one is interested in studying the behaviour of some given arithmetic function along congruence classes, for instance to determine whether a set of integers has finite or infinite intersection with some congruence class. The analysis of such  problem, which may involve quite sophisticated manipulations, often involves certain specific classes of functions on $\Zz/q\Zz$. 

When studying such functions, it is natural to invoke the 
 {\em Chinese Remainder Theorem}
$$\Zz/q\Zz\simeq \prod_{p^\alpha\|q}\Zz/p^\alpha\Zz$$
which largely  reduces the study to the case of prime power moduli; then, in many instances, the deepest case is when $q$ is a prime; the ring $\Zz/q\Zz$ is then a finite field, denoted $\Fq$, and often the functions that occur are what we will call {\em trace functions.}

The objective of these lectures is utilitarian: our aim is to describe these trace functions, many examples, their theory and most importantly how they are handled when they occur in analytic number theory. Indeed the mention of "\'etale" or "$\ell$-adic cohomology", "sheaves", "purity", "functors", "local systems" or "vanishing cycles" sounds  forbidding to the working analytic number theorist and often prevents him/her to embrace the subject and make full use of the powerful methods that Deligne, Katz, Laumon have developed for us. It is our hope that after these introductory lectures, any of the remaining readers will feel ready for and at ease with more serious activities such as the reading of the wonderful series of orange books by Katz, and eventually will be able to tackle by him/herself any trace function that nature has laid in front of him/her.


\subsection*{Acknowledgements} These expository notes are an expanded version of a series of lectures given by Ph.M. and W.S. during the 2016 Arizona Winter School and based on our recent joint works.

We would like to thank the audience for its attention and its numerous questions during the daily lectures, as well as the teams of student, who  engaged in the research activities that we proposed during the evening sessions, for their enthusiasm. Big thanks are also due to Alina Bucur, Bryden Cais and David Zureick-Brown for the perfect organization, making this edition of the AWS a memorable experience. We would also like to thank the referees for correcting many mistakes and typosin earlier versions of this text.

\section{Examples of trace functions}

Unless stated otherwise, we now assume that $q$ is a prime number.

\subsection{Characters}
{\em Trace functions modulo $q$} are special classes of $\Cc$-valued functions on $\Fq$ of geometric origin. Perhaps the first significant example, beyond the constant function $1$, is the {\em Legendre symbol} (for $q\geq 3$)
$$\left(\frac{\cdot}q\right): x\in\Fq\to \begin{cases}0& \hbox{ if $x=0$}\\
+1&\hbox{ if $x\in(\Fqt)^2$}\\
 	-1&\hbox{ if $x\in\Fqt-(\Fqt)^2$}
 \end{cases}
$$
which detects the squares modulo $q$, and whose arithmetic properties (especially the {\em quadratic reciprocity law}) were studied by Gauss in the {\em Disquisitiones}.

The class of trace functions was further enriched by P. G. Dirichlet: on his way to  proving his famous theorem on primes in arithmetic progressions, he introduced what are  now called {\em Dirichlet characters}, i.e.~the homomorphisms of the multiplicative group
$$\chi: (\Zz/q\Zz)^\times\to \Ct$$
(with $\chi(0)$ defined to be $0$ for $\chi$ non-trivial).

Another significant class of trace functions are the additive characters
$$\psi:(\Zz/q\Zz,+)\to \Ct.$$
These are all of the shape
$$x\in\Zz/q\Zz\mapsto \eq(ax):=\exp\left(2\pi i\frac{\tilde a\tilde x}q\right)$$
(say) for some $a\in\Zz/q\Zz$, where $\tilde a$ and $\tilde x$ denote elements (lifts) of the congruence classes $a\mods q$ and $x\mods q$. Both additive and multiplicative  characters satisfy the important {\em orthogonality relations}
$$\frac{1}{q}\sum_{x\in\Fq}\psi(x)\ov{\psi'(x)}=\delta_{\psi=\psi'},\ \frac{1}{q-1}\sum_{x\in\Fqt}\chi(x)\ov{\chi'(x)}=\delta_{\chi= \chi'};$$
and we will see later a generalization of these relations to arbitrary trace functions.

Additive and multiplicative characters can be combined together (by means of a Fourier transform) to form the (normalized) {\em Gauss sums}
$$\eps_\chi(a)=\frac{1}{q^{1/2}}\sum_{x\in\Fqt}\chi(x)\eq(ax),$$
but these are not really new functions of $a$: by a simple change of variable, one has 
$$\eps_\chi(a)=\ov{\chi}(a)\eps_\chi(1)$$
for $a\in\Fqt$. For $\chi$ non-trivial,  Gauss proved that
$$|\eps_\chi(1)|=1.$$

\subsection{Algebraic exponential sums}
Another important source of trace functions comes from the study of the diophantine equations
\begin{equation}\label{eqdiophantine}
Q(\bfx)=0,\ \bfx=(x_1,\ldots,x_n)\in\Zz^n,\ Q(X_1,\ldots,X_n)\in\Zz[X_1,\ldots,X_n].
\end{equation}
 
For instance, the analysis of the {\em major arcs} in the  {\em circle method} of Hardy--Littlewood (cf.~\cite[Chap. 4]{Vaughan}) leads to the following algebraic exponential sums on $(\Zz/q\Zz)^n$ obtained by Fourier transform
$$(a,\bfx)\in (\Zz/q\Zz)^{n+1}\mapsto \frac{1}{q^{n/2}}\sum_{\bfy\in(\Zz/q\Zz)^n}\eq(aQ(\bfy)+\bfx.\bfy).$$
In the 1926's, while studying the case of a positive definite homogeneous polynomial $Q$ of degree $2$ in four variables (a positive definite integral quaternary quadratic form),  and introducing  a new  variant of the circle method, Kloosterman \cite{Kloost},  defined the so-called (normalized) {\em Kloosterman sums}
$$\Kl_2(a;q)=\frac{1}{q^{1/2}}\sum_\stacksum{x,y\in\Fqt}{xy=a}\eq(x+y).$$
This is another example of a trace function, and indeed one that is defined via  Fourier transform.

By  computing their fourth moment  (see \cite[(4.26)]{IwTopics}), Kloosterman was able to obtain the first non-trivial bound for Kloosterman sums, namely
$$|\Kl_2(a;q)|\leq 2q^{1/4}.$$ This estimate proved crucial for the study of equation \eqref{eqdiophantine} in the case of quaternary positive definite quadratic forms. In the 1940's, this bound was improved  by A. Weil, who as a consequence of his proof of the Riemann hypothesis for curves over finite fields  proved the best individual upper bound (see \cite[\S 11.7]{IwKo}):
$$|\Kl_2(a;q)|\leq 2.$$

 In 1939, Kloosterman sums appeared again in the work of Petersson who related them to Fourier coefficients of modular forms.\footnote{In fact, Poincar\'e had already written them down in one of his last papers, published posthumously.} Since then, via the works of Selberg, Kuznetsov, Deshouillers-Iwaniec and many others, Kloosterman sums play a fundamental role in the analytic theory of automorphic forms\footnote{The double occurence of Kloosterman sums in the context of quadratic forms and of modular forms is explained by the theta correspondence}.
 
 A further important example of trace functions are the (normalized) {\em hyper-Kloosterman sums}. These are higher dimensional generalisations of Kloosterman sums, and are given, for any integer $k\geq 1$ by
 $$\Kl_k(a;q)=\frac{1}{q^{(k-1)/2}}\sum_\stacksum{x_1,\ldots, x_{k}\in\Fqt}{x_1.x_2.\ldots.x_k=a}\eq(x_1+x_2+\ldots +x_k).$$
Hyper-Kloosterman sums were introduced by P. Deligne, who also established the following generalization of the Weil bound:
 $$|\Kl_k(a;q)|\leq k.$$
 
Hyper-Kloosterman sums can be interpreted as inverse (discrete)  Mellin transforms of powers of Gauss sums, and therefore can be used to study the distribution of Gauss sums. As was denoted by Katz in \cite{Sommes}, this fact and Deligne's bound imply the following\footnote{See \cite{KatzConvol} for a considerable generalisation of this theorem.}
 \begin{theorem} As $q\ra\infty$, the set of (normalized)  Gauss sums
 $$\{\eps_\chi(1),\ \chi\mods q\hbox{ non trivial }\}$$
 become equidistributed on the unit circle  $\Ss^1\subset \Ct$ with respect to the uniform (Haar) probability measure.
 \end{theorem}
Hyper-Kloosterman sums also occur in the theory of automorphic forms; for instance, Luo, Rudnick and Sarnak used the fact that powers of Gauss sums occur in the root number of the functional equation of certain automorphic $L$-functions, the inverse Mellin transform property and Deligne's bound, to obtain non-trivial estimates for the Langlands parameters of automorphic representations on $\GL_n$ (giving in particular the first improvement of Selberg's famous $3/16$ bound for the Laplace eigenvalues of Maass cusp forms).

In addition, just as for the classical Kloosterman sums, hyper-Kloosterman sums also occur in the spectral theory of $\GL_k$ automorphic forms.

There are many more examples of trace functions, and we will describe some below along with ways to construct new trace functions from older ones.

\section{Trace functions and Galois representations}

Let $\Pp^1_{\Fq}$ be the projective line and $\Aa^1_{\Fq}\subset \Pp^1_{\Fq}$ be the affine line and $K=\Fq(X)$ be the field of functions of $\Pp^1_{\Fq}$. 

In the sequel we fix some prime $\ell\not=q$, $\ovQl$ an algebraic closure of the field of $\ell$-adic numbers $\Ql$ and an embedding $\iota\colon\ovQl\hookrightarrow \Cc$  into the complex numbers. Trace functions modulo $q$ are $\ov\Ql$-valued functions\footnote{Hence $\Cc$-valued via the fixed embedding $\iota$} defined on the set of $\Fq$-points of the affine line $\Aa^1(\Fq)\simeq \Fq$. They are obtained from {\em constructible} $\ell$-adic sheaves (often denoted $\mcF$) for the \'etale topology on $\Pp^1_{\Fq}$. All these notions are quite forbidding at first; fortunately the category of {\em constructible} $\ell$-adic sheaves on $\Pp^1_{\Fq}$ can be rather conveniently described in terms of the category of representations of the Galois group of $K$. Following \cite{Sommes,GKM}, we will start from this viewpoint.

 Let $\Ksep\supset K$ be a separable closure of $K$,  and $\ov\eta$ the associated geometric generic point (i.e.~$\Spec(\Ksep)=\ov\eta$). Let $\ov\Fq\subset \Ksep$ denote the separable (or algebraic) closure of $\Fq$ in $\Ksep$. We denote
$$\Ggeom:=\Gal(\Ksep/\ov\Fq.K)\subset \Garith=\Gal(\Ksep/K), $$
the {\em geometric}, resp.~{\em arithmetic}, Galois group. By restricting the action of an element of $\Garith$ to $\ov\Fq$ we have the exact sequence
\begin{equation}
1\rightarrow \Ggeom\rightarrow \Garith\rightarrow \Gal(\ov\Fq/\Fq)\rightarrow 1.\label{exactgalois}
\end{equation}

\begin{definition} Let $U\subset \Aa^1_\Fq$ be a non-empty open subset of $\Aa^1_\Fq$ that is defined over $\Fq$. An $\ell$-adic sheaf lisse on $U$, say $\mcF$, is a continuous finite-dimensional Galois representation
$$\rho_\mcF:\Garith\to \GL(V_\mcF)$$
where $V_\mcF$ is a finite dimensional $\ovQl$-vector space, which is unramified at every closed point $x$ of $U$.
The dimension $\dim V_\mcF$ is called the rank of $\mcF$ and is denoted $\rk(\mcF)$. The vector space $V_{\mcF}$ is also denoted $\mcF_{\ov\eta}$.
\end{definition}

\subsection{Closed points on the affine line}

In this section we spell-out the meaning of the sentence "unramified at every closed point $x$ of $U$".

Let us recall that the datum of closed point of $\Pp^1_{\Fq}$ is equivalent to the datum of an embedding $\mcO_x\hookrightarrow K$ of a local ring\footnote{A PID with a unique prime ideal \cite[Chap. 1]{Serre}} $\mcO_x$ (the ring of rational functions defined in a neighborhood of $x$) whose field of fractions is $K$. Given such an embedding, we denote by $\mfp_x$ its unique prime ideal, $\pi_x$ a generator of $\pi_x$ (an uniformizer) and by $v_x:K\ra\Zz\cup\{\infty\}$ the associated discrete valuation (normalized so that $v_x(\pi_x)=1$): we have
$$\mcO_x=\{f\in K,\ v_x(f)\geq 0\}\supset \mfp_x=\{f\in K,\ v_x(f)> 0\}.$$
We denote by $k_x=\mcO_x/\mfp_x$ its residue field and by $q_x=|k_x|=:q^{\deg x}$ the size of $k_x$ and $\deg x$ its degree

 The set of closed points of the projective line $\Pp^1_{\Fq}$ is the union of the set of closed points of the affine line $\Aa^1_{\Fq}$ which is indexed by the set of monic, irreducible (non-constant) polynomials of $\Fq[X]$
and the point $\infty.$ 
\begin{itemize}
\item For $\pi$ irreducible, monic and not constant, the local ring $\mcO_\pi$ is the localization of $\Fq[X]$ at the prime ideal $(\pi) \subseteq \Fq[X]$:
$$\mcO_\pi=\{P/Q,\ P,Q\in\Fp[X],\ \pi\!\!\not|Q\}\supset \mfp_\pi=\{P/Q,\ P,Q\in\Fp[X],\ \pi|P,\ \pi\!\!\not|Q\},$$
the valuation $v_\pi$ is the usual valuation: for any  polynomial $P\in \Fq[X]$, $v_x(P)=v_\pi(P)$ is the exponent of the highest power of $\pi$ dividing $P$ which is extended to $K$ by setting $v_x(P/Q)=v_\pi(P)-v_\pi(Q)$, and the degree is $\deg \pi$.
\item For $\infty$, 
$$\mcO_\infty=\{P/Q,\ P,Q\in\Fp[X],\ \deg P\leq\deg Q\}\supset \mfp_\infty=\{P/Q,\ P,Q\in\Fp[X],\ \deg P<\deg Q\},$$
the valuation is minus the degree of the rational fraction
 $v_\infty(P/Q)=\deg(Q)- \deg(P),$ and the degree of $\infty$ is $1$.
\end{itemize}

\begin{remark}\label{closedpointsident}
We denote by $\Pp^1(\Fq)$ the set of closed points of degree $1$ and by $\Aa^1(\Fq)=\Pp^1(\Fq)-\{\infty\}$. Note that $\Aa^1(\Fq)$ is identified with $\Fq$ by identifying $x\in\Fq$ with the degree $1$ (irreducible) polynomial $X-x$.

Similarly a non-empty open set $U\subset \Aa^1_{\Fq}$ is the open complement of the closed set $Z_Q\subset \Aa^1_{\Fq}$ of zeros of some (non-zero) polynomial $Q\in\Fq[X]$, i.e.~defined by the equation $Q(x)=0$. 

The "closed points of $U$" are the closed point associated with the irreducible monic  polynomials $\pi\in\Fq[X]$ coprime to $Q$ and the set of closed points of degree $1$, is identified with the complement of the set of roots of $Q$ contained in $\Fq$:
$$U(\Fq)\simeq \{x\in\Fq,\ Q(x)\not=0\}\subset\Fq.$$
\end{remark}

\subsubsection{Decomposition group, inertia and Frobenius}
 The valuation $v_x$ can be extended (in multiple ways) to a ($\Qq$-valued) valuation on $\Ksep$ and the choice of one such extension (denoted $v_{\{x\}}$) determines a decomposition and an inertia subgroup in the arithmetic Galois group
 $$I_{\{x\}}\subset D_{\{x\}}\subset \Garith$$ fitting in the exact sequence
\begin{equation}\label{eqlocalgal}
1\rightarrow I_{\{x\}}\rightarrow D_{\{x\}}\rightarrow \Gal(\ov\Fq/k_x)\rightarrow 1.	
\end{equation}
Let also us recall that $\Gal(\ov\Fq/k_x)$ is topologically generated by the {\em arithmetic Frobenius} 
$$\Frob^{\mathrm{arith}}_{k_x}\colon\map{\ov\Fq}{\ov\Fq}{u}{u^{q_x}}.$$ In the sequel we will denote by $\Frob^{\mathrm{geom}}_{k_x}$ its inverse, also called the {\em geometric Frobenius}. The lifts of the (geometric) Frobenius therefore define a (left) $I_{\{x\}}$-class in the decomposition subgroup which we denote by $$\Fr_{\{x\}}\subset D_{\{x\}}$$ and
which we call the Frobenius class at $\{x\}$.

\begin{remark}\label{abuse} The choice of a different extension $v_{\{x\}'}$ of $v_x$ yields a priori another decomposition, inertia subgroups and Frobenius class, $D_{\{x\}'}, I_{\{x\}'},\ Fr_{\{x\}'}$, but these are conjugate to $D_{\{x\}}, I_{\{x\}},\ Fr_{\{x\}}$ because $\Garith$ acts transitively on the set of extensions. As we will see the various quantities that we will discuss in relation to these sets will be conjugacy-invariant and therefore depend only on $x$ but not of a choice of $\{x\}$ and will use the indice $x$ instead of $\{x\}$. Sometimes, to simplify notations, we will implicitly assume the choice of an $\{x\}$ without mentioning it and will simply write $D_x, I_x, \Fr_x$ 
\end{remark}

We can now explain the term unramified.

\begin{definition} Given $x$ a closed point of $\Pp^1_{\Fq}$, a $\Garith$-module $V$ is unramified (or lisse) at $x$ at if for one (or equivalently any) extension $\{x\}$, the corresponding inertia subgroup $I_{\{x\}}$ acts trivially  on $V$. Otherwise $V$ is ramified at $x$.

If $V$ is unramified at $x$, all the elements in the Frobenius class $\Fr_{\{x\}}$
 act by the same automorphism of $V$ and we will denote this automorphism by $(\Fr_{\{x\}}|V)$. 
 
 Moreover if we change the extension $\{x\}$ we obtain an automorphism which is $\Garith$-conjugate to $(\Fr_{\{x\}}|V)$. We denote by $(\Fr_x|V_\mcF)$ this conjugacy class.

 \end{definition}
 
 It follows from this discusion that for any sheaf $\mcF$ there is a non-empty open subset on which $\mcF$ is unramified and maximal for this property. We will note this open set $U_\mcF$.

\subsection{The trace function attached to a lisse sheaf}

Let $\mcF$ be an $\ell$-adic sheaf lisse on $U\subset\Aa^1_{\Fq}$ and $$\rho_\mcF\colon \Garith\to \GL(V_\mcF)$$ the corresponding representation.

For $x\in U(\Fq)$ a closed point of degree $1$ at which the representation $\rho_\mcF$ is unramified, we have, in the previous section, associated a Frobenius conjugacy class $(\Fr_{x}|V_\mcF)$
namely the union of all the $(\Fr_{\{x\}}|V_\mcF)$. By conjugacy, the trace of all these automorphisms $(\Fr_{\{x\}}|V_\mcF)$
 is constant within that class: we denote this common value by
 $$\tr(\Fr_{x}|V_\mcF)$$
 and call it the Frobenius trace of $\mcF$ at $x$.

\begin{definition}
Given an $\ell$-adic sheaf $\mcF$ lisse on $U\subset \Aa^1_{\Fq}$; the trace function $K_\mcF$ associated to this situation is the function on $U(\Fq)$ given by
$$x\in U(\Fq)\mapsto K_\mcF(x)=\tr(\Frob_x|V_\mcF).$$
This is a priori a $\ovQl$-valued function but it can be considered complex-valued via the fixed embedding $\iota\colon\ovQl\hookrightarrow \Cc$.
\end{definition}

\begin{remark}As we have seen in Remark \ref{closedpointsident} $U(\Fq)$ is identified with
$$\{x\in\Fq,\ Q(x)\not=0\}\subset \Fq$$
and therefore $K_\mcF$ can be considered as a function defined on a subset of $\Fq$.
\end{remark}

\begin{remark}\label{remextension} There are several ways by which one could extend $K_\mcF$ to the whole of $\Aa^1(\Fq)$. The simplest way is the extension by zero outside $U(\Fq)$; another possible extension (called the {\em middle extension}) would be to set for any $x\in\Aa^1(\Fq)$,
$$K_\mcF(x):=\tr(\Frob_{\{x\}}|V_\mcF^{I_{\{x\}}})$$
where $V_\mcF^{I_{\{x\}}}\subset V_\mcF$ is the subspace of $I_{\{x\}}$-invariant vectors: the action of the Frobenius class $\Frob_{\{x\}}$ on $V_\mcF^{I_{\{x\}}}$ is well-defined and its trace does not depend on $\{x\}$. For our purpose, any of the two extensions would work (cf.~ Remark \ref{remdeligneext}).
\end{remark}

\subsection{Trace functions over $U(\Fqn)$}
In fact, an $\ell$-adic sheaf, lisse on $U_\Fq$ give rise to a whole family of trace functions.

For any $n\geq 1$, let us consider the finite extension $\Ff_{q^n}$ let us and base change the whole situation to that field: this amounts to replace $\Pp^1_{\Fq}$ by $\Pp^1_{\Fqn}$, $K=\Fq(X)$ by $K_n=\Fqn(X)$, and the arithmetic Galois group $\Garith$ by  $\Garithn=\Gal(\Ksep/K_n)$ (notice that the geometric Galois group does not change). 

The group $\Garithn$ is a normal subgroup of $\Garith$ (whose quotient is $\Gal(\Fqn/\Fq)$, so we may restrict our initial Galois representation to it: in that way we obtain another $\ell$-adic sheaf denoted $\mcF_n$
 $$\rho_{\mcF_n}\colon\Garithn\ra \GL(V_\mcF)$$
 and another trace function
$$K_{\mcF,n}\colon\map{U(\Fqn)}{\Cc}{x}{\tr(\Fr_x|V_\mcF)}$$
where $U(\Fqn)$ denotes now the set of closed points of $\Pp^1_{\Ff_{q^n}}$ of degree $1$ which are contained in $U$: this set is identified with the set of irreducible monic polynomials of degree $1$ coprime with $Q$ and is therefore identified with
$$\{x\in\Ff_{q^n},\ Q(x)\not=0\}.$$  As we will see below, the existence of this sequence of auxiliary functions is very important: for instance (the Chebotareff density theorem) the full sequence $(K_{\mcF,n})_{n\geq 1}$ characterizes the representation $\rho_\mcF$ up to semi-simplification.

\begin{remark}\label{remwarning} As we have remarked already one has the identifications
$$U(\Fq)\simeq \{x\in\Fq,\ Q(x)\not=0\},\ U(\Fqn)\simeq \{x\in\Fqn,\ Q(x)\not=0\}.$$ 
However the inclusion
$$\{x\in\Fq,\ Q(x)\not=0\}\subset \{x\in\Fqn,\ Q(x)\not=0\}$$
 does NOT imply that the function $K_\mcF$ is "the restriction" of $K_{\mcF,n}$ to $U(\Fq)$. More precisely, if we denote by 
$x$ the closed point in $U(\Fq)$ associated with the polynomial $X-x\in\Fq[X]$ and by $x_n$ the closed point in $U(\Fqn)$ associated with the same polynomial $X-x\in\Fqn[X]$ one has the formula
$$K_{\mcF,n}(x_n)=\tr(\Frob_{x_n}|V_\mcF)=\tr(\Frob^n_{x}|V_\mcF).$$

More generally, for $d$ dividing $n$ let $\pi\in\Fq[X]$ be a monic irreducible polynomial of degree $d$ and coprime to $Q$. Then $\pi$ defines a closed point $x_\pi$ of $U$ of degree $d$. Since $d|n$, the polynomial $\pi$ splits in $\Fqn$
$$\pi(X)=\prod_{i=1}^d(X-x_i)$$ and any of its roots $x_i$ defines a closed point in $U(\Fqn)$ (corresponding to the polynomial $X-x_i\in\Fqn[X]$); we then have for $i=1,\ldots,d$
\begin{equation}\label{eqwarning}
	K_{\mcF,n}(x_i)=\tr(\Frob_{x_i}|V_\mcF)=\tr(\Frob_{\pi}^{n/d}|V_\mcF).
\end{equation}

\end{remark}

\begin{remark} There is, a priori, no reason to limit ourselves to the affine line: if $\mcC_{\Fq}$ is any smooth geometrically connected curve over $\Fq$ with function field $K_\mcC$ (which is a finite extension of $\Fq(X)$) and any dense open subset  $U\subset\mcC$  defined over $\Fq$, an $\ell$-adic sheaf $\mcF$ on $\mcC$  lisse on some non-empty open set $U$ is a continuous representation$$\rho_\mcF\colon\Gal(K^{sep}_\mcC/K_\mcC)\ra \GL(V_\mcF)$$
which is unramified at every closed point of $U$. \end{remark}

\subsection{The language of representations}

The definition of sheaves and trace functions in terms of Galois representations enable to use consistently the vocabulary from representation theory. For instance

\begin{itemize}
\item An $\ell$-adic sheaf is {\em irreducible} (resp.~{\em isotypic}) if the representation $\rho_\mcF$ is.
\item An $\ell$-adic sheaf is {\em geometrically irreducible} (resp.~{\em geometrically  isotypic}) if the {\em restriction} of $\rho_\mcF$ to the {\em geometric Galois group} $\Ggeom$ is.
\item An $\ell$-adic sheaf is  {\em trivial} if the representation $\rho_\mcF$ is. The trace function is constant, equal to $1$.
\item An $\ell$-adic sheaf is  {\em geometrically trivial} if the {\em restriction} of $\rho_\mcF$ to the {\em geometric Galois group} $\Ggeom$ is. In view of \ref{exactgalois} its trace function is a constant, say $K_\mcF(x)=\alpha$ and for any $n\geq 1$, $$K_{\mcF,n}(x)=\alpha^n.$$
 \end{itemize}
 One can also create new sheaves and trace function from other sheaves.
 
 \begin{itemize}
\item The {\em dual sheaf} $D(\mcF)$ is the contragredient representation $D(\rho_\mcF)$ acting on the dual space $\Hom(V_\mcF,\ovQl)$. This   sheaf is also lisse on $U$ and its trace function is given for $x\in U(\Fq)$ by
$$K_{D(\mcF)}(x)=\tr(\Frob_x^{-1}|V_\mcF).$$
\item Given another sheaf $\mcG$ lisse on some $U'$, one can form the {\em direct sum sheaf} $\mcF\oplus\mcG$ whose representation is $\rho_{\mcF\oplus\mcG}=\rho_\mcF\oplus\rho_\mcG$; the sheaf is lisse (at least) on $U\cap U'$, of rank the sum of the ranks,  and its trace function is given, for $x\in U(\Fq)\cap U'(\Fq)$ by the sum
$$K_{\mcF\oplus\mcG}(x)=K_{\mcF}(x)+K_{\mcG}(x).$$
\item Given another sheaf $\mcG$ lisse on some $U'$, one can form the {\em tensor product sheaf} $\mcF\otimes\mcG$ whose representation is $\rho_{\mcF\otimes\mcG}=\rho_\mcF\otimes\rho_\mcG$; the sheaf is lisse (at least) on $U\cap U'$, of rank the product of the ranks, and its trace function is given, for $x\in U(\Fq)\cap U'(\Fq)$ by the product
$$K_{\mcF\otimes\mcG}(x)=K_{\mcF}(x)K_{\mcG}(x).$$
 \item As a special case, one construct the {\em sheaf of homomorphisms} between $\mcF$ and $\mcG$ and the {\em sheaf of endomorphisms} of $\mcF$, 
$$\mathrm{Hom}(\mcF,\mcG)=D(\mcF)\otimes\mcG ,\ \mathrm{End}(\mcF)=D(\mcF)\otimes \mcF.$$
\item Let $H\subset\GL(V_\mcF)$ be an algebraic group containing $\rho_\mcF(\Garith)$ and let $r\colon H\to \GL(V')$ be a finite-dimensional continuous $\ell$-adic representation; the composite representation $r\circ \rho_\mcF$ defines an $\ell$-adic sheaf, denoted $r\circ\mcF$, which is lisse on $U$ and has rank $\dim V'$. Its trace function  is given, for $x\in U(\Fq)$ by 
$$K_{r\circ \mcF}(x)=\tr(r(\Frob_x|V_\mcF)|V').$$

 \item Let $f\in\Fq(X)$ be non-constant; we can view $f$ as a non-constant morphism $\Pp^1_{\Fq}\to \Pp^1_{\Fq}$. The Galois subgroup  corresponding to this covering $$\Gal(\Ksep/\Fq(f(X)))\subset \Garith$$ is isomorphic to $\Garith$  and therefore the restriction of $\rho_{\mcF}$ to $\Gal(\Ksep/\Fq(f(X)))$ defines an $\ell$-adic sheaf on $\Pp^1_{\Fq}$ lisse on $f^{-1}(U)$ which is denoted $f^*\mcF$ and is called the {\em pull-back} of $\mcF$ by $f$. Its rank equals the rank of $\mcF$ and its trace function  is given, for $x\in f^{-1}(U)(\Fq)-\{\infty\}$ by 
 $$K_{f^*\mcF}(x)=K_\mcF(f(x)).$$ 
 
 \item If the sequel, we will use this pull-back sheaf construction for the following morphisms:  This are special cases of {\em fractional linear transformations}: given $\gamma=\begin{pmatrix}a&b\\c&d \end{pmatrix}\in\PGLd(\Fq)$ (the group of automorphisms of $\Pp^1_{\Fq}$) one defines the automorphism
$$[\gamma]\colon x\to \frac{ax+b}{cx+d}.$$
We denote the pull-back sheaf by $\gamma^*\mcF$. In particular, for $\gamma=n(b)=\begin{pmatrix}1&b\\0&1 \end{pmatrix}$ we obtain the 
the additive translation map 
 $[+b]\colon x\to x+b$, and for $\gamma=t(a)=\begin{pmatrix}a&0\\0&1 \end{pmatrix},\ a\not=0$ we obtain the multiplicative translation map $[\times a]\colon x\to ax.$

 \end{itemize}

\subsection{Purity}
We will be interested in the size of trace functions. For this we need the notion of {\em purity}.

\begin{definition} Let $w\in\Zz$. an $\ell$-adic sheaf $\mcF$, lisse on $U$ is punctually pure of weight $w$ if, for any $x\in U_{\Fq}$, the eigenvalues of $(\Frob_x|V_\mcF)$ are complex numbers\footnote{via the fixed embedding $\ov{\Ql}\hookrightarrow\Cc$.} of modulus $q_x^{w/2}$. It is mixed of weights $\leq w$ if (as a representation) it is a successive extension of sheaves punctually pure of weights $\leq w$.

In particular, if $\mcF$ is mixed of weights $\leq w$, one has for any $x\in U(\Fq)$
\begin{equation}\label{infinitybound}
|K_\mcF(x)|\leq \rk(\mcF)q^{w/2}.	
\end{equation}

\end{definition}

\begin{remark}
 It is always possible to reduce to the case of $\ell$-adic sheaves mixed of weight $w=0$. For any $w\in\Zz$ there exist an $\ell$-adic sheaf denoted $\ovQl({w/2})$ of rank $1$, lisse on $\Pp^1_\Fq$, whose restriction to $\Ggeom$ is trivial and such that $\Fr_x$ acts by multiplication by $q_x^{-w/2}$ (in particular $\ovQl({w/2})$ is pure of weight $-w$). Given $\mcF$ of some weight $w'$, the tensor product
 $$\mcF(w/2):=\mcF\otimes \ovQl({w/2})$$ 
 has weight $w'-w$ and has trace function given by
 $$x\mapsto q^{-w/2}K_\mcF(x).$$
 \end{remark}

 In the sequel, unless stated otherwise, we will always assume that trace functions are associated with sheaves which are mixed of weights $\leq 0$.

 \begin{remark}\label{remdeligneext} Deligne proved (\cite[Lemme (1.8.1)]{WeilII}) that for a sheaf punctually pure of weight $w$, for any closed point $x\in\Pp^1_\Fq$, the eigenvalues of 
$(\Frob_x|V_\mcF^{I_x})$ have modulus $\leq q_x^{w/2}$. In particular
$$|\tr(\Frob_x|V_\mcF^{I_x})|\leq \rk(\mcF)q_x^{w/2}.$$
In particular (assuming that $w=0$) $\ell^\infty$-norm of the difference between the extension by $0$ of $K_\mcF$ from $U(\Fq)$ to $\Aa^1(\Fq)$ and the middle-extension (described in Remark \ref{remextension}) is bounded by
$$\rk(\mcF)|\Aa^1(\ov\Fq)-U(\ov\Fq)|.$$
As we will see, we will be interested in situations where this quantity is  bounded by an absolute constant (independent of $q$) the consequence being that whatever extension we choose between the two, it won't make much of a difference.
\end{remark}

\subsection{Other functions} There are other  functions on $\Fq$ of great interest which do not qualify as trace functions under our current definition. For instance the Dirac function at some point $a\in\Fq$
$$\delta_a(n)=\begin{cases}1&\hbox{ if }n\equiv a\mods q\\0&\hbox{ otherwise }.	
\end{cases}
$$
which, extended to $\Zz$ is the characteristic function of the arithmetic progression $a+q\Zz$ (obviously of considerable interest for analytic number theory). It turns out that such functions can be related to trace functions in our sense by very natural transformations and this will allow us to make some progress on problems from "classical" analytic number theory.

\begin{remark} In fact this function could be interpreted as the trace function of a {\em skyscraper sheaf} supported at the closed point $a$ but we will not do this here.
\end{remark}

\subsection{Local monodromy representations}
Given $\mcF$ some $\ell$-adic sheaf, let $$D^{ram}_{\mcF}\subset \Pp^1(\ov{\Fq})-U(\ov\Fq)$$ be the set of geometric points where the representation $\rho_\mcF$ is ramified, that is the inertia group $I_x$ acts non-trivially. The restricted representation
$$\rho_{\mcF|I_x}=\rho_{\mcF,x}$$ is called the local monodromy representation of $\mcF$ at $x$ (cf.~Remark \ref{abuse} for the abuse of notation). Although $D^{ram}_{\mcF}$ is disjoint from $U(\ov\Fq)$, this finite set of representations is fundamental to study $\mcF$ and its trace function. Let us recall  the exact sequence \cite[Chap. 1]{GKM}
$$1\ra P_x\ra I_x\ra I_x^{tame}\ra 1$$
where $I_x^{tame}$ is the {\em tame inertia quotient} and is isomorphic to $\prod_{p\not= q}\Zz_p$, while $P_x$ is the $q$-Sylow subgroup of $I_x$ and is called the wild inertia subgroup.
\begin{definition} The sheaf is tamely ramified at $x$ if $P_x$ acts trivially on $V_\mcF$ (so that $\rho_{\mcF,x}$ factors through $I_x^{tame}$) and is called wildly ramified otherwise.
\end{definition}
\subsubsection{The Swan conductor}
If the representation is wildly ramified one can measure how deep it is by means of a numerical invariant: {the Swan conductor}. The wild inertia subgroup $I_x$ is equipped with the decreasing {\em  upper numbering filtration} $I_x^{(\lambda)}$ indexed by non-negative real numbers $\lambda\geq 0$, such that $$P_x=I_x^{(>0)}=\bigcup_{\lambda>0}I_x^{\lambda}.$$ Given $V=V_\mcF$ as above there is a $P_x$-stable direct sum decomposition 
$$V=\bigoplus_{\lambda\in \mathrm{Break}(V)}V(\lambda)$$ indexed by a finite set of rational numbers $ \mathrm{Break}(V)\subset \Qq_{\geq 0}$ (the set of {\em breaks} of the $I_x$-module $V$) such that 
$$V(0)=V^{P_x},\ V(\lambda)^{I_x^{(\lambda)}}=0,\ V(\lambda)^{I_x^{(\lambda')}}=V(\lambda),\ \lambda'>\lambda$$
(see \cite[Chap. 1]{GKM}). The {\em Swan conductor} is defined as
$$\swan_x(\mcF)=\sum_{\lambda\in \mathrm{Break}(V)}\lambda\dim V(\lambda)$$ and turns out to be an integer \cite[Prop. 1.9]{GKM}.

In the decomposition
$$V=V(0)\oplus \bigoplus_\stacksum{\lambda\in \mathrm{Break}(V)}{\lambda>0}V(\lambda)=V(0)\oplus V(>0):=V^{tame}\oplus V^{wild}$$
the first summand is called the {\em tame part} and the remaining one the {\em wild part}.

\section{Summing trace functions over $\Fq$}
Let $K_\mcF$ be the trace function associated to a sheaf $\mcF$ lisse on $U_\Fq$. It is a function on $U(\Fq)$ which we may extend by zero to $\Aa^1(\Fq)\simeq\Fq=\Zz/q\Zz$. 

The Grothendieck-Lefschetz trace formula provides an alternative expression for the sum of $K_\mcF$ over the whole $\Aa^1(\Fq)$.

\begin{theorem}[Grothendieck-Lefschetz trace formula] Let $\mcF$ be lisse on $U$; there exists three finite dimensional $\ell$-adic representations of $\Gal(\ovFq/\Fq)$, $H^i_c(U_{\ovFq},\mcF)$ such that
\begin{equation}\label{GLformula}
\sum_{x\in U(\Fq)}K_\mcF(x)=\sum_{x\in U(\Fq)}\tr(\frob_x|\mcF)=\sum_{i=0}^2(-1)^i\tr(\Frob_q|H^i_c(U_{\ovFq},\mcF)).
\end{equation}

More generally, for any $n\geq 1$,
$$\sum_{x\in U(\Fqn)}K_{\mcF,n}(x)=\sum_{x\in U(\Fqn)}\tr(\frob_x|\mcF)=\sum_{i=0}^2(-1)^i\tr(\Frob^n_q|H^i_c(U_{\ovFq},\mcF)).$$
\end{theorem}
The $\ovQl$-vector spaces $H^i_c(U_{\ovFq},\mcF)$ are the so-called compactly supported étale cohomology groups of $\mcF$ and can also be considered as $\ell$-adic sheaves over the point $\Spec(\Fq)$.

The above formula reduces the evaluation of averages of trace functions to that of the three summands
$$\tr(\Frob_q|H^i_c(U_{\ovFq},\mcF)),\ i=0,1,2,$$
we need therefore to control the dimension of these spaces as well as the size of the eigenvalues. We start with the former.

\subsection{Bounding the dimension of the cohomology groups}

 The extremal cohomology groups have a simple interpretation. First
$$H^0_c(U_{\ovFq},\mcF)=\begin{cases}0&\hbox{ if $U\not=\Pp^1_{\Fq}$}\\
V_\mcF^{\Ggeom}&\hbox{ if $U=\Pp^1_{\Fq}$}.	
\end{cases}
$$
As a $\Gal(\ovFq/\Fq)$-representation, one has an isomorphism
\begin{equation}\label{H2isom}
H^2_c(U_{\ovFq},\mcF)\simeq V_{\mcF,\Ggeom}(-1)	
\end{equation}
(ie $H^2_c(U_{\ovFq},\mcF)$ is isomorphic to the quotient of $\Ggeom$-coinvariants of $V_{\mcF}$ twisted by $\ovQl(-1)$). In particular, if $\mcF$ is geometrically irreducible (non geometrically trivial) or more generally geometrically isotypic (the underlying geometric irreducible representation being non trivial) one has
$$H^2_c(U_{\ovFq},\mcF)=0.$$
In any case, one has $$\dim H_c^0(U_{\ovFq},\mcF),\ \dim H_c^2(U_{\ovFq},\mcF)\leq\rk(\mcF).$$

The dimension of the middle cohomology group is now determined by the

\begin{theorem}[The Grothendieck-Ogg-Shafarevich formula] One has the following equality
$$\chi(U_{\ov\Fq},\mcF)=\sum_{i=0}^2(-1)^i\dim H^i_c(U_{\ovFq},\mcF)=\rk(\mcF)(2-|\Pp^1(\ov\Fq)-U(\ov\Fq)|)-\sum_{x\in D^{ram}_{\mcF}(\ov\Fq)}\swan_x(\mcF).$$ 
\end{theorem}

Observe that the quantities that occur are local geometric data associated to the sheaf  yet this collection of local data provides global informations.

We then define the following ad-hoc numerical invariant which serves as a measure of the complexity of the sheaf $\mcF$:
\begin{definition} The conductor of $\mcF$ is defined via the following formula
$$C(\mcF)=\rk(\mcF)+|\Pp^1(\ov\Fq)-U(\ov\Fq)|+\sum_{x\in D^{ram}_{\mcF}(\ov\Fq)}\swan_x(\mcF).$$
\end{definition}
In view of this definition we have
\begin{equation}\label{eqdimbound}
\sum_{i=0}^2\dim H^i_c(U_{\ovFq},\mcF)\ll C(\mcF)^2.
\end{equation}

\subsection{Examples}
\subsubsection{The trivial sheaf} The trivial representation $\ovQl$ is everywhere lisse, pure of weight $0$, of rank $1$ and conductor $1$ and
$$K_{\ovQl}(x)=1.$$
\subsubsection{Kummer sheaf \cite{sga4h}} For any non-trivial Dirichlet character 
$\chi\colon(\Fqt,\times)\to \Ct$ there exists an $\ell$-adic sheaf (a Kummer sheaf) denoted $\mcL_\chi$ which is of rank $1$, pure of weight $0$, lisse on $\GmFq=\Pp^1_\Fq-\{0,\infty\}$ with trace function
$$K_{\mcL_\chi}(x)=\chi(x),\ K_{\mcL_\chi,n}(x)=\chi(\nr_{\Fqn/\Fq}(x))=:\chi_n(x)$$
and conductor
$$C(\mcL_\chi)=3;$$
indeed $\swan_0(\mcL_\chi)=\swan_\infty(\mcL_\chi)=0$.

\subsubsection{Artin-Schreier sheaf \cite{sga4h}} For any additive character $\psi\colon(\Fq,+)\to \Ct$ there exists an $\ell$-adic sheaf (an Artin-Schreier sheaf) denoted $\mcL_\psi$ which is of rank $1$, pure of weight $0$, lisse on $\Aa^1_{\Fq}=\Pp^1_\Fq-\{\infty\}$ with trace function
$$K_{\mcL_\psi}(x)=\psi(x),\ K_{\mcL_\psi,n}(x)=\psi(\tr_{\Fqn/\Fq}(x))=:\psi_n(x)$$
and conductor (if $\psi$ is non-trivial)
$$C(\mcL_\psi)=3.$$
(indeed $\swan_\infty(\mcL_\psi)=1$). If $f\in\Fq(X)-\Fq$, the pull-back sheaf $\mcL_{\psi(f)}$ is geometrically irreducible and has conductor
$$1+\hbox{ number of poles of $f$}+\hbox{ sum of multiplicities of the poles of $f$}.$$
More generally any character $\psi$ of $(\Fqn,+)$ is of the shape
$$x\mapsto \psi_1(\tr_{\Fqn/\Fq}(ax))$$
for $\psi_1$ a non-trivial character of $(\Fq,+)$ and $a\in\Fqn$,  and associated to each such character is an Artin-Schreier sheaf $\mcL_\psi$.

\subsubsection{(hyper)-Kloosterman sheaves \cite{GKM}} Hyper-Kloosterman sums are formed by multiplicative convolution out of additive characters. 

Given $K_1,K_2\colon\Fqt\ra\Cc$ one defines their (normalized) multiplicative convolution:
$$K_1\star K_2\colon x\in\Fqt\to \frac{1}{q^{1/2}}\sum_\stacksum{x_1,x_2\in\Fqt}{x_1.x_2=x}K_1(x_1)K_2(x_2)=\frac{1}{q^{1/2}}\sum_{x_1\in\Fqt}K_1(x_1)K_2(x/x_1).$$
Similarly for any $n\geq 1$ one defines the multiplicative convolution of $K_{1,n},K_{2,n}\colon \Fqnt\ra\Cc$ as
$$K_{1,n}\star K_{2,n}\colon x\in\Fqnt\to \frac{1}{q^{n/2}}\sum_\stacksum{x_1,x_2\in\Fqnt}{x_1.x_2=x}K_{1,n}(x_1)K_{2,n}(x_2).$$ 
Now, given  a non-trivial additive character $\psi$ of $\Fq$ and $k\geq 2$, the hyper-Kloosterman sums can be expressed as the $k$-fold multiplicative convolutions of $\psi$:
$$\Kl_{k,\psi}(x;q)=\star_{k\hbox{ times}}\psi(x)=\frac{1}{q^{({k-1})2}}\sum_\stacksum{x_1,\ldots,x_k\in\Fqt}{x_1.\ldots.x_k=x}\psi(x_1+\ldots+x_k)$$
and more generally, one defines hyper-Kloosterman sums over $\Fqnt$
$$\Kl_{k,\psi}(x;q^n)=\star_{k\hbox{ times}}\psi_n(x)=\frac{1}{q^{n({k-1})2}}\sum_\stacksum{x_1,\ldots,x_k\in\Fqnt}{x_1.\ldots.x_k=x}\psi_n(x_1+\ldots+x_k).$$
These are in fact trace functions: their underlying sheaves were constructed by Deligne and were subsequently studies in depth by Katz \cite{GKM}:
\begin{theorem} For any $k\geq 2$, there exists an $\ell$-adic sheaf (the Kloosterman sheaf) denoted $\KL_{k,\psi}$, of rank $k$, pure of weight $0$, geometrically irreducible, lisse on $\GmFq$ with trace function
$$K_{\KL_{k,\psi}}(x)=\Kl_{k,\psi}(x;q)$$
and more generally, for any $n\geq 1$
$$K_{\KL_{k,\psi},n}(x)=\Kl_{k,\psi}(x;q^n).$$
One has $\swan_0(\KL_{k,\psi})=0$ and $\swan_\infty(\KL_{k,\psi})=1$ so that the conductor of that sheaf equals
$$C(\KL_{k,\psi})=k+2+1.$$
The Kloosterman sheaves have trivial determinant
$$\det\KL_k=\ovQl$$
and if (and only if) $k$ is even, the Kloosterman sheaf $\KL_k$ is self-dual:
$$D(\KL_k)\simeq \KL_k.$$
\end{theorem}
\begin{remark}
	When $\psi(\cdot)=\eq(\cdot)$ we will not mention the additive character $\eq$ in the notation.
\end{remark}

\subsection{Deligne's Theorem on the weight}
Now that we control the dimension of the cohomology groups occurring in the Grothendieck-Lefschetz trace formula, it remains to control the size of their Frobenius eigenvalues. Suppose that $\mcF$ is pure of weight $0$ so that
$$|K_\mcF(x)|\leq \rk(\mcF).$$
As we have seen, as long as $U\not=\Pp^1$, $H^0_c(U_{\ovFq},\mcF)=0$.

By \eqref{H2isom}, the eigenvalues of $\Frob_q$ acting on $H^2_c(U_{\ovFq},\mcF)$ are of the form $$q\alpha_i,\ i=1,\ldots,\dim(V_{\mcF,\Ggeom})\hbox{ with }|\alpha_i|=1.$$
The trace of the Frobenius on the middle cohomology group $\tr(\Frob_q|H^1_c(U_{\ovFq},\mcF))$ is much more mysterious but fortunately we have the following theorem of Deligne \cite{WeilII}.
\begin{theorem}[The Generalized Riemann Hypothesis for finite fields]\label{delignethm} The eigenvalues of $\Frob_q$ acting on $H^1_c(U_{\ovFq},\mcF)$ are complex numbers of modulus $\leq q^{1/2}$.
\end{theorem}
We deduce from this
\begin{corollary}\label{delignecor} Let $\mcF$ be an $\ell$-adic sheaf lisse on some $U$ pure of weight $0$; one has
$$\sum_{x\in\Fq}K_\mcF(x)-{\tr(\Frob_q|H^2_c(U_{\ovFq},\mcF))}\ll C(\mcF)^2q^{1/2}.$$
More generally for any $n\geq 1$
$$\sum_{x\in\Fqn}K_{\mcF,n}(x)-{\tr(\Frob^n_q|H^2_c(U_{\ovFq},\mcF))}\ll C(\mcF)^2q^{n/2}.$$
In particular if $\mcF$ is geometrically irreducible or isotypic with no trivial component, one has
$$\sum_{x\in\Fq}K_\mcF(x)\ll C(\mcF)^2q^{1/2}.$$	
Here, the implied constants are absolute.
\end{corollary}
In practical applications we will be faced with situations where we have a sequence of sheaves $(\mcF_q)_q$ indexed by an infinite set of primes (with $\mcF_q$ a sheaf over  the field $\Fq$) such that the sequence of conductors $(C(\mcF_q))_q$ remains uniformly bounded (by $C$ say). In such situation, the above formula represents an asymptotic formula as $q\ra\infty$ for the sum of $q-O(1)$ terms
$$\sum_{x\in U(\Fq)}K_\mcF(x)$$
with main term ${\tr(\Frob_q|H^2_c(U_{\ovFq},\mcF))}$ (possibly $0$) and an error term of size $\ll C^2q^{1/2}$.

\section{Quasi-orthogonality relations}

We will often apply the trace formula and Deligne's theorem to the following sheaf: given $\mcF$ and $\mcG$ two $\ell$-adic sheaves both lisse on some non-empty open set $U\subset\Aa^1_{\Fq}$ and both pure of weight $0$; consider the tensor product $\mcF\otimes D(\mcG)$. This sheave is also lisse on $U$ and pure of weight $0$, moreover from the definition of the conductor (see \cite[Chap. 1]{GKM}) one sees that
\begin{equation}\label{conductortensor}
C(\mcF\otimes D(\mcG))\leq C(\mcF)C(\mcG).	
\end{equation}
The trace functions of $\mcF\otimes  D(\mcG)$ are given for $x\in U(\Fqn)$ by
$$x\mapsto K_{\mcF\otimes D(\mcG),n}(x)=K_{\mcF,n}(x)\ov{K_{\mcG,n}(x)}.$$
Therefore the trace formula can be used to evaluate the correlation sums between the trace function of $\mcF$ and $\mcG$,
$$\mcC(\mcF,\mcG):=\frac{1}{q}\sum_{x\in \Fq}K_{\mcF}(x)\ov{K_{\mcG}(x)};$$
more generally for any $n\geq 1$ we set
$$\mcC_n(\mcF,\mcG):=\frac{1}{q^n}\sum_{x\in \Fqn}K_{\mcF,n}(x)\ov{K_{\mcG,n}(x)}.$$
Indeed, by Corollary \ref{delignecor}, one has
\begin{equation}\label{corrn}
\mcC_n(\mcF,\mcG)=\tr(\Frob^n_q|V_{\mcF\otimes D(\mcG),\Ggeom})+O(\frac{C(\mcF)C(\mcG)}{q^{n/2}}).
\end{equation}
In particular if $C(\mcF)C(\mcG)$ are bounded while $q^n\ra\infty$, one obtains an asymptotic formula whose main term is given by the trace of the powers of Frobenius acting on the coinvariants of $\mcF\otimes D(\mcG)\simeq \Hom(\mcG,\mcF)$. 

\subsection{Decomposition of sheaves and trace functions}
Using first a weaker version of the formula (with an error term converging to $0$ as $n\ra\infty$), Deligne, on his way to the proof of Theorem \ref{delignethm}, established that any $\ell$-adic sheaf pure of weight $0$ is geometrically semi-simple (the representation $\rho_{\mcF|\Ggeom}$ decomposes into a direct sum of irreducible representations (of $\Ggeom$)) \cite[Th\'eor\`eme (3.4.1)]{WeilII}; the irreducible components occurring in the decomposition of $\rho_{\mcF|\Ggeom}$ are called the {\it geometric irreducible components of $\mcF$.}	 

This is not exactly valid for the arithmetic representation, but considering its semi-simplification, one obtains a decomposition
$$\rho^{ss}_{\mcF}=\bigoplus_{i\in I}\rho_{\mcF_i}$$
where the $\rho_{\mcF_i}$ are arithmetically irreducible (and pure) and lisse on $U$. Regarding geometric reducibility, each $\rho_{\mcF_i}$ is either geometrically isotypic or is induced from a representation of $\Gal(\Ksep/k.K)$ for $k$ some finite extension of $\Fq$. Since semi-simplification does not change the trace function, we obtain a decomposition of the trace function
$$K_\mcF=\sum_i K_{\mcF_i}.$$
Moreover a computation shows that whenever $\mcF_i$ is induced one has $K_{\mcF_i}\equiv 0$ on $U(\Fq)$. Therefore we obtain
\begin{proposition}\label{propdecomp} The trace function associated to some punctually pure sheaf $\mcF$ lisse on $U$ can be decomposed into the sum of $\leq C(\mcF)$  trace functions associated to sheaves $\mcF_i$, that are lisse on $U$, punctually pure of weight $0$, geometrically isotypic with conductors $C(\mcF_i)\leq C(\mcF)$.
\end{proposition}
This proposition   reduces the study of trace functions to trace functions associated to geometrically isotypic or (most of the time) geometrically irreducible sheaves. From now  on (unless stated otherwise) we will assume that the trace functions are associated to sheaves that are punctually pure of weight $0$ and geometrically isotypic. To ease notations, we say that such sheaves are "isotypic" or "irreducible" omitting the mention "geometrically" and likewise will speak of isotypic or irreducible trace functions. In such situation, using Schur lemma, the formula for \eqref{corrn} specializes to the
\begin{theorem}[Quasi-orthogonality relations]\label{thmcorrelation} Supppose that $\mcF$ and $\mcG$ are both geometrically isotypic with $n_\mcF$ copies of the irreducible component $\ov\mcF_{irr}$ for $\mcF$ and $n_\mcG$ copies of the irreducible component $\ov\mcG_{irr}$ for $\mcG$. There exists $n_\mcF.n_\mcG$ complex numbers $\alpha_{i,\mcF,\mcG}$ of modulus $1$ such that
\begin{equation}\label{quasiorth1}
\mcC_n(\mcF,\mcG)=(\sum_{i=1}^{n_\mcF n_\mcG}\alpha^n_{i,\mcF,\mcG})\delta_{\ov\mcF\sim_{geom}\mcG}+O(C(\mcF)^2C(\mcG)^2q^{-n/2}).
\end{equation} 
In particular if $\mcF$ and $\mcG$ are both geometrically irreducible there exist $\alpha_{\mcF,\mcG}\in\Ss^1$ such that
\begin{equation}\label{quasiorth2}
\mcC_n(\mcF,\mcG)=\alpha^n_{\mcF,\mcG}\delta_{\ov\mcF\sim_{geom}\mcG}+O(C(\mcF)^2C(\mcG)^2q^{-n/2}).
\end{equation}
In both \eqref{quasiorth1} and  \eqref{quasiorth2} the implicit constants are independent of $n$.
\end{theorem}
\begin{remark}
 Observe that for $\mcF$ and $\mcG$ either the Kummer or Artin-Schreier sheaves these equalities correspond to the orthogonality relations of characters.
\end{remark}
\begin{remark} If two geometrically irreducible sheaves $\mcF,\mcG$ are geometrically isomorphic, then their trace functions are proportional: more precisely one has for any $n$
$$K_{\mcF,n}=\alpha^n_{\mcF,\mcG}K_{\mcG,n}$$
where $\alpha_{\mcF,\mcG}$ is the complex number of modulus $1$ introduced in the previous statement.
	
\end{remark}

When $q^n$ is large compared to $C(\mcF)^2C(\mcG)^2$, the above formula gives a useful criterion to detect whether $\mcF$ and $\mcG$ have geometric irreducible components in common. While our focus is on the case $n=1$ and $q\ra\infty$ (while $C(\mcF)^2C(\mcG)^2$ remains bounded), the case $n\ra\infty$ will also prove useful. We start with the following easy lemma
\begin{lemma} Given $\alpha_1,\ldots,\alpha_d\in\Ss^1$, arbitrary  complex numbers of modulus $1$, one has
$$\limsup_{n\ra\infty}(\alpha_1^n+\ldots+\alpha_d^n)=d.$$
\end{lemma}
Using this lemma together with the decomposition into irreducible representations, one obtains the following
\begin{corollary}[Katz's Diophantine criterion for irreducibility] Let $\mcF$ be an $\ell$-adic sheaf lisse on $U$ pure of weight $0$ with decomposition into geometrically irreducible subsheaves denoted
$$\mcF^{geom}=\bigoplus_i\ov{\mcF}_{i}^{\oplus n_{i}}.$$
Then 
$$\limsup_{n\ra\infty}\mcC_n(\mcF,\mcF)=\sum_{\ov\mcF_{i}}n_{i}^2.$$
In particular,
$\mcF$ is geometrically irreducible if and only if
$$\limsup_{n\ra\infty}\mcC_n(\mcF,\mcF)=1.$$	
\end{corollary}

\subsection{Counting trace functions}
The above orthogonality relations lead to upper bounds for the number of geometric isomorphism classes of $\ell$-adic sheaves of bounded conductor (see \cite{MRL} for the proof):
\begin{theorem} Let $C\geq 1$, the number of geometric isomorphism classes of geometrically irreducible $\ell$-adic sheaves of conductor $\leq C$ is finite and bounded by
$$q^{O(C^6)}$$
where the implied constant is absolute.
\end{theorem}

\proof The principle of the proof is as follows: the sheaf-to-trace-function map $\mcF\to t_\mcF$ associates to the geometric isomorphism class of some sheaf a line in the $q$-dimensional Hermitian space $\Cc^{\Fq}$ of complex-valued functions on $\Fq$ with inner product
$$\peter{K,K'}=\frac{1}{q}\sum_{x\in\Fq}K(x)\ov{K'}(x).$$
The quasi-orthogonality relations show that these different lines are almost orthogonal to one another and so one obtains a number of almost orthogonal (circles of) unit vectors in the corresponding unit sphere. A sphere-packing argument for high-dimensional hermitian spaces (see \cite{kale}) implies that the number of such vectors cannot be too large.\qed

\section{Trace functions over short intervals}\label{Secshort}
In the next few sections, we discuss the correlations between trace functions and  other classical arithmetic functions. Indeed given a trace function
$$K_\mcF\colon \Aa^1(\Fq)=\Fq\ra \Cc$$
(extended from $U(\Fq)$ to $\Aa^1(\Fq)$ either by zero or by the middle-extension) we obtain a $q$-periodic function on $\Zz$ (which we also denote by $K_\mcF$) via the $\mods q
$-map 
$$K=K_\mcF\colon \Zz\to \Zz/q\Zz=\Aa^1(\Fq)\ra \Cc.$$ Given some other arithmetic function
$\lambda\colon \Nn\ra \Cc$ it is natural to compare them by evaluating their correlation sums
$$\sum_{n\leq N} K(n)\ov{\lambda(n)}$$
as $N\ra\infty$ (in suitable ranges of interest depending on $C(\mcF)$ and $\lambda$).

\subsection{The P\'olya-Vinogradov method}
 We start with the basic case where $\lambda=1_I$ is the characteristic function of an interval $I$ of $\Zz$ (which we may assume is contained in $[0,q-1]$). We want to evaluate non-trivially the sum
$$S(K;I):=\sum_{n\in I}K(n).$$
Remember that we mayand do  assume that $\mcF$ is geometrically isotypic and that if $I= [0,q-1]$ such sum can be dealt with by Deligne's theorem.

By Parseval, one has
$$S(K;I)=\sum_{y\in \Fq}\what K(y)\ov{\what{1_I}}(y)$$
where
\begin{equation}\label{eqfourier}
\what K(y)=\frac{1}{q^{1/2}}\sum_{x\in\Fq}K(x)\eq(xy)	
\end{equation}
 and
$$\what{1_I}(y)=\frac{1}{q^{1/2}}\sum_{x\in I}\eq(xy)$$
are the (normalized) Fourier transforms of $K$ and $1_I$ (for the abelian group $(\Fq,+)$).
One has
$$|\what{1_I}(y)|\ll \frac{1}{q^{1/2}}\min(|I|,\|\frac yq\|^{-1})
\ll \frac{1}{q^{1/2}}\min(|I|,\frac q{|y|})$$
(here $\|y/q\|$ denote the distance to the nearest integer) which implies that
$$\|\what{1_I}\|_1\ll \frac{|I|}{q^{1/2}}+q^{1/2}\log q.$$
Therefore one has
$$\sum_{n\in I}K(n)\ll \|\what K\|_\infty q^{1/2}\log q.$$
This leads us to look at the size of the Fourier transform $y\mapsto \what K(y)$. If $K$ is of the shape $\eq(ax)$ for some $a\in\Fq$, its Fourier transform is a Dirac function
$$\what K(y)=q^{1/2}\delta_{y=a\mods q}$$
and is therefore highly concentrated. To avoid this we make the following
\begin{definition} An isotypic sheaf $\mcF$ is Fourier if its geometric irreducible component is not (geometrically) isomorphic to any Artin-Schreier sheaf $\mcL_\psi$.
\end{definition}
In particular, if $K$ is Fourier of conductor $C(\mcF)$, it follows from Theorem \ref{thmcorrelation} that for any $y\in\Fq$
$$\what K(y)\ll C(\mcF)^2.$$
In that way we obtain the 
\begin{theorem}[P\'olya-Vinogradov bound] Let $\mcF$ be a Fourier sheaf of conductor $C(\mcF)$ and $K$ its associated trace function. For any interval $I$ of length $\leq q$, one has
	$$\sum_{x\in I}K(x)\ll C(\mcF)^2q^{1/2}\log q;$$
	here the implicit constant is absolute.
\end{theorem}
\begin{remark} This statement was obtained for the first time by P\'olya and Vinogradov, independently, in the case of Dirichlet characters $\chi$. In that case the Fourier transform is the normalized Gauss sum
$$\what\chi(y)=\eps_\chi(y)=\frac{1}{\qde}\sum_{x\in\Fq}\chi(x)\eq(xy)$$	
which is bounded in absolute value by $1$.
\end{remark}

Observe that this bound is better than the trivial bound
$$|\sum_{x\in I}K(x)|\leq C(\mcF)|I|$$
as long as $$|I|\gg_{C(\mcF)}q^{1/2}\log q.$$

This range is called the {\em P\'olya-Vinogradov range} and the question of bounding non-trivially for as many trace functions as possible over shorter intervals is a   fundamental problem in analytic number theory with many striking applications. At this moment, the problem is solved only in a very limited number of cases. One important example is the celebrated work of Burgess on Dirichlet characters \cite{Bur} which we  discuss in \S \ref{Bursec}. A lot of the forthcoming lectures will indeed be concerned with breaking this barrier in specific cases or in different contexts, and to give some  applications. 

\subsubsection{Bridging the P\'olya-Vinogradov range}

The following argument of Fouvry, Kowalski, Michel, Rivat, Soundararajan and Raju improves slightly the P\'olya-Vinogradov range:
\begin{theorem}\cite{FKMRRS} Let $\mcF$ be a Fourier sheaf of conductor $C(\mcF)$ and $K$ its associated trace function. For any interval $I$ of length $\sqrt q<|I|\leq q$, we have
$$\sum_{x\in I}K(x)\ll C(\mcF)^2\qde(1+\log(|I|/q^{1/2})).$$
\end{theorem}
\proof
Given $r\in\Zz$, let $I_r=r+I$; this is again an interval and  $S(K;I)$ and $S(K;I_r)$ differ only by $O(\|K\|_\infty r)$, which is a useful bound when $r$ is not too large.  Moreover 
$$\what{1_{I_r}}(y)=\eq(ry)\what{1_I}(y).$$
We have therefore
$$S(K;I)=\sum_{|y|\leq q/2}\what K(y)\ov{\what{1_I}(y)}\frac{1}{R}\sum_{0\leq r\leq R-1}\eq(-ry).$$
We choose $R=[\qde]+1$; using the bounds 
$$|\what{1_I}(y)|\ll \qmd\min(|I|,q/|y|),\ \sum_{0\leq r\leq R-1}\eq(-ry)\ll \min(R,q/|r|)$$
and 
$$\|K\|_\infty+\|\what K\|_\infty\ll C(\mcF)^2$$ we obtain the result.
\qed

\subsection{A smoothed version of the P\'olya-Vinogradov method}
Often in analytic number theory one is not faced with summing a trace function over an interval but instead against some smooth compactly supported function, for instance one has to evaluate sums of the shape
$$\sum_{n\in\Zz}K(n)V(\frac{n}N),\ V\in C_c^\infty(\Rr)\hbox{ fixed}.$$
By the Poisson summation formula  one has the identity
\begin{equation}\label{eqpoisson}
\sum_{n\in\Zz}K(n)V(\frac{n}N)=\frac{N}{q^{1/2}}\sum_{n\in\Zz}\what K(n)\what V(\frac{n N}q)	
\end{equation}
where
$$\what V(y)=\int_\Rr V(x)e(xy)dx$$
is the Fourier transform of $V(x)$ (over $\Rr$).	

Observe that $\what V(y)$ is not compactly supported but at least is of rapid decay:
$$\forall A\geq 0,\ \what V(y)\ll_{V,A}(1+|y|)^{-A}.$$ Therefore the dual sum in \eqref{eqpoisson} decays rapidly for $n\gg q/N$ and we obtain
\begin{proposition} We have
\begin{equation}\label{eqsmoothPV}
\sum_{n\in\Zz}K(n)V(\frac{n}N)\ll_{V} q^{1/2}\|\what K\|_\infty\ll_{V,\CF} q^{1/2}.	
\end{equation}
\end{proposition}

\subsection{The Deligne-Laumon Fourier transform}

The Fourier transform
$$K\mapsto \what K\colon  y\to \frac{1}{q^{1/2}}\sum_{x\in\Fq}K(x)\eq(-xy)$$
is a well-known and very useful operation on the space of function on $(\Zz/q\Zz,+)$. It serves to realize the spectral decomposition of the functions on $\Zz/q\Zz$ in terms of eigenvectors of the irreducible representations (characters) of $\Zz/q\Zz$. Let us recall that 
\begin{itemize}
\item The Fourier transform is essentially involutive: $$\what{\what K}(x)=K(-x);$$ stated otherwise, one has the Fourier inversion formula:
$$K(x)=\sum_{y\in\Fq}\what K(y)\eq(yx).$$
\item The Fourier transform is an isometry on $L^2(\Zz/q\Zz)$; stated otherwise, one has the Plancherel formula
$$\sum_{x\in\Fq}K(x)\ov{K'(x)}=\sum_{y\in\Fq}\what K(y)\ov{\what{K'}(y)}.$$
\item The Fourier transform behaves well with respect to  to additive and  multiplicative shifts: for $a\in\Fq,\ z\in\Fqt$,
$$\what{[+a]K}(y)=\eq(ay)\what K(y),\ \what{[\times z]K}(y)=[\times z^{-1}]\what K(y)=\what K(y/z).$$	
\end{itemize}

A remarkable fact, due to Deligne is that, to the Fourier transform for trace functions corresponds a "geometric Fourier transform" for sheaves. The following theorem is due to G. Laumon \cite{laumon87}:
 \begin{theorem} Let $\mcF$ be a Fourier sheaf, lisse on $U$ and pure of weight $0$. There exists a Fourier sheaf $\what \mcF$, lisse on some open set $\what U$, pure of weight $0$, such that if $K_{\mcF,n}$ denotes the (middle-extension of the) trace function of $\mcF$, the (middle extension of the) trace function of $\what\mcF$ is given by the Fourier transform $\what{K_{\mcF,n}}$ where
	  $$\what {K_{\mcF,n}}(x)=\frac{1}{q^{n/2}}\sum_{y}K_{\mcF,n}(y)\eq(\tr_{\Ff_{q^n}/\Fq}(xy)).$$
	  The map\footnote{This is in fact a functor in the derived category of constructible $\ell$-adic sheaves.} $\mcF\mapsto \what\mcF$ is called the geometric Fourier transform. The geometric Fourier transform satisfies (for $a\in\Fq,\ z\in\Fqt$)

	 $$\what{\what\mcF}=[\times -1]^*\mcF,\ \what{[+a]^*\mcF}=\mcL_{\eq(a).}\otimes \what\mcF,\ \what{[\times z]^*\mcF}=[\times z^{-1}]^* \what\mcF.$$
\end{theorem}

In addition, Laumon also defined local versions of the geometric Fourier transform making possible the computation of the local monodromy representations of $\what\mcF$ in terms of those of $\mcF$; using these results one deduces
\begin{proposition} Given $\mcF$ as above, one has
$$C(\what \mcF)\leq 10 C(\mcF)^2.$$
\end{proposition}
 
 Also the Fourier transform preserves irreducibility:
 
 \begin{proposition}
 The Fourier transform maps irreducible (resp.~isotypic) sheaves to irreducible (resp.~isotypic) sheaves.	
 \end{proposition}
\proof
 
 Given $\mcF$ a geometrically irreducible sheaf pure of weight $0$, to prove that $\what\mcF$ is irreducible, it is enough to show (by Katz's irreducibility criterion) that
 $$\limsup_n \mcC_n(\what\mcF,\what\mcF)=\limsup_n \frac{1}{q^n}\sum_{x\in\Fqn}|\what {K_{\mcF,n}}(x)|^2=1$$
 but by the Plancherel formula
 $$\frac{1}{q^n}\sum_{x\in\Fqn}|\what {K_{\mcF,n}}(x)|^2=
 \frac{1}{q^n}\sum_{y\in\Fqn}| {K_{\mcF,n}}(y)|^2$$
 and 
 $$\limsup_n \frac{1}{q^n}\sum_{y\in\Fqn}| {K_{\mcF,n}}(y)|^2=1$$
  by Katz's irreducibility criterion applied in the reverse direction.
 
\qed
 
 \begin{xca} Prove that the hyper-Kloosterman sheaves are geometrically irreducible ( hint: observe that the hyper-Kloosterman sums $\Kl_{k+1}$ can be expressed in terms of the Fourier transform of $\Kl_{k}$).	
 \end{xca}

\section{Autocorrelation of trace functions; the automorphism group of a sheaf}

The next couple of appplications we are going to discuss involve a special type of correlation sums between a trace function and its transform by an automorphism
 of the projective line. 

Let $\mcF$ be an $\ell$-adic sheaf lisse on $U\subset\Pp^1_\Fq$, pure of weight $0$, geometrically irreducible but non trivial, with conductor $C(\mcF)$. Let $\gamma$ be an automorphism of $\Pp^1_\Fq$: $\gamma$ is a fractional linear transformation:
$$\gamma\colon  z\to \gamma\cdot z=\frac{az+b}{cz+d},\ \begin{pmatrix}a&b\\c&d	
 \end{pmatrix}\in\PGL_2(\Fq).
$$
Let $\gamma^*\mcF$ be the associated pull-back sheaf; it is lisse on $\gamma^{-1}\cdot U$ and its trace function is
$$\gamma^*K(z)=K(\gamma\cdot z)=K(\frac{az+b}{cz+d}).$$ Moreover since $\gamma$ is an automorphism of $\Pp^1_\Fq$, one has $C(\gamma^*\mcF)=C(\mcF)$.

The correlations sums we will consider are those of $K$ and $\gamma^*K(z)$
$$\mcC(\mcF,\gamma):=\mcC(K,\gamma^*K)=\frac{1}q\sum_{z}K(z)\ov{K(\gamma\cdot z)}$$
and
$$\mcC_n(\mcF,\gamma):=\mcC_n(K,\gamma^*K)=\frac{1}{q^n}\sum_{z\in\Fqn}K_n(z)\ov{K_n(\gamma\cdot z)}$$
which are associated to the tensor product sheaf
$$\mcF\otimes \gamma^*D(\mcF)$$
which is lisse on $U_\gamma=U\cap \gamma^{-1}\cdot U.$

\subsection{The automorphism group}
The question of the size of the sums $\mcC(\mcF,\gamma)$ is largely determined by the following invariant of $\mcF$ (see \cite{FKM1,FKM2})
\begin{definition} Given $\mcF$ as above, the group of automorphisms of $\mcF$, denoted $\Aut_{\mcF}(\Fq)\subset\PGL_2(\Fq)$, is the group of $\gamma\in \PGL_2(\Fq)$ such that
$$\gamma^*\mcF\simeq_\geom\mcF.$$
The group $\Aut_{\mcF}(\Fq)$ is the group of $\Fq$-points of an algebraic subgroup, $\Aut_{\mcF}\hookrightarrow\PGL_{2}$ defined over $\Fq$. 
Let $B\subset \PGL_2$ the subgroup generated by upper-triangular matrices; we define
$$B_\mcF:=\Aut_{\mcF}\cap B$$ the subgroup of upper-triangular matrices of $\Aut_{\mcF}$ and $B_\mcF(\Fq)$	the group of $\Fq$-points.
\end{definition}
The relevance of this notion for the above correlations sums is the following
\begin{proposition} For $\gamma\not\in \Aut_{\mcF}(\Fq)$, one has
$$\mcC(K,\gamma)=O_{C(\mcF)}(q^{-1/2}).$$
\end{proposition}
In view of this proposition it is important to determine  $\Aut_{\mcF}(\Fq)$ and $B_\mcF(\Fq)$.

\begin{example}
Obviously any element of $\Aut_{\mcF}$ has to leave $\Pp^1(\ov\Fq)-U(\ov\Fq)$ invariant and all the points in the same orbit have isomorphic local monodromies. This may impose very strong constraints on $\Aut_{\mcF}$.

\begin{itemize}
\item If $\mcF$ is geometrically trivial then $\Aut_{\mcF}=\PGL_2$.
\item If $\psi\colon (\Fq,+)\ra\Ss^1$ is non trivial then $G_{\mcL_\psi}=N=\{\begin{pmatrix}1&x\\&1	
\end{pmatrix}\subset\PGL_2
 \}$.
 \item If $\chi\colon (\Fq,+)\ra\Ss^1$ is non trivial, then $$G_{\mcL_\chi}=T^{0,\infty}=\{\begin{pmatrix}a&0\\0&d	
\end{pmatrix}\subset\PGL_2
 \}$$
 is the diagonal torus,  unless $\chi$ is quadratic in which case $G_{\mcL_\chi}=N(T^{0,\infty})$ is the normalizer of the diagonal torus. 
\item For the Kloosterman sheaves, one can show that $\mcG_{\KL_k}$ is trivial: since $\KL_k$ is not lisse at $0$ and $\infty$, with Swan conductor $0$ at $0$ and $1$ at $\infty$, one has $\mcG_{\KL_k}\subset T^{0,\infty}$. One can then show (see \cite{MiDMJ}) that $[\times a]^*{\KL_k}\simeq_\geom{\KL_k}$ iff $a=1$.
\end{itemize}

\end{example}

Given $x\not=y\in\Pp^1(\ov\Fq)$, we denote by $T^{x,y}$ the pointwise stabilizer  of the pair $(x,y)$ (this is a maximal torus defined over some finite extension of $\Fq$) and $N(T^{x,y})$ its normalizer. The torus $T^{x,y}$ is defined over $\Fq$ if $x,y$ belong to $\Pp^1(\Fq)$ or if $x,y$ belong to $\Pp^1(\Ff_{q^2})$ and are Galois conjugates.

\begin{proposition}\label{thmautogroup} Suppose $q\geq 7$. Given $\mcF$ as above, at least one of the following holds:
\begin{itemize}
\item $\CF>q$.
\item $q$ does not divide $|\Aut_{\mcF}(\Fq)|$ and either $\Aut_{\mcF}(\Fq)$ is of order $\leq 60$ or is a subgroup of the normalizer of some maximal torus $N(T^{x,y})$ defined over $\Fq$.
\item $q$ divides $|\Aut_{\mcF}(\Fq)|$ and then $\mcF\simeq \sigma^*\mcL_\psi$ for some $\psi$ and $K(x)=\alpha\psi(\sigma.x)$ for 
for some $\sigma\in\PGL_2(\Fq)$ and $\Aut_{\mcF}(\Fq)=\sigma N\sigma^{-1}$.
\end{itemize}
\end{proposition}
\begin{remark}
Observe that in the last case $$\mcC(K,\gamma)=|K(0)|^2\mcC(\psi(\sigma.x),\gamma)$$
\end{remark}

Concerning the size of the group $B_\mcF(\Fq)$, one can show that
\begin{theorem}\label{thmautogroupB} Let $\mcF$ be an isotypic sheaf whose geometric components are not  isomorphic to $[+x]^*\mcL_\chi$ for some $x\in\Fq$ and some multiplicative character $\chi$ and such that 
$$C(\mcF)< q.$$
Then
$$|B_\mcF(\Fq)|\leq C(\mcF).$$
\end{theorem}
The proof of this theorem involves the following rigidity statements \cite[Lemma 2.6.13]{KatzRLS}:
\begin{proposition}
Let $\mcL $ be geometrically irreducible.
\begin{itemize}
\item If for some $x\in\Fqt$, $[+x]^*\mcL\simeq\mcL$, then either
$$C(\mcL)>q\hbox{ or }\mcL\simeq \mcL_\psi \hbox{ for some } \psi.$$
\item If $\Aut_\mcL(\Fq)$ contains a subgroup of	 order $m$ of diagonal matrices then either
 $$c(\mcL)>m\hbox{ or }\mcL\simeq\mcL_\chi \hbox{ for some } \chi.$$
\end{itemize}
	
\end{proposition}

\section{Trace functions vs.~primes}\label{Secprimes}

Another possible question to consider (natural from the viewpoint of analytic number theory at least) is how trace functions correlate with the characteristic function of the primes. In this section, we discuss the structure of the proof of the following result: 

\begin{theorem}[Trace function vs.~primes, \cite{FKM2}]\label{thmprimesumthm}
  Let $\mcF$ be a geometrically isotypic sheaf of conductor $C(\mcF)$ whose geometric components are not of the shape $\mcL_\psi\otimes\mcL_\chi$ and let $K$ its associated trace function. For any $V\in C^\infty_c(\Rr_{>0})$, one has \par
\begin{align}
\label{primesuminterval}
\sum_\stacksum{p\ \text{prime}}{p \leq X}K(p)&\ll
X(1+q/X)^{1/12}p^{-\eta/2},\\
\label{primesumsmooth}\sum_{p\ \text{prime}}K(p)V\Bigl(\frac{p}X\Bigr)&\ll
X(1+q/X)^{1/6}q^{-\eta},
\end{align}
for $X\ll q$ and $\eta<1/24$. The implicit constants depend only on $\eta$,
$\CF$ and $V$. Moreover, the dependency on $\CF$ is at most polynomial.
\end{theorem}

\begin{remark} This result exhibits cancellations when summing trace functions  along the primes in intervals of length larger than $q^{3/4}$. It is really a pity that Dirichlet characters are excluded by our hypotheses: such a bound in that case would amount to a quasi generalized Riemann hypothesis for the corresponding Dirichlet character $L$-function !
\end{remark}
We discuss the proof for $X=q$.
 
\subsection{Combinatorial decomposition of the characteristic function of the primes}
As is well-known, the problem is equivalent to bounding the sum
$$\sum_{n}\Lambda(n) K(n)V\Bigl(\frac{n}q\Bigr)$$
where $$\Lambda(n)=\begin{cases}\log p&\hbox{ if }n=p^\alpha\ \alpha\geq 1\\
0&\hbox{ otherwise, } 	
\end{cases}$$
 is the von Mangoldt function. A standard method in analytic number theory is a combinatorial decomposition of this function as a sum of Dirichlet convolutions; one way to achieve this is to use the celebrated Heath-Brown identity:
\begin{lemma}[Heath-Brown]  \label{lemHB}
  For any integer  $J\geq 1$ and $n< 2X$, we have
$$
\Lambda(n)=-\sum_{j=1}^J(-1)^j\binom{J}{j} \sum_{m_1,\ldots, m_j\leq
  Z}\mu(m_1)\ldots\mu(m_j) \sum_{m_1\ldots m_jn_{1}\ldots
  n_{j}=n}\log n_1,
$$
where $Z=X^{1/J}$.
\end{lemma}
Hence splitting the range of summation of the various variables appearing (using partition of unity) and separating these variables, our preferred sum decomposes (essentially) into $O((\log X)^{2J})$ sums of the shape
$$\Sigma(M_1,\ldots,M_{2j})=\sumsum_{m_1,\ldots m_{2j}}\mu(m_1)\ldots\mu(m_j) K(m_1.\ldots.m_{2j})V_{1}	\Bigl(\frac{m_{1}}{M_{1}}\Bigr)\ldots V_{2j}	\Bigl(\frac{m_{2j}}{M_{2j}}\Bigr)$$
for $j\leq J$; here $V_i,\ i=1,\ldots 2j$ are smooth  functions compactly supported in $]1,2[$, and $(M_1,\ldots,M_{2j})$ is a tuple satisfying
$$M_i=:q^{\mu_i},\ \forall i\leq j,\ \mu_i\leq 1/J,\ \sum_{i\leq 2j}{\mu_i}=1+o(1);$$
The objective is to show that
$$\Sigma(M_1,\ldots,M_{2j})\ll q^{1-\eta}$$
for some fixed $\eta>0$. We will take $J=3$ so that $Z=q^{1/3}$. We may assume that
$$\mu_1\leq\ldots\leq \mu_j\leq 1/3,\ \mu_{j+1}\leq \ldots\leq\mu_{2j}.$$
We will bound these sums differently depending on the vector $(\mu_1,\ldots,\mu_{2j})$.

Let $0<\delta<1/6$ be some small but fixed parameter to be chosen optimally later.

\begin{enumerate}

\item Suppose that  $\mu_{2j}\geq 1/2+\delta$. Then $m_{2j}$ is a long "smooth variable" (because the weight attached to it is smooth); therefore using  \eqref{eqsmoothPV} to sum over $m_{2j}$ while fixing the other variables, we get
$$\Sigma(M_1,\ldots,M_{2j})\ll q^{\mu_1+\ldots\mu_{2j-1}}q^{1/2+o(1)}=q^{1-\delta+o(1)}.$$
(In the literature, sum of that shape are called "type I" sums).

\item  We may therefore assume that $$m_{j+1}\leq \ldots\leq \mu_{2j}\leq 1/2+\delta;$$
in other words, there is no  long smooth variable. 
What one can then do is to group variables together to form longer ones: for this one partitions the indexing set  into two blocks
$$\{1,\ldots,2j\}=\mcI\sqcup \mcI',$$ and form the variables 
$$m=\prod_{i\in \mcI} m_i,\ n=\prod_{i'\in\mcI'}m_{i'}$$
so that denoting by $\alpha_m$  the Dirichlet convolutions
of either $\mu(\cdot)V(\frac{\cdot}{M_i})$ or $V(\frac{\cdot}{M_i})$ for $i\in\mcI$ and similarly for $\beta_n$ for $i'\in\mcI'$, we are led to bound  bilinear sums of the shape
\begin{equation}\label{eqbilineardef} B(K;\alpha,\beta)=\sumsum_{m\ll M,n\ll N}\alpha_m\beta_n K(mn).	
\end{equation}
where
$$M=q^\mu,\ \mu=\sum_{i\in \mcI}\mu_i,\ N=q^\nu,\ \nu=\sum_{i'\in \mcI'}\mu_{i'}.$$
The weights $\alpha_m,\beta_n$ are rather irregular and it is difficult to exploit their structure (such sums are called "type II"). 

Assuming that the irreducible component of $\mcF$ is not of the shape $\mcL_\chi\otimes\mcL_\psi$, we will prove in Theorem \ref{thmbilinear} below the following bound
	$$\Sigma(M_1,\ldots,M_{2j})=B(K;\alpha,\beta)\ll_{C(\mcF)}\|\alpha_M\|_2\|\beta_N\|_2(MN)^{1/2}(\frac{1}M+\frac{q^{1/2}\log q}{N})^{1/2}.$$

Assuming that $$\mu\geq {\delta}\hbox{ and }\nu\geq {1/2+\delta}$$ we obtain that
$$B(K;\alpha,\beta)\ll q^{1-\delta/2+o(1)}.$$

\item  It remains to treat the sums for which neither $\mu_{2j}\leq 1/2+\delta$ nor a decomposition as in (2) exist. This necessarily implies that $\sum_{i\leq j}\mu_i\leq 1/3,$ $j\geq 2$ and $\mu_{2j-1}+\mu_{2j}\geq 1-\delta.$
Setting $M=M_{2j-1}$ and $N=M_{2j}$  , denoting
$$a=m_1\ldots m_{2j-2}\ll q^{\delta},$$ it will be sufficient to obtain a bound of the shape
$$\sum_{m,n\geq 1}K(amn)V(\frac{m}M)W(\frac{n}N)\ll_{V,W} (MN)^{1-\eta}$$
for some $\eta>0$ whenever $MN$ is sufficiently close to $q$. 
What we have are is a sum involving two smooth variables which are however too short for the P\'olya-Vinogradov method to work, but whose product is rather long. We call these sums "type I$1/2$". We will then use Theorem \ref{cortypeI1/2} below whose proof is discussed in \S \ref{sectracemodular}. Observe that this theorem provides a bound which is non trivial as long as $MN\geq q^{3/4}$.

\item Optimizing parameters in these three approaches leads to Theorem \ref{thmprimesumthm}.
\end{enumerate}

\begin{theorem}\label{cortypeI1/2} Let $\mcF$ be a geometrically  isotypic Fourier sheaf of conductor $C(\mcF)$ and $K$ its associated trace function. For any $V,W\in C^\infty_c(\Rr_{>0})$, any $M,N\geq 1$ and any $\eta<1/8$, one has
$$\sum_{m,n\geq 1}K(mn)V(\frac{m}M)W(\frac{n}N)\ll_{V,W,C(\mcF)} MN(1+\frac{q}{MN})^{1/2}q^{-\eta/2}.$$
\end{theorem}

\section{Bilinear sums of trace functions}\label{secbilinear}
Let $K$ be a trace function associated to some isotypic sheaf $\mcF$, pure of weight $0$ and let $(\alpha_m)_{m\leq M}$, $(\beta_n)_{n\leq N}$ be arbitrary complex numbers. In this section, we  bound the "type II" bilinear sums  encountered in the previous section : 
$$B(K;\alpha,\beta)=\sumsum_{m\leq M,n\leq N}\alpha_m\beta_n K(mn).$$

Using the Cauchy-Schwarz inequality, the trivial bound is
$$|B(K;\alpha,\beta)|\ll_{C(\mcF)}\|\alpha_M\|_2\|\beta_N\|_2(MN)^{1/2}.$$
We wish to improve over this bound.
\begin{theorem}[Bilinear sums of trace functions]\label{thmbilinear} Notations as above; assume that $1\leq M,N<q$ and that the irreducible component of $\mcF$ is not of the shape $\mcL_\chi\otimes\mcL_\psi$. Then
	$$B(K;\alpha,\beta)\ll_{C(\mcF)}\|\alpha_M\|_2\|\beta_N\|_2(MN)^{1/2}(\frac{1}M+\frac{q^{1/2}\log q}{N})^{1/2}.$$
\end{theorem}
\begin{remark}This bound is non-trivial as soon as $M\gg 1$ and $N\gg \qde\log q$.
\end{remark}
\proof
By Cauchy-Schwarz, we have
\begin{equation}\label{eqfirstCS}
|B(K;\alpha,\beta)|^2\leq  \|\beta_N\|_2^2\sum_{m_1,m_2\leq M}\alpha_{m_1}\ov{\alpha_{m_2}}\sum_{n\leq N}K(m_1n)\ov K(m_2n).
\end{equation}
We do not expect to gain anything from the diagonal terms $m_1\equiv m_2\mods q$ (equivalently, $m_1=m_2$ since $M<q$) and the contribution of such terms is bounded trivially by
\begin{equation}\label{eqdiagbound}
\ll_{C(\mcF)}\|\alpha_M\|_2^2\|\beta_N\|^2_2N.	
\end{equation}
As for the non-diagonal terms, their contribution is
$$ \|\beta_N\|_2^2\sum_{m_1\not=m_2\mods q}\alpha_{m_1}\ov{\alpha_{m_2}}\sum_{n\leq N}K(m_1n)\ov K(m_2n).$$
Using the P\'olya-Vinogradov method, we are led to evaluate the Fourier transform of
$$n\mapsto K(m_1n)\ov K(m_2n).$$
By the Plancherel formula, this Fourier transform equals
\begin{eqnarray*}
y\mapsto\frac{1}{\qde}\sum_{x\in\Fq}K(m_1x)\ov K(m_2x)\eq(-yx)&=&
\frac{1}{\qde}\sum_{z\in\Fq}\what K((z-y)/m_1)\ov {\what K}(z/m_2)
\\&=&\frac{1}{\qde}\sum_{z\in\Fq}\what K((m_2z-y)/m_1)\ov {\what K}(z)\\
&=&	\frac{1}{\qde}\sum_{z\in\Fq}\what K(\gamma z)\ov {\what K}(z)
\end{eqnarray*}
with
$$\gamma=\begin{pmatrix}
m_2/m_1&-y/m_1\\0&1	
\end{pmatrix}\in B(\Fq).
$$
This sum is $q^{1/2}$ times $\mcC(\what\mcF,\gamma)$, the correlation sum associated to the isotypic sheaves $\what\mcF$ and $\gamma^*\what\mcF$, whose conductors are controlled in terms of $C(\mcF)$.

If $\gamma\not\in B_\mcF(\Fq)$ we have
\begin{equation}\label{eqcorB}
\mcC(\what\mcF,\gamma)\ll_{C(\mcF)}\frac{1}{\qde}.	
\end{equation}
The condition that the irreducible component of $\mcF$ is not of the shape $\mcL_\chi\otimes\mcL_\psi$ translates into the irreducible component of $\what\mcF$ not being of the shape $[+x]^*\mcL_{\ov\chi}$. In that case, by Theorem \ref{thmautogroupB}, there is a set $S_\mcF\subset \Fqt$ such that for any $(m_1,m_2,y)\in\Fqt\times\Fqt\times\Fq$ for which $m_2/m_1\not\in S_\mcF$ one has
$$\mcC(\what\mcF,\gamma)\ll_{\CF} q^{-1/2}.$$

Returning to \eqref{eqfirstCS}, we bound trivially (by \eqref{eqdiagbound}) the contribution of the $O_\mcF(M)$ $(m_1,m_2)$ such that the ratio $m_2/m_1\mods q$ is in $S_\mcF$. For the other terms, we may use the P\'olya-Vinogradov method  and bound these terms by
$$\ll_{C(\mcF)}\|\alpha_M\|_2^2\|\beta_N\|_2^2Mq^{1/2}\log q.$$
Combining these bounds leads to the final result.
\qed

\section{Trace functions vs.~modular forms}\label{sectracemodular}
In this section we discuss the proof of Theorem \ref{cortypeI1/2}. This theorem is a special case of the resolution of another problem: the question of the correlation between trace functions and the Fourier coefficients $(\rf(n))_{n}$ of some modular Hecke eigenform (cf.~\cite[Chap. 14\&15]{IwKo} and references herein for a quick introduction to the theory modular forms). Given some trace function, we  consider the correlation sum
$$\tsum(K,f;X):=\sum_{n\leq X}\rf(n)K(n)$$ or its smoothed version
$$\tsum_V(K,f;X):=\sum_{n}\rf(n)K(n)V(\frac{n}X).$$
These sums are bounded (using the Rankin-Selberg method) by
$$O_{C(\mcF),f}(X\log^{3} X).$$
It turns out that the problem of bounding $\tsum(K,f;X)$ and $\tsum_V(K,f;X)$ non-trivially is most interesting when $N$ is of size $q$ or smaller. 

In this section, we sketch the proof of the following
\begin{theorem}[Trace function vs.~modular forms, \cite{FKM1}]\label{thmKmodular} Let $\mcF$ be an irreducible Fourier sheaf of weight $0$ and $K$ its associated trace function. Let $(\rf(n))_{n\geq 1}$ be the sequence of Fourier coefficients of some modular form $f$ with trivial nebentypus and $V\in C^\infty_c(\Rr_{>0})$. For $X\geq 1$ and any $\eta<1/8$, we have
$$\tsum(K,f;X)\ll X(1+\frac{q}X)^{1/2}q^{-\eta/2},$$
and
$$\tsum_V(K,f;X)\ll X(1+\frac{q}X)^{1/2}q^{-\eta}.$$
The implicit constants depend only on $\eta$, $f$,
$\CF$ and $V$. Moreover, the dependency on $\CF$ is
at most polynomial.
\end{theorem}

 This result shows the absence of correlation when $X\gg q^{1-1/8}$.
 The proof, which uses the amplification method and the Petersson-Kuznetzov trace formula, will ultimately be a consequence of Theorem \ref{thmautogroup}.

We give below an idea of the proof. To simplify matters, we will assume that $X=q$ and we wish to bound non-trivially the sum
\begin{equation}\label{SVdefnoq}
\tsum_V(K,f):=\sum_{n\geq 1}\rf(n)K(n)V(\frac{n}q)	
\end{equation}
for $V$ a fixed smooth function. Moreover, to simplify things further, we will assume that $f$ has level $1$ and is cuspidal and holomorphic of very large (but fixed) weight.

\subsection{Trace functions vs.~the divisor function}
An important special case of Theorem \ref{thmKmodular} is when $f$ is an Eisenstein series, for instance when
$$f(z)=\frac{\partial}{\partial s}E(z,s)_{|s=1/2}\hbox{ for }E(z,s)=\frac{1}2\sum_{(c,d)=1}\frac{y^s}{|cz+d|^{2s}}$$
is the non-holomorphic Eisenstein series at the central point. In that case we have 
$$\rf(n)=d(n)$$ the divisor function, and so one has
\begin{equation}\label{eqbounddiv}
\sum_{m,n\geq 1}K(mn)V(\frac{mn}X)\ll_{V,\CF} X(1+\frac{q}X)^{1/2}q^{-\eta}	
\end{equation}
whenever $K$ is the trace function of a Fourier sheaf. This bound holds similarly for the unitary Eisenstein series $E(z,s)$ at any $s=\frac12+it$, where the divisor function is replaced by
$$d_{it}(n)=\sum_{ab= n}(a/b)^{it}.$$
Such general bounds make it possible to separate the variables $m,n$ in \eqref{eqbounddiv} and eventually to prove Theorem \ref{cortypeI1/2}.

\begin{remark} As we will see below, the proof of Theorem \ref{thmKmodular} is not a "modular form by modular form" analysis; instead the proof is global, involving the full automorphic spectrum, and establishes the required bound "for all modular forms $f$ at once", including Eisenstein series and therefore proving Theorem \ref{cortypeI1/2} on the way.
\end{remark}

\subsection{Functional equations}
Our first objective is to understand why the range $X=q$ is interesting. This come from the functional equations satisfied by modular forms as a consequence of their automorphic properties. These equations present themselves in various shapes. One is the Voronoi summation formula, which in its
simplest form is the following:
\begin{proposition}[Voronoi summation formula]\label{Voronoi} Let $f$ be a holomorphic modular form of weight $k$ and level $1$ with Fourier coefficients $(\rf(n))_n$. Let $V$ be a smooth compactly supported function, $q\geq 1$ and $(a,q)=1$. We have
  for $X>0$
$$
\sum_{n\geq 1} \rf(n)V\Bigl(\frac{n}X\Bigr)e\Bigl(\frac{an}{q}\Bigr) = 
\eps(f)\frac{X}{q} \sum_{n\geq
  1}\rf(n)e\Bigl(-\frac{\overline{a}n}{q}\Bigr)
  \widetilde V\Bigl(\frac{Xn}{q^2}\Bigr) 
$$
where $\eps(f)=\pm 1$ denotes the sign of the functional equation of $L(f,s)$, and
$$\widetilde V(y)=\int_{0}^\infty V(u)\mathcal{J}_k(4\pi\sqrt{ uy})du,$$
with
\begin{equation*}
  \mathcal{J}_k(u)  =
  2\pi i^kJ_{k-1}(u),
  \end{equation*}
  where $$J_{k-1}(x)=\sum_{l=0}^\infty\frac{(-1)^l}{l !(l+k-1) !}(\frac{x}{2})^{2l+k-1}$$
is the Bessel function of order $k-1$.
\end{proposition}
There are several possible proofs of this proposition: one can proceed classically from the Fourier expansion of the modular form $f$ using automorphy relations (see \cite[Theorem A.4]{KMVDMJ}). Another more conceptual approach is to use the Whittaker model of the underlying automorphic representation; this approach extends naturally to higher rank automorphic forms (see \cite{IT}). One could also point out other related works like \cite{MilSch} as well as the recent paper \cite{KirZhou}. 
We can extend this formula to general functions modulo $q$. Given $K\colon \Zz\ra\Cc$ a $q$-periodic function, we define its \emph{Voronoi transform} $\bessel{K}$ of $K$ as
$$
\bessel{K}(n) = \frac{1}{\sqrt{q}}\sum_{\substack{h\bmod q\\(h,q) =1}}
\fourier{K}(h) \eq ({\overline h n} )=\frac{1}{\sqrt{q}}\sum_{\substack{h\bmod q\\(h,q) =1}}
\fourier{K}(h^{-1}) \eq ({h n} ).
$$ 
Combining the above formula with the Fourier decomposition
$$K(n)=\frac{1}{q^{1/2}}\sum_{a\mods q}\what K(a)\eq(-an),$$
we get
\begin{corollary}\label{corvoronoi} Notations are above, given $K$ a $q$-periodic arithmetic function, we have  for $X>0$ 
\begin{eqnarray*}
\sum_{n\geq 1}\rf(n)K(n)V\Bigl(\frac nX\Bigr)&=& \frac{\what K(0)}{q^{1/2}}\sum_{n\geq 1} \rf(n)V\Bigl(\frac nX\Bigr)+\\ 
&&\ \eps(f)\frac{X}{q} \sum_{n\geq
  1}\rf(n)\widecheck K(- n)\widetilde V\Bigl(\frac{nX}{q^2}\Bigr).	
\end{eqnarray*}

\end{corollary}
\begin{remark} Another way to obtain such result is to consider the Mellin transform of (the restriction to $\Fqt$ of) $K$:
$$\tilde K(\chi)=\frac{1}{(q-1)^{1/2}}\sum_{x\in\Fqt}K(x)\chi(x)$$
so that for $x\in\Fqt$
$$K(x)=\frac{1}{(q-1)^{1/2}}\sum_{\chi}\tilde K(\chi)\chi^{-1}(x).$$
One can then use the (archimedean) inverse-Mellin transform and the functional equation satisfied by the Hecke $L$-function
$$L(f\otimes\chi,s)=\sum_{n\geq 1}\frac{\rf(n)\chi(n)}{n^{s}}$$
to obtain the formula. For this, one observes that the Mellin transform of $\widecheck{K}_{|\Fqt}$ is proportional to
$$\chi\mapsto \eps(\chi)\tilde K(\chi^{-1})$$
where $\eps(\chi)$ is the normalized Gauss sum. This method extends easily to automorphic forms of higher rank but uses the fact that $q$ is prime (so that $\Fqt$ is not much smaller that $\Fq$).
\end{remark}
The identity of Corollary \ref{corvoronoi} is formal and has nothing to do whether $K$ is a trace function or not. In particular applying it to the Dirac function $\delta_a(n)=\delta_{n\equiv a\mods q}$, for some $a\in\Fqt$ we obtain
$$\widehat{\delta_a}(h)=\frac{1}\qde \eq(ah),\ \widecheck \delta_a(n)=\frac{1}{\qde}\Kl_2(an;q)$$
so that
\begin{eqnarray}\label{eqdeltacase}
q^{1/2}\sum_{n\equiv a\mods q}\rf(n)V\Bigl(\frac nX\Bigr)&=& \frac{1}{q^{1/2}}\sum_{n\geq 1} \rf(n)V\Bigl(\frac nX\Bigr)+\\ 
&&\ \eps(f)\frac{X}{q} \sum_{n\geq
  1}\rf(n)\Kl_2(-an;q)\widetilde V\Bigl(\frac{nX}{q^2}\Bigr). 	\nonumber
\end{eqnarray}
This is an example of a natural transformation which, starting from the elementary function $\delta_a$ produces a genuine trace function ($\Kl_2$).

Besides this case we would like to use the formula for $K$ a trace function. We observe that the Voronoi transform $\widecheck K$ is "essentially" the Fourier transform of the function $$h\in\Fqt\mapsto \what K(h^{-1})=\what K(w\cdot h)$$ with $w=\begin{pmatrix}0&1\\1&0	
\end{pmatrix}$; it is therefore essentially involutive. It would be useful to know that $\widecheck K$ is a trace function. Suppose that $K$ is associated to some isotypic Fourier sheaf $\mcF$, then $\widecheck K$ is a (isotypic) trace function as long as  $w^*\what\mcF$ is a Fourier sheaf. This means that $\what\mcF$ has no irreducible constituent of the shape $w^*\mcL_\psi$ which (by involutivity of the Fourier transform means that $\mcF$ has no irreducible constituent isomorphic to some Kloosterman sheaf $\KL_2$. This reasoning\footnote{by involutivity of the Voronoi transform} is essentially the reverse of the one leading to \refs{eqdeltacase}.

Let us assume that $\widecheck K$ is also a trace function. Then, integration by parts show that for $V$ smooth and compactly supported, $\widetilde V(x)$ has rapid decay for $x\gg 1$.  Hence Corollary \ref{corvoronoi} is an equality between a sum of length $X$ and a sum of length about $q^2/X$ (up to the term $\frac{\what K(0)}{q^{1/2}}\sum_{n\geq 1} \rf(n)V\Bigl(\frac nX\Bigr)$ which is easy to understand). The two lengths are the same when $X=q$.

\subsection{The amplification method}

As mentioned above Theorem \ref{thmKmodular} is proven "for all modular forms at one" as a consequence of the amplification method. 

The principle of the amplification method (invented by H. Iwaniec and which  in the special case $K=\chi$ was used first by Bykovskii) consist, in  the following. For $L\geq 1$ and $(x_l)_{l\leq L}$ real numbers we consider the following average over orthogonal bases of modular forms (holomorphic or general) of level $q$:
\begin{equation}\label{momentholk}
M_k(K ):=\sum_{g\in\mcB_k(q)}|A(g )|^2|\tsum_V(g,K)|^2	
\end{equation}
(cf.~\eqref{SVdefnoq} for the definition of $\tsum_V(g,K)$) and
\begin{multline}\label{momentdef} 
M(K ):=\sum_{k \equiv 0 \mods{2},\ k>0} \dot{\phi}(k)(k-1)
\sum_{g\in\mcB_k(q)}|A(g )|^2|\tsum_V(g,K)|^2
  \\
  + \sum_{g\in\mcB(q)} \tilde{\phi}(t_g)\frac{4 \pi }{\cosh(\pi
    t_g)}|A(g )|^2|\tsum_V(g,K)|^2\\
  + \,\sum_{g\in
    \mcB_E(q)}\int_{-\infty}^{\infty}\tilde{\phi}(t)\frac{1}{\cosh(\pi
    t)} |A(g,t)|^2|\tsum_V(E_{g}(t),K)|^2\,dt,
\end{multline}
where $\mcB_{k}(q),\ \mcB(q),\ \mcB_E(q)$ denote orthonormal bases of Hecke eigenforms of level $q$ (either holomorphic of weight $k$ or Maass or Eisenstein series), $\dot{\phi},\ \tilde{\phi}$ are  weights constructed from some smooth function, $\phi$, rapidly decreasing at $0$ and $\infty$,  which depend only on the spectral parameters of the forms and for each form $g$, $A(g )$ ("A" is for amplifier) is the linear form in the Hecke eigenvalues $(\lamg(n))_{(n,q)=1}$ given by
$$A(g )=\sum_{l\leq L}x_l\lamg(l).$$
 
The weights $\tilde{\phi}$ are positive while the weight $\dot{\phi}(k)$ is positive at least for $k$ large enough; one can then add to this quantity a finite linear combination of the $M_k(K),\ k\ll 1$ from which one can bound  
\begin{multline}\label{momentdef2} 
|M|(K ):=\sum_{k \equiv 0 \mods{2},\ k>0} |\dot{\phi}(k)|(k-1)
\sum_{g\in\mcB_k(q)}|A(g )|^2|\tsum_V(g,K)|^2
  \\
  + \sum_{g\in\mcB(q)} \tilde{\phi}(t_g)\frac{4 \pi }{\cosh(\pi
    t_g)}|A(g )|^2|\tsum_V(g,K)|^2\\
  + \,\sumsum_{g\in
    \mcB_E(q)}\int_{-\infty}^{\infty}\tilde{\phi}(t)\frac{1}{\cosh(\pi
    t)} |A(g,t)|^2|\tsum_V(E_{g}(t),K)|^2\,dt.
\end{multline}
As we explain below one will be able to prove the following bound
\begin{equation}\label{Mbound}
M(K ), M_{k}(K )\ll_{\CF} q^{o(1)}(q\sum_{l\leq L}|x_l|^2+q^{1/2}L(\sum_{l\leq L}|x_l|)^2).	
\end{equation}
Now if $f$ is a Hecke-eigenform of level $1$ (of $L^2$ norm $1$ for the usual inner product on the level one modular curve) then $f/(q+1)^{1/2}$ embeds in an orthonormal basis of forms of level $q$. 

Since all the terms in $|M|(K )$ are non-negative, this sums bounds any of its terms occurring discretely  (i.e.~when $f$ is a cusp form). Therefore we obtain
$$\frac{1}{q+1}|A(f )|^2|\tsum_V(f,K)|^2\ll_{\CF,f} q^{o(1)}(q\sum_{l\leq L}|x_l|^2+q^{1/2}L(\sum_{l\leq L}|x_l|)^2).$$
Now we perform amplification by choosing some  bounded sequence $(x_l)_{l\leq L}$ tailor made for $f$ such that $A(f)$ is "large".
Specifically, choosing
$$x_l=\mathrm{sign}(\lf(l)),$$
we obtain 
$$|A(f)|\gg L^{1+o(1)}.$$
Dividing by $L$ we obtain
$$|\tsum_V(f,K)|^2\ll q^{o(1)}(q^2/L+q^{3/2}L^2)$$
and the optimal choice is $L=q^{1/6}$ giving us
$$\tsum_V(f,K)\ll q^{1-1/12+o(1)}.$$

\subsection{Computing the moments}
We now bound $M(K)$. Opening squares and using the multiplicative properties of Hecke eigenvalues, we are essentially reduced to bounding sums of the shape
\begin{equation}\label{eqmodelsummodular}
\sumsum_{m,n}V(\frac{m}q)V(\frac{n}q)K(m)\ov{K(n)}\Delta_{q,\phi}(lm,n)	
\end{equation}
and
\begin{equation}\label{eqmodelsummodular}
\sumsum_{m,n}V(\frac{m}q)V(\frac{n}q)K(m)\ov{K(n)}\Delta_{q,k}(lm,n)	
\end{equation}
where $1\leq l\leq L^2$,
$$\Delta_{q,k}(lm,n)=\sum_{g\in\mcB_k(q)}\rhg(lm)\ov{\rhg(n)}$$
and
\begin{eqnarray*}\label{moment2} 
\Delta_{q,\phi}(lm,n)&&=\sum_{k \equiv 0 \mods{2},\ k>0} \dot{\phi}(k)(k-1)
\sum_{g\in\mcB_k(q)}\rhg(lm)\ov{\rhg(n)} \\
&&\ \ \   + \sum_{g\in\mcB(q)} \tilde{\phi}(t_g)\frac{4 \pi }{\cosh(\pi
    t_g)}\rhg(lm)\ov{\rhg(n)}\\
  &&\ \ \ + \,\sum_{g\in
    \mcB_E(q)}\int_{-\infty}^{\infty}\tilde{\phi}(t)\frac{1}{\cosh(\pi
    t)} \rho_g(lm,t)\ov{\rho_g(n,t)}\,dt.
\end{eqnarray*}

The Petersson-Kuznetzov formula expresses $\Delta_{q,k}(m,n)$ $\Delta_{q,\phi}(m,n)$ as sums of Kloosterman sums:
\begin{equation}
\label{Pet}
\Delta_{q,k}(m,n)=\delta_{m=n}+2\pi i^{-k}\sum_{c} \frac{1}{cq}
S(m, n;cq)J_{k-1}\left(\frac{4\pi\sqrt{mn}}{cq}\right).
\end{equation}
and
\begin{equation}
\label{Kuz}
\Delta_{q,\phi}(m,n)=\sum_{c} \frac{1}{cq}
S(m, n;cq)\phi\left(\frac{4\pi\sqrt{mn}}{cq}\right),
\end{equation}
where
$$S(m, n;cq)=\sum_{(x,cq)=1}e\left(\frac{mx+n\ov x}{cq}\right)$$
is the non-normalized Kloosterman sum of modulus $cq$ (where $x.\ov x\equiv 1\mods{cq}$). In \eqref{eqmodelsummodular}, because $m$ and $n$ are of size $q$ and $\phi$ is rapidly decreasing at $0$, the contribution of the $c\gg l^{1/2}$ is small. We will simplify further by evaluating only the contribution of $c=1$, that is
$$\frac1q\sumsum_{m,n}V(\frac{m}q)V(\frac{n}q)K(m)\ov{K(n)}S(lm,n;q)\phi(\frac{4\pi\sqrt{lmn}}{q}).$$
Our next step is  to open the Kloosterman sum and apply the Poisson summation formula on the $m$ and $n$ variables. We obtain
$$\frac{1}q\frac{q^2}{(q^{1/2})^2}\sumsum_{m^*,n^*}\what{W}(m^*,n^*)\sum_{x\in\Fqt}\what K(lx+m^*)\ov{\what K(x^{-1}+n^*)}$$
where 
$$W(x,y)=V(x)V(y)\phi(4\pi\sqrt{lxy}).$$
In particular, the Fourier transform $\what{W}(m^*,n^*)$ is very small unless $m^*+n^*\ll l$ so the above sum is over $m^*,n^*\ll l$.
Setting $$\gamma_1=\begin{pmatrix}l&m^*\\&1
\end{pmatrix},\ \gamma_2=\begin{pmatrix}n^*&1\\1&0
\end{pmatrix}
$$
we see that the $x$-sum is the correlation sum
$q\mcC(K,\gamma_2.\gamma_1^{-1})$ which is $\ll q^{1/2}$ if $\gamma_2.\gamma_1^{-1}$ does not belong to the group of automorphism of ${\what\mcF}$. Using Theorem \ref{thmautogroup} one show that if $l$ is a sufficiently small fixed (positive) power of $q$, the bound $$\sum_{x\in\Fqt}\what K(lx+m^*)\ov{\what K(x^{-1}+n^*)}\ll_{\CF} q^{1/2}$$ holds for most pairs  $(m^*,n^*)$. From this we deduce \eqref{Mbound}.

\section{The ternary divisor function in  arithmetic progressions to large moduli}\label{Secternary}
Given some arithmetic function $\lambda=(\lambda(n))_{n\geq 1}$,  a natural question in analytic number theory is to understand how well $\lambda$ is distributed in arithmetic progressions: given $q\geq 1$ and $(a,q)=1$ one would like to evaluate the sum
$$\sum_\stacksum{n\leq X}{n\equiv a\mods q}\lambda(n)$$
as $X\ra\infty$ and for $q$ as large as possible with respect to $X$. It is natural to evaluate the difference
$$E(\lambda;q,a):=\sum_\stacksum{n\leq X}{n\equiv a\mods q}\lambda(n)-\frac{1}{\vphi(q)}\sum_\stacksum{n\leq X}{(n,q)=1}\lambda(n)$$
and assuming that $\lambda$ is "essentially" bounded the target would be to obtain a bound of the shape
\begin{equation}\label{eqEbound}
E(\lambda;q,a)\ll_A \frac{X}q(\log X)^{-A}	
\end{equation}
for any $A\geq 0$, as $X\ra+\infty$ and for $q$ as large as possible compared to $X$.

The emblematic case is when $\lambda=1_\mcP $ is the characteristic function of the primes. In that case the problem can be approached through the analytic properties of Dirichlet $L$-functions and in particular the localization of their zeros. The method of Hadamard-de la Vallee-Poussin (adapted to this setting by Landau) and the Landau-Siegel theorem show that \eqref{eqEbound} is satisfied for $q\leq (\log X)^{B}$ for any given $B$, while the validity of the generalized Riemann hypothesis would give \eqref{eqEbound} for $q\ll X^{1/2-\delta}$ for any fixed $\delta>0$.  Considering averages over $q$, it is possible to reach the GRH range and this is the content of the Bombieri-Vinogradov theorem
\begin{theorem}[Bombieri-Vinogradov] For any $A\geq 0$ there exists $B=B(A)$ such that for $Q\leq X^{1/2}/\log^BX$
$$\sum_{q\leq Q}\max_{(a,q)=1}|E(1_\mcP;q,a)|\ll X/\log^AX.$$
\end{theorem}
Passing the GRH/Bombieri-Vinogradov  range and reaching the inequality $Q\leq x^{1/2+\eta}$ for some $\eta>0$ is a fundamental problem in analytic number theory with many major applications. For instance, Y. Zhang's breakthrough on the existence of bounded gaps between primes proceeded by establishing a version of the Bombieri-Vinogradov theorem going beyond the $Q=X^{1/2}$ range on average over smooth moduli. \cite{YZhang}; we will discuss some of the techniques entering his proof below. 

Several arithmetic functions are of interest besides the characteristic function of the primes or other sequences. One of the simplest are the divisor functions
$$d_k(n)=\sum_{n_1.\ldots n_k=n}1.$$
For $k=2$, Selberg and others established the following (still unsurpassed)
\begin{theorem} [The divisor function in  arithmetic progressions to large moduli]\label{thmd2} For every non-zero integer $a$, every $\eps,A>0$, every $X\geq 2$ and  every prime $q$, coprime with $a$, satisfying $$q\leq X^{2/3-\eps},$$ we have
$$
E(d_2;q,a)\ll \frac{X}{q}(\log X)^{-A},
$$
where the implied constant only depends on $\eps$ and $A$ (and not on $a$).
\end{theorem}
\proof(Sketch)
To simplify matters we consider the problem of evaluating the model sum
$$\sum_{n_1n_2\equiv a\mods q}V(\frac{n_1}{N_1})V(\frac{n_2}{N_2})$$
for $N_1N_2=X$ and $V\in\mcC^\infty_c(]1,2[)$. We apply the Poisson summation formula to the $n_1$ variable and to the $n_2$ variable. The condition $n_1n_2\equiv a\mods q$ get transformed into
$$\delta_{n_1n_2\equiv a\mods q}\ra q^{-1/2}\eq(an_1/n_2)\ra q^{-1/2}\Kl_2(an_1n_2;q).$$
The ranges the ranges $N_1,N_2$ are transformed into $$N_1^*=q/N_1,N_2^*=q/N_2$$ and the whole model sum is transformed into a sum of the shape
$$MT(a;q)+ET(a;q)$$
where $MT(a;q)$ is a main term which we will not specify (but is of the right order of magnitude), and $ET(a;q)$ is an error term of the shape
$$ET(a;q)=\frac{1}{q^{1/2}}\frac{N_1}{q^{1/2}}\frac{N_2}{q^{1/2}}\sum_{n_1,n_2}\Kl_2(an_1n_2;q)\tilde V(\frac{n_1}{N_1^*})\tilde V(\frac{n_2}{N_2^*})$$ 
where $\tilde V$ is a rapidly decreasing function. By Weil's bound for Kloosterman sums, the error term is bounded by $q^{1/2+\epsilon}$ which is smaller that $X(\log X)^{-A}/q$ as long as $X\leq q^{2/3-2\eps}$.
\qed

\begin{remark}\label{remselberg} Improving the exponent $2/3$ is tantamount to detect cancellation in the sum of Kloosterman sums above. We have given such an improvment in \eqref{eqbounddiv}; unfortunately in the present case the range of the variable $n_1n_2$ is $N_1^*N_2^*=q^2/X\leq q^{1/2}$ which is too short with current technology. See however the \cite{FoIw} for an improvement beyond the $q=x^{2/3}$ limit on average over a family of moduli $q$ admitting a specific factorisation.
\end{remark}

We now show how to go beyond the Bombieri-Vinogradov range for the specific case of the ternary divisor function
$$d_3(n)=\sum_{n_1n_2n_3=n}1$$
(in fact in a stronger form because it is not even necessary to average over the modulus $q$ !). The very first result of that kind is due to Friedlander-Iwaniec \cite{FrIw} (with $\frac12+\eta=\frac12+\frac1{231}$) and was later improved by Heath-Brown (with $\frac12+\eta=\frac12+\frac1{81}$) \cite{HBActa}. When the modulus $q$ is prime, the best result to date is to be found in \cite{FKM3}:
\begin{theorem}[The ternary divisor function in  arithmetic progressions to large moduli]\label{thmd3} For every non-zero integer $a$, every $A>0$, every $X\geq 2$ and  every prime $q$, coprime with $a$, satisfying $$q\leq X^{\frac{1}2+\frac{1}{47}},$$ we have
$$
E(d_3;q,a)\ll \frac{X}{q}(\log X)^{-A},
$$
where the implied constant only depends on $A$ (and not on $a$).
\end{theorem} 
\begin{remark} One may wonder why these higher order divisor functions are so interesting: one reason is that these problems can be considered as approximations for the case of the von Mangoldt function. Indeed, the Heath-Brown identity (Lemma \ref{lemHB}) expresses the von Mangoldt function as a linear combination of arithmetic functions involving higher divisor functions, therefore studying higher divisor functions in  arithmetic progressions  to large moduli will enable to progress on the von Mangoldt function.\footnote{This was formalised by Fouvry \cite{FouCrelle}.}

\end{remark}
\proof
We consider again a model sum of the shape
$$\sum_{n_1n_2n_3\equiv a\mods q}V(\frac{n_1}{N_1})V(\frac{n_2}{N_2})V(\frac{n_3}{N_3})$$
for $N_1N_2N_3=X$ and $V\in\mcC^\infty_c(]1,2[)$. We apply the Poisson summation formula to the variables $n_1$  $n_2$ and $n_3$. The condition $n_1n_2n_3\equiv a\mods q$ is this time transformed into the hyper-Kloosterman sum
$$\frac{1}{q^{1/2}}\Kl_3(an_1n_2n_3;q).$$
The model sum is transformed into a main term (of the correct order of magnitude) and an error term
$$ET_3(a;q)=\frac{1}{q^{1/2}}\frac{N_1}{q^{1/2}}\frac{N_2}{q^{1/2}}\frac{N_3}{q^{1/2}}\sum_{n_1,n_2,n_3}\Kl_2(an_1n_2n_3;q)\tilde V(\frac{n_1}{N_1^*})\tilde V(\frac{n_2}{N_2^*})\tilde V(\frac{n_3}{N_3^*})$$
with $$N_i^*=q/N_i,\ i=1,2,3.$$
The objective is to obtain a bound of the shape
\begin{equation}\label{d3goal}
\Sigma_3:=\sum_{n_1,n_2,n_3}\Kl_3(an_1n_2n_3;q)\tilde V(\frac{n_1}{N_1^*})\tilde V(\frac{n_2}{N_2^*})\tilde V(\frac{n_3}{N_3^*})\ll\frac{q}{\log^Aq}	
\end{equation}
for  $X=q^{2-\eta}$ for some fixed $\eta>0$ (small), or equivalently for
$$N_1^*N_2^*N_3^*=q^{1+\eta}.$$
We will show that when $\eta=0$, \eqref{d3goal} holds with the stronger bound $\ll q^{1-\delta}$
for some $\delta>0$. A variation of this argument will show \eqref{d3goal} for some positive $\eta$. Write
$$N_i^*=q^{\nu_i},\ i=1,2,3,\ \nu_1+\nu_2+\nu_3=1;$$
we assume that
$$0\leq \nu_1\leq\nu_2\leq\nu_3.$$
Suppose that $\nu_3\geq 1/2+\delta$. Then the P\'olya-Vinogradov method, applied to the $n_3$ variable, leads to a bound of the shape
$$\Sigma_3\ll q^{1-\nu_3+1/2}\log q\ll q^{1-\delta}\log q.$$

Otherwise we have
$\nu_3\leq 1/2+\delta$. We assume now that $\nu_1\geq 2\delta$; then $\nu_1\leq 1/3$, so that grouping the variables $n_2n_3$ into a single variable $n$ of size $\geq q^{2/3}$ (weighted by a divisor like function) and applying Theorem \ref{thmbilinear}, we obtain the bound
$$\Sigma_3\ll  q^{1-\delta}\log^3 q.$$
We may therefore assume that
$$\nu_1\leq 2\delta,\ \nu_2+\nu_3\geq 1-2\delta.$$
The $n_2n_3$-sum is similar to the sum in \eqref{eqbounddiv} (for $K(n)=\Kl_3(an_1n;q)$) and indeed the same bound holds, so that for any $\eps>0$, we have
$$\Sigma_3\ll_\eps q^{\nu_1+\frac{\nu_2+\nu_3}{2}+\frac12-\frac18+\epsilon}\ll_\eps q^{2\delta+1-\frac18+\epsilon}$$
which gives the required bounds if $\delta$ is chosen $<1/24$. 
\qed

\section{The geometric monodromy group and Sato-Tate laws}

In this section we discuss an important invariant attached an $\ell$-adic sheaf: its geometric  monodromy group. This will be crucial in the next section to study more advanced sums of trace functions (multicorrelation sums). Another rather appealing outcome of this notion are the {\em Sato-Tate} type laws which describe the distribution of the set of values of trace functions as $q^n$ grows.

\subsection{Sato-Tate laws for elliptic curves}

The term "Sato-Tate law" comes from the celebrated {\em Sato-Tate Conjecture} for elliptic curves over $\Qq$ which is now a theorem established in a series of papers principally by Clozel, Harris, Shepherd-Barron and Taylor \cite{CHT,HSBT,Tay,BGHT}. Let $E/\Qq$ be an elliptic curve defined over $\Qq$ with a model over $\Zz$ --for instance given by the Weierstrass equation
$$E\colon zy^2=x^3-azx^2-bz^3,\ a,b\in\Zz,\ \Delta(a,b)=4a^3-27b^2\not=0.$$
For any prime $q$, we denote by $E(\Fq)$ the reduction modulo $q$ of $E$; we have (Hasse bound) 
$$a_q(E):=q+1-|E(\Fq)|\in[-2q^{1/2},2q^{1/2}];$$
we can then define the angle $\theta_{E,q}\in [0,\pi]$ of $E$ at the prime $q$ by the formula
$$a_q(E)/q^{1/2}=2\cos(\theta_{E,q}).$$

\begin{theorem}[Sato-Tate law for an elliptic curve]\label{origST} Let $E/\Qq$ be a non-CM elliptic curve. As $X\ra\infty$, the multiset of angles
$\{\theta_{E,q},\ q\leq X,\ q\ prime\}$  becomes equidistributed  on $[0,\pi]$ with respect to the so-called Sato-Tate measure $\mu_{ST}$ whose density is given by $$d\mu_{ST}=\frac{2}\pi\sin^2(\theta)d\theta.$$
In other words, for any interval $I\subset[0,\pi]$, we have
$$\frac{|\{q\leq X,\ q\ prime,\ \theta_{E,q}\in I\}|}{\pi(X)}\ra \mu_{ST}(I)=\frac{2}\pi\int_I\sin^2(\theta)d\theta$$
as $X\ra\infty$.
\end{theorem}
The Sato-Tate measure $\mu_{ST}$ introduced in this statement has a more conceptual description: let $\SU_2(\Cc)$ be the special unitary group in two variables and let $\SU_2(\Cc)^\natural$ be its space of conjugacy classes, that space is identified with $[0,\pi]$ via the map
$$\begin{pmatrix}e^{i\theta}&0\\0&e^{-i\theta}	
\end{pmatrix}^\natural \mapsto \theta\mods \pi.
$$
The Sato-Tate measure $\mu_{ST}$ then corresponds to the direct image of the Haar measure on $\SU_2(\Cc)$ under the natural projection $\SU_2(\Cc)\mapsto \SU_2(\Cc)^\natural$: this follows from the Weyl integration formula. Now let us recall that attached to the elliptic curve $E$ is a Galois representation on its $\ell$-adic  Tate module\footnote{which is an $\ell$-adic sheaf over $\Spec(\Zz)$}
$$\rho_E\colon \Gal(\ov\Qq/\Qq)\to \GL(V_\ell(E))$$
which is unramified at every prime $q$ not dividing the discriminant (of the integral model) of $E$ and for such a prime, the Frobenius conjugacy class satisfies
$$\tr(\Frob_q|V_\ell(E))=a_q(E)=2q^{1/2}\cos(\theta_{E,q})$$
hence defines a complex conjugacy class
$$\begin{pmatrix}e^{i\theta_{E,q}}&0\\0&e^{-i\theta_{E,q}}	
\end{pmatrix}^\natural.$$
The Sato-Tate law for non-CM elliptic curves then states that this collection of Frobenius conjugacy classes becomes equidistributed relative to this measure. 

\begin{remark} For CM-elliptic curves there is also a (different) Sato-Tate law which was established by Hecke much earlier: the angles $\theta_{E,q}$ are equidistributed with respect to the uniform measure.
\end{remark}

The proof of the Sato-Tate conjecture in the non-CM case is one of the crowning achievements of the Langlands program; several decades before its proof, several variants of this conjecture have been established for {\em families} of elliptic curves over finite fields: given $a,b\in\Fq$ such that $\Delta(a,b):=4a^3-27b^2\not=0$ the  Weierstrass equation
$$E_{a,b}\colon y^2=x^3-ax^2-b$$
defines an elliptic curve over $\Fq$ and let
$$a_q(a,b)=q+1-|E_{a,b}(\Fq)|=2q^{1/2}\cos(\theta_{a,b,q}).$$
Using the Selberg trace formula, Birch \cite{BirchST}, established the following variant of the Sato-Tate law for elliptic curves
\begin{theorem} As $q\ra\infty$ the multiset of angles
$\{\theta_{a,b,q},\ (a,b)\in\Ff_q^2,\ \Delta(a,b)\not=0\}$   becomes equidistributed  on $[0,\pi]$ with respect to  $\mu_{ST}$:
for any interval $I\subset[0,\pi]$, we have
$$\frac{|\{(a,b)\in\Ff_q^2,\ \Delta(a,b)\not=0, \theta_{a,b,q}\in I\}|}{|\{(a,b)\in\Ff_q^2,\ \Delta(a,b)\not=0\}|}\ra \mu_{ST}(I),\ q\ra\infty.$$
\end{theorem}

There is another variant, spelled out by Katz and which is consequence of Deligne's work \cite{WeilII}; it concerns one parameter families of elliptic curves: let $a(T),b(T)\in\Zz[T]$ be polynomials such that $\Delta(T):=4a(T)^3+27b(T)^2\not=0$; for $q$ a sufficiently large prime, the equation over $\Fq$,
$$E_t\colon y^2=x^3-a(t)x^2-b(t)$$
defines a family of elliptic curves indexed by the set $U(\Fq):=\{t\in\Fq,\ \Delta(t)\not=0\}.$
For any $t\in U(\Fq)$ we set $$\theta_{t,q}:=\theta_{a(t),b(t),q}\in[0,\pi].$$
\begin{theorem}\label{Ellequid} Assume that the $j$-invariant $j(T)=-1728\frac{4a(T)^3}{\Delta(T)}$ is not constant, then the multiset
$\{\theta_{t,q},\ t\in U(\Fq)\}$   becomes equidistributed  on $[0,\pi]$ with respect to  $\mu_{ST}$ as $q\ra\infty$. In other words, for any interval $I\subset[0,\pi]$, we have
$$\frac{|\{t\in U(\Ff_q),\ \theta_{t,q}\in I\}|}{|U(\Fq)|}\ra \mu_{ST}(I),\ q\ra\infty.$$
\end{theorem}
\begin{remark} Deligne \cite[Proposition 3.5.7]{WeilII} proved another variant of the Sato-Tate law  when the parameter set is $U(\Fqn)$ with $q$ fixed (large enough) and $n\ra\infty$; this is in fact a special case of "Deligne's equidistribution theorem" \cite[Theorem 3.5.3]{WeilII}
\end{remark}

Theorem \ref{Ellequid} is a special case of very general Sato-Tate laws for $\ell$-adic sheaves: indeed the function
$$t\in U(\Fq)\mapsto \frac{a_q(t)}{q^{1/2}}$$
is the trace function of some geometrically irreducible $\ell$-adic sheaf $\mcE_{a,b}$ whose associated trace function is given by
\begin{equation}\label{ellsheaf}
t\mapsto -\frac{1}{q^{1/2}}\sum_{x\in\Fq}(\frac{x^3+a(t)x+b(t)}q),	
\end{equation}
 where $\left(\frac{\cdot}{q}\right)$ is the Legendre symbol. A key player for such Sato-Tate law is the

\subsection{The geometric monodromy group of a sheaf}

\begin{definition}[\cite{GKM}{ Chap. 3}] Let $\mcF$ be a sheaf pure of weight $0$ and let $\rho_\mcF$ be the associated Galois representation. The geometric (resp.~arithmetic) monodromy group $\GgeomF$ (resp.~$\GarithF$) is the Zariski closure of $\rho_\mcF(\Ggeom)$ (resp.~$\rho_\mcF(\Garith)$) inside $\GL(V_\mcF)$; in particular $$\GgeomF\subset\GarithF.$$ It follows from \cite[Th\'eor\`eme (3.4.1)]{WeilII} that the connected component  $\GgeomF^0$ of $\GgeomF$ is semisimple.
\end{definition}
\begin{example}
\begin{itemize}
\item In the case of the trace function \eqref{ellsheaf}, Deligne  showed \cite[Lemme 3.5.5]{WeilII}, that if $q>2$ and the $j$-invariant $j(T)\mods q$ is not constant, one has
$$\Ggeomd{\mcE_{a,b}}=\Garithd{\mcE_{a,b}}=\SL_2.$$
\item In his numerous books \cite{GKM,ESDE,Katzbull,MMP,TLM,KatzConvol} Katz computed the monodromy groups of various classes of sheaves: for instance, he proved in \cite[Theorem 11.1]{GKM} that for Kloosterman sheaves one has (for $q>2$)
$$\Ggeomd{\KL_k}=\Garithd{\KL_k}=\begin{cases}\SL_k&\hbox{ if $k$ is odd}\\ \Sp_{k}&\hbox{ if $k$ is even}.	
\end{cases}
$$
	
\end{itemize}

\end{example}

\subsection{Sato-Tate laws}\label{secSTlaw}

In the sequel we make the simplifying hypothesis that
\begin{equation}\label{eqGarithincluded} \Ggeomd{\mcF}=\Garithd{\mcF}.	
\end{equation}

\subsubsection{Moments of trace functions}

Before presenting the Sato-Tate laws in general, let us consider the very specific concrete problem of evaluating the {\em moments} of a trace function $K$. For $l\geq 0$ an integer, the $2l$-th moment of $K$ is the average 
$$\mcM_{2l}(K)=\frac{1}{q}\sum_{x\in\Fq}|K(x)|^{2l}.$$
The possibility of evaluating these comes from the fact that $x\mapsto |K(x)|^{2l}$ is indeed a trace function (not necessarily and in fact almost never irreducible). Indeed let $\Std\colon \Ggeomd{\mcF}\hookrightarrow \GL(V_\mcF)$ be the standard representation of the group $\Ggeomd{\mcF}$ and let
 $\rho_{l,l}$ be the representation $$\rho_{l,l}=(\Std\otimes \Std^{*})^{\otimes l}.$$
 Because of our assumption \eqref{eqGarithincluded} , the composition $$\rho_{l,l}(\mcF)"="\rho_{l,l}\circ\rho_\mcF$$ is a representation of $\Garithd{\mcF}$ hence defines an $\ell$-adic sheaf pure of weight $0$ whose trace function is\footnote{at least at the $x$ where it is lisse}
  $x\mapsto |K(x)|^{2l}.$
  
 The decomposition of this representation into irreducible representations of $\Ggeomd{\mcF}$
$$\rho_{l,l}=m_1(\rho_{l,l}).1\oplus\bigoplus_{1\not=r\in\Irr(\Ggeomd{\mcF})}m_r(\rho_{l,l}).r$$
yields a decomposition of $\rho_{l,l}(\mcF)$ into a sum of geometrically irreducible sheaves
$$\rho_{l,l}\circ \mcF=m_1(\rho_{l,l})\ov{\Ql}\oplus\bigoplus_{1\not=r\in\Irr(\Ggeomd{\mcF})}m_r(\rho_{l,l})r\circ \mcF$$ and a decomposition of
$|K(x)|^{2l}$ as a sum of trace functions
$$|K(x)|^{2l}=m_1(\rho_{l,l})+\sum_{1\not=r}m_r(\rho_{l,l})K_{r\circ\mcF}(x).$$
From Deligne's Theorem (Cor. \ref{delignecor}) one deduce that
$$\frac{1}{q}\sum_{x}|K(x)|^{2l}=m_1(\rho_{l,l})+O_{\CF,l}(q^{-1/2})$$
where $m_1(\rho_{l,l})$ is the multiplicity of the trivial representation
 in the representation $(\Std\otimes \Std^{*})^{\otimes l}$ of $\Ggeomd{\mcF}$. In the same way, we could evaluate (in terms of the representation theory of the group $\Ggeomd{\mcF}$) more general moments like
$$ \frac{1}{q}\sum_{x\in\Fq}|K(x)|^{2l}K(x)^{l'}$$
for integers $l,l'\geq 0$.

\subsubsection{Equidistribution of Frobenius conjugacy classes}

There is a more conceptual interpretation of these moments. For any $x\in U(\Fq)$, the Frobenius at $x$  acting on $V_\mcF$ produces a $\rho_\mcF(\Garith)$-conjugacy class $$\rho_\mcF(\Frob_x)\subset\Garithd{\mcF}(\Cc)=\GgeomF(\Cc).$$ The {\em Frobenius conjugacy class} of $\mcF$ at $x$ is by definition the $\GgeomF(\Cc)$-conjugacy class of its semisimple part (in the sense of Jordan decomposition) and is denoted $\theta_{x,\mcF}.$
 Let $K$ be any maximal compact subgroup of $\Ggeomd{\mcF}(\Cc)$ and $K^\natural$ its space of conjugacy classes. As explained in \cite{GKM}(Chap. 3), the conjugacy class 
$\theta_{x,\mcF}$ defines a unique conjugacy class in $K$, also denoted $\theta_{x,\mcF}\in K^\natural$.
 The Sato-tate laws describe the distribution of the set
$\{\theta_{x,\mcF},\ x\in U(\Fq)\}\hbox{ inside $K^\natural$}$ as $q\ra\infty$.

More precisely, let $G$ be a connected semisimple algebraic group over ${\ovQl}$ and $K\subset G(\Cc)$ a maximal compact subgroup. Let $\mu^\natural$ be the direct image of the Haar probability measure on $K$ under the projection $K\mapsto K^\natural$.

\begin{theorem}[Sato-Tate law] Let $G$ and $K\subset G(\Cc)$ as above. Suppose we are given a sequence of primes $q\ra\infty$ and for each such prime some $\ell$-adic sheaf $\mcF$ over $\Fq$, satisfying \eqref{eqGarithincluded},  whose conductor $\CF$ is bounded independently of $q$, such that $$\Ggeomd{\mcF}=\GarithF =G.$$
For any such $q$ and $x\in U(\Fq)$ let $\theta_{x,\mcF}\in K^\natural$ be the conjugacy class of $\mcF$ at $x$ relative to $K$.
	
	As $q\ra\infty$ the sets of conjugacy classes
	$$\{\theta_{x,\mcF},\ x\in U(\Fq)\}$$ become equidistributed with respect to the  measure $\mu^\natural$:
the probability measure
	$$\frac{1}{|U(\Fq)|}\sum_{x\in U(\Fq)}\delta_{\theta_{x,\mcF}}$$
	 converges weakly to $\mu^\natural$.	 In other words, for any $f\in\mcC(K^\natural)$
\begin{equation}\label{eqequid}
\frac{1}{|U(\Fq)|}\sum_{x\in\Fq}f(\theta_{x,\mcF})\ra \int_{K^\natural}f(\theta)d\mu^\natural(\theta),\ q\ra\infty.	
\end{equation}

\end{theorem}

\proof By the Peter-Weyl theorem, the functions
$$\tr(r)\colon \theta\in K^\natural\to \tr(r(\theta))\in\Cc$$
when $r$ ranges over all the irreducible representations of $G$, form an orthonormal basis of $L^2(K^\natural,\mu^\natural)$ and generates a dense subspace of the space of continuous functions on $K^\natural$. By Weyl equidistribution criterion it is therefore sufficient to show that for any $r$ irreducible and non-trivial, one has
$$\frac{1}{|U(\Fq)|}\sum_{x\in U(\Fq)}\tr(r(\theta_{x,\mcF}))\ra \mu^\natural(\tr(r))=0.$$
The function
$$K_{r,\mcF}\colon x\in U(\Fq)\to r(\theta_{x,\mcF})$$ is the trace function associated to the sheaf $r\circ\mcF$ corresponding to the representation of $\GarithF$, $r\circ\rho_{\mcF}$ (because of \eqref{eqGarithincluded} this composition is well defined). That sheaf is by construction geometrically irreducible, non-trivial and its conductor is bounded in terms of $C(\mcF)$ and $r$ only, so it follows from Deligne's Theorem that
$$\frac{1}{|U(\Fq)|}\sum_{x\in U(\Fq)}\tr(r(\theta_{x,\mcF}))\ll_{C(\mcF),r}q^{-1/2}\ra 0.$$ 
\qed

\subsubsection{The case of Kloosterman sums} As we have seen above, for the Kloosterman sums $\Kl_2(x;q)$, we have
$$G=\Sp_2=\SL_2,\  K=\mathrm{SU}_2(\Cc)$$ and, via the identification $K^\natural\simeq[0,\pi]$,  the  measure $\mu^{\natural}$ is identified with the { Sato-Tate} measure $\mu_{ST}$.

For $x\in\Fqt$, we define the angle $\theta_{q,x}\in [0,\pi]$ of the Kloosterman sum $\Kl_2(x;q)$ as 
$$\Kl_2(x;q)=\tr\begin{pmatrix}e^{i\theta_{q,x}}&0\\0&e^{-i\theta_{q,x}}
\end{pmatrix}=2\cos(\theta_{q,x}).
$$

The Sato-Tate law becomes the following explicit statement (due to Katz):
\begin{theorem}[Sato-Tate law for Kloosterman sums] For any interval $I\subset[0,\pi]$
$$\frac{1}{q-1}|\{x\in\Fqt,\ \theta_{q,x}\in I\}|\ra \frac{2}\pi\int_I\sin^2(\theta)d\theta,\ q\ra\infty.$$
\end{theorem}
The above Sato-Tate law is called "vertical" as it describes  the distribution of Kloosterman sums with varying parameters $x\in\Fqt$ as $q\ra\infty$; such law is analogous to the Sato-Tate law of Theorem \ref{Ellequid}.

In \cite{Sommes}, Katz in analogy with the original Sato-Tate conjecture (Theorem \ref{origST})  asked for the distribution of the Kloosterman sums for a fixed value of the parameter (say $x=1$) and for  a varying prime modulus $q$. Katz made the following
\begin{conjecture}[Horizontal Sato-Tate law for Kloosterman sums]\label{katzconj} As $X\ra\infty$, the multiset of Kloosterman angles
$\{\theta_{q,1},\ q\leq X,\ prime\}$  becomes equidistributed with respect to the Sato-Tate measure:
for any $[a,b]\subset[0,\pi]$, we have
$$\frac{1}{\pi(X)}|\{q\leq X,\ q\ prime,\ \theta_{q,1}\in[a,b]\}|\ra \frac{2}\pi\int_a^b\sin^2(\theta)d\theta$$
as $X\ra\infty$.
\end{conjecture}

\begin{remark} There are other variants of this vertical equidistribution conjecture that have been established recently:

\begin{itemize}
\item 	Heath-Brown and Patterson \cite{HBPatt} have proven that the angles of cubic Gauss sums of varying prime moduli are equidistributed with respect to the uniform measure. 
\item Even closer to the current discussion, Duke, Friedlander and Iwaniec \cite{DFISalie} have proven the vertical equidistribution of the angles $\theta^S_{q,1}$ of {\em Sali\'e} sums defined by
$$S(1;q):=\frac{1}{q^{1/2}}\sum_\stacksum{x,y\in\Fqt}{xy=1}(\frac{x}{q})e\left(\frac{x+y}{q}\right)=:2\cos(\theta^S_{q,1})$$
again with respect to the uniform measure.
\end{itemize}
 \end{remark}

\subsection{Towards the horizontal Sato-Tate conjecture for almost prime moduli}

Unlike the original Sato-Tate conjecture the prospect for a proof of Conjecture \ref{katzconj} seem very distant at the moment.
Even the following very basic consequences of this conjecture seem today completely out of reach:
\begin{itemize}
\item There exist infinitely many primes $q$ such that $|\Kl_2(1;q)|\geq 2017^{-2017}$,
\item There exist infinitely many primes $q$ such that $\Kl_2(1;q)>0$ (resp.~$\Kl_2(1;q)<0$)
\end{itemize}
In this section we will explain how some of the results discussed so far enable to say something non-trivial as the cost of replacing the prime moduli $q$ by {\em almost prime} moduli (that is squarefree-integers with an absolutely bounded number of prime factors). 

Recall that for $c\geq 1$ a squarefree integer and $(a,c)=1$ the normalized Kloosterman sum of modulus $c$ and parameter $a$ is
$$\Kl_2(a;c)=\frac{1}{c^{1/2}}\sum_{x\in(\Zz/c\Zz)^\times}e\left(\frac{\ov x+ax}c\right).$$
By the Chinese remainder theorem, Kloosterman sums satisfy the {\em twisted multiplicativity} relation: for $c=c_1c_2$, $(c_1,c_2)=1$ one has
\begin{equation}\label{twistedmult}
\Kl_2(a;c)=\Kl_2(a\ov{c_2}^2;c_1)\Kl_2(a\ov{c_1}^2;c_2)
\end{equation}
so that by Weil's bound one has
$$|\Kl_2(a;c)|\leq 2^{\omega(c)}$$
where $\omega(c)$ is the number of prime factors of $c$. We can then define the corresponding Kloosterman angle by
$$\cos(\theta_{c,a})=\frac{\Kl_2(a;c)}{2^{\omega(c)}}.$$
It is then natural to make the following
\begin{conjecture}[Horizontal Sato-Tate law for Kloosterman sums with composite moduli]\label{katzconjk} Given  $k\geq 1$ un integer, let $\pi_k(X)$ be the number of squarefree integers $\leq X$ with exactly $k$ prime factors and let $\mu_{ST,k}$ be the Sato-Tate measure of order $k$, defined as the push-forward of the measure $\mu_{ST}^{\otimes k}$ on $[0,\pi]^k$ by the map
$$(\theta_1,\ldots,\theta_k)\in [0,\pi]^k\mapsto \arccos(\cos(\theta_1)\times\ldots\times\cos(\theta_k)))\in[0,\pi].$$
 for any $k\geq 1$, the multiset of Kloosterman angles
$$\{\theta_{c,1},\ c\leq X,\ c\hbox{ is squarefree with $k$ prime factors}\}$$  becomes equidistributed with respect to $\mu_{ST,k}$ as $X\ra\infty$.
\end{conjecture}
This conjecture for any $k\geq 2$ seem as hard as the original one (and is not implies by it). On the other hand it is possible to establish some of its consequences:

\begin{theorem} There exists $k\geq 2$  such that
\begin{enumerate} 
\item 	for infinitely many square-free integers $c$ with at most $k$ prime factors,$$|\Kl_2(1;c)|\geq 2017^{-2017};$$
 \item 	for infinitely many square-free integers $c$ with at most $k$ prime factors,$$\Kl_2(1;c)>0;$$
\item 	for infinitely many square-free integers $c$ with at most $k$ prime factors,$$\Kl_2(1;c)<0.$$
\end{enumerate}
\end{theorem}
The first statement above was proven in \cite{Mi1} for $k=2$ (with $2017^{-2017}$ replaced by $4/25$; the second and the third were first proven in \cite{FouMiAnn} for $k=23$; this value was subsequently improved by Sivak, Matom\"aki and Ping who holds the current record with $k=7$ \cite{Sivak,Matomaki,Xi,Xi2}. 

\subsubsection{Kloosterman sums can be large} We start with the first statement which we prove for $c=pq$ a product of two distinct primes. The main idea is to use the twisted multiplicativity relation
$$\Kl_2(1;pq)=\Kl_2(\ov p^2;q)\Kl_2(\ov q^2;p)$$
and to establish the existence of some $\kappa$ for which there exist infinitely many pairs of distinct primes $(p,q)$ such that
$$|\Kl_2(\ov p^2;q)|\ |\Kl_2(\ov q^2;p)|\geq \kappa.$$
Indeed, for such pairs we have
$$|\Kl_2(1;pq)|\geq \kappa^2.$$
Given $X$ large, we will consider pairs $(p,q)$ such that $p,q\in[X^{1/2},2X^{1/2}[$ and will show that for $\kappa$ small enough the two sets
$$\{(p,q),\ p\not=q\in[X^{1/2},2X^{1/2}[,\ p,q\hbox{ primes}\ |\Kl_2(\ov p^2;q)|\geq \kappa\}$$
$$\{(p,q),\ p\not=q\in[X^{1/2},2X^{1/2}[,\ p,q\hbox{ primes}\ |\Kl_2(\ov q^2;p)|\geq \kappa\}$$
are large enough to have a non-empty (and in fact large) intersection as $X\ra\infty$.
This is a consequence of the following equidistribution statement
\begin{proposition} Given $X\geq 1$, and a prime $q\in[X^{1/2},2X^{1/2}]$, the (multi)-set of Kloosterman angles
$$\{\theta_{q,\ov p^{2}},\ p\in[X^{1/2},2X^{1/2}[,\ p\hbox{ prime,}\ p\not=q\}$$ is equidistributed with respect to  the Sato-Tate measure: for any interval $[a,b]\subset [0,\pi]$
$$\frac{|\{p\in[X^{1/2},2X^{1/2}[,\ p\not=q\hbox{ prime,}\ \theta_{q,\ov p^{2}}\in[a,b]\}|}{|\{p\in[X^{1/2},2X^{1/2}[,\ p\not=q\hbox{ prime}\}|}\ra \frac{2}\pi\int_a^b\sin^2(\theta)d\theta$$
 as $X\ra\infty$. 
\end{proposition}

\proof We consider the pull-back sheaf $\mcK:=[x\ra x^{-2}]^*\KL_2$ whose trace function is given by
$x\ra \Kl_2(\ov x^2;q)$. As a representation of the geometric Galois group, it corresponds to restricting the representation $\KL_2$ to a subgroup of index $2$. Since the geometric monodromy group of $\KL_2$ is $\SL_2$, the same is true for the pull-back (the algebriac group $\SL_2$ has no non-trivial finite-index subgroups); therefore
$$\Ggeomd{\mcK}=\Garithd{\mcK}=\SL_2.$$
The non-trivial irreducible representations of $\SL_2$ are the symmetric powers of the standard representation, $\Sym_k(\Std),\ k\geq 1$. Given $k\geq 1$ the composed sheaf $$\mcK_k=\Sym_k\circ\mcK$$ is by construction geometrically irreducible, has rank $k+1$ with conductor bounded in terms of $k$ only and its trace function equals
$$K_k(x)=\tr(\Sym_k\begin{pmatrix}e^{i\theta_{q,\ov x^{2}}}&0\\0&e^{-i\theta_{q,\ov x^{2}}}
\end{pmatrix})=\sum_{j=0}^ke^{i(k-j)\theta_{q,\ov x^{2}}}e^{-ij\theta_{q,\ov x^{2}}}=\frac{\sin((k+1)\theta_{q,\ov x^{2}})}{\sin(\theta_{q,\ov x^{2}})}.$$
 In particular $\mcK_k$ cannot be geometrically isomorphic to any tensor product of an Artin-Schreier sheaf and a Kummer sheaf (as they have rank $1$). Hence by a simple variant of Theorem \ref{thmprimesumthm} we obtain that
 $$\frac{1}{\pi(2X^{1/2})-\pi(X^{1/2})}\sum_\stacksum{p\not=q}{p\sim X^{1/2}}K_k(p)\ra 0=\frac{2}\pi\int_0^\pi\frac{\sin((k+1)\theta)}{\sin(\theta)}\sin^2(\theta)d\theta$$
\qed

Averaging over $q$, we deduce the existence of some $\kappa>0$ ($\kappa=0,4$) such that for $X$ large enough
$$\frac{|\{(p,q),\ p\not=q\in[X^{1/2},2X^{1/2}[,\ p,q\hbox{ primes},\ |\Kl_2(\ov p^2;q)|\geq \kappa\}|}{|\{(p,q),\ p\not=q\in[X^{1/2},2X^{1/2}[,\ p,q\hbox{ primes}\}|}\geq 0,51$$
hence
\begin{equation}\label{lowerKl1}
{|\{(p,q),\ p\not=q\in[X^{1/2},2X^{1/2}[,\ p,q\hbox{ primes}\ |\Kl_2(1;pq)|\geq \kappa^2\}|}\geq (0,01+o(1))\frac{X}{(\frac12\log X)^2}.	
\end{equation}

\subsubsection{Kloosterman sums change sign}
We now discuss briefly the proof of the remaining two statements: to establish the existence of sign changes, it suffices to prove that given $V\in\mcC_c^\infty(]1,2[)$ some non-zero non-negative smooth function, there exists $u>0$ such that, for $X$ large enough
\begin{equation}\label{compareKl}
\bigl|\sum_\stacksum{c\geq 1}{p|c\Rightarrow p\geq X^{1/u}}\Kl_2(1;c)V(\frac{c}X)\bigr|<\sum_\stacksum{c\geq 1}{p|c\Rightarrow p\geq X^{1/u}}|\Kl_2(1;c)|V(\frac{c}X).	
\end{equation}
which will prove the existence of sign changes for Kloosterman sums $\Kl_2(1;c)$ whose modulus has at most $1/u$ prime factors.
Using sieve methods and the Petersson-Kuznetzov formulas to express sums of Kloosterman sums in terms of Fourier coefficients of modular forms (\eqref{Pet} and \eqref{Kuz}) and using the theory of automorphic forms, one can show that (see \cite{FouMiAnn} for a proof)

\begin{proposition} For any $\eta>0$, there exists $u=u(\eta)>0$ such that
$$\bigl|\sum_\stacksum{c\geq 1}{p|c\Rightarrow p\geq X^{1/u}}\Kl_2(1;c)V(\frac{c}X)\bigr|\leq \eta\frac{X}{\log X}$$
for $X$ large enough (depending on $\eta$ and $V$).
\end{proposition}

To conclude, it is sufficient to show that for some $u=u_0$, one has
\begin{equation}\label{lowerwish}
\sum_\stacksum{c\geq 1}{p|c\Rightarrow p\geq X^{1/u}}|\mu^2(c)\Kl_2(1;c)|V(\frac{c}X)\gg_{V} \frac{X}{\log X}	
\end{equation}
(the left-hand side is an increasing function of $u$ so the above inequality remains valid for any $u\geq u_0$). The inequality \eqref{lowerKl1} points in the right direction (for $u_0=2$), however as stated it is off by a factor $\log X\log\log X$. 
One can however recover this factor $\log X$ entierely and prove the lower bound 
$$\sum_\stacksum{c\geq 1}{p|c\Rightarrow p\geq X^{3/8}}\mu^2(c)|\Kl_2(1;c)|V(\frac{c}X)\gg_{V} \frac{X}{\log X}.$$
The reason is that Theorem \ref{thmprimesumthm} applies also when $p$ is significantly smaller than $q$ ( if $q\simeq X^{1/2+\delta}$ one can obtain a non-trivial bound  in \eqref{primesumsmooth} for $p$ of size $X^{1/2-\delta}$ for $\delta\in[0,1/8[$). The details involve making a partition of unity and we leave it to the interested reader. Another possibility (the one followed originally in \cite{FouMiAnn}) is to establish the lower bound \eqref{lowerwish} for a suitable $u$ by restricting to moduli $c$ which are products of exactly three prime factors, using the techniques discussed so far.

 \section{Multicorrelation of trace functions}\label{multisec}
So far we have mainly discussed the evaluation of correlation sums associated to two trace functions $K_1$ and $K_2$ (especially the case $K_1=K$ and $K_2=\gamma^*K$), namely
$$\mcC(K_1,K_2)=\frac{1}q\sum_{x}K_1(x)\ov{K_2(x)}.$$

In many applications, multiple correlation sums occur: sums of the shape
$$\mcC(K_1,K_2,\ldots, K_{L}):=\frac{1}q\sum_{x}K_1(x)K_2(x)\ldots K_L(x)$$
where the $K_i,\ i=1,\ldots, L$ are trace functions; of course rewriting the inner term of the sum above as a product of two factors reduces to evaluating a double correlation sum, say associated to the sheaves 
$$\mcF=\mcK_1\otimes\ldots\mcK_{l},\ \mcG=\mcK_{l+1}\otimes\ldots\mcK_{L}$$
but it would remain to determine if $\mcF$ and $\mcG$ share a common irreducible component and this may be a hard task. In practice, the multicorrelation sums that occur (due to the application of some H\"older inequality and of the P\'olya-Vinogradov method) are often of the shape
$$\mcC(K,\uple{\gamma},h)=\frac{1}q\sum_{x}K(\gamma_1\cdot x)\ldots K(\gamma_l\cdot x)\ov{K(\gamma'_1\cdot x)\ldots K(\gamma'_l\cdot x)}\eq(xh)$$
for $K$ the trace function of some geometrically irreducible sheaf $\mcF$, pure of weight $0$, 
$${\uple{\gamma}}=(\gamma_1,\ldots,\gamma_l,\gamma'_1,\ldots,\gamma'_l)\in\PGL_2(\Fq)^{2l}$$   and some $h\in\Fq$.

 This sum is the correlation associated to the trace functions of the sheaves
$$\gamma_1^*\mcF\otimes\ldots\otimes\gamma_l^*\mcF\ \hbox{and }{\gamma'}_1^*\mcF\otimes\ldots\otimes{\gamma'}_l^*\mcF\otimes\mcL_\psi$$
whose conductors are bounded polynomially in terms of $C(\mcF)$. If $\mcF$ has rank one, the two sheaves above have rank one and it is usually not difficult to determine whether these sheaves are geometrically isomorphic or not.

For $\mcF$ of higher rank, we describe a method due to Katz which has been axiomatized in \cite{FKMSP}: this method rests on the notion of geometric monodromy group which we discussed in the previous section.

 \subsection{A theorem on sums of products of trace functions}
 In this section we discuss some general result making it possible to evaluate multicorrelations sums of trace functions of interest for analytic number theory. The method is basically due to Katz and was used on several occasions, for instance in \cite{Mi1,FoMi}. The general result presented here is a special case of the results of \cite{FKMSP}. For this we need to introduce the following variants of the group of automorphism of a sheaf: one is the group of projective automorphisms
$$\Aut^p_\mcF(\Fq)=\{\gamma\in\PGLd(\Fq),\ \exists \hbox{ some rank one sheaf $\mcL$ s.t. }\gamma^*\mcF\simeq_{geom}\mcF\otimes\mcL\},$$
the other is the right-$\Aut^p_\mcF(\Fq)$-orbit
$$\Aut^d_\mcF(\Fq)=\{\gamma\in\PGLd(\Fq),\ \exists \hbox{ some rank one sheaf $\mcL$ s.t. }\gamma^*\mcF\simeq_{geom}D(\mcF)\otimes\mcL\}.$$
Let $\mcF$ be a weight $0$, rank $k$, irreducible sheaf. We assume that
\begin{itemize}
\item the geometric monodromy group equals $\Ggeomd{\mcF}=\SL_k\hbox{ or }\Sp_k,$ (we then say that $\mcF$ is of $\SL$ or $\Sp$-{\em type}),
\item the equality \eqref{eqGarithincluded} holds,
\item 	$\Aut^p_\mcF(\Fq)=\{\Id\}$; in particular $\Aut^d_\mcF(\Fq)$ is either empty or is reduced to a single element, $\xi_\mcF$ which is a possibly trivial involution ($\xi_\mcF^2=\Id$) and is called the {\em special involution.}
\end{itemize}

\begin{example} The Kloosterman sheaves $\KL_k$ have this property \cite{GKM}. The special involution is either $\Id$ if $k$ is even ($\KL_k$ is self-dual) or the matrix 
$\xi=\begin{pmatrix}-1&\\&1	
\end{pmatrix}$ for $k$ odd.
\end{example}

Finally we introduce the following ad-hoc definition:
\begin{definition} Given $${\uple{\gamma}}=(\gamma_1,\ldots,\gamma_l,\gamma'_1,\ldots,\gamma'_l)\in\PGL_2(\Fq)^{2l},$$ one says that
\begin{itemize}
\item ${\uple{\gamma}}$ is normal if there is $\gamma\in\PGLd(\Fq)$ such that $$|\{i,\ \gamma_i=\gamma\}|+|\{j,\ \gamma'_j=\gamma\}|\equiv 1\mods 2.$$
\item For $k\geq 3$,	${\uple{\gamma}}$ is $k$-normal if there exists $\gamma\in\PGLd(\Fq)$ such that
$$|\{i,\ \gamma_i=\gamma\}|-|\{\gamma'_j=\gamma\}|\not\equiv 0\mods k.$$
\item For $k\geq 3$, and $\xi\in\PGLd(\Fq)$ a non-trivial involution, ${\uple{\gamma}}$ is $k$-normal w.r.t. $\xi$ if there exist $\gamma\in\PGLd(\Fq)$ such that
 $$|\{i,\ \gamma_i=\gamma\}|+|\{j,\ \gamma'_j=\xi\gamma\}|-|\{j,\ \gamma'_j=\gamma\}|-|\{i,\ \gamma_i=\xi\gamma\}|\not\equiv 0\mods{k}.$$
\end{itemize}

\end{definition}

\begin{theorem}\label{cor-concrete}
  Let $K$ be the trace function of a
   sheaf $\mcF$ as above, $l\geq 1$, $\uple{\gamma}\in\PGL_2(\Fq)^{2l}$ and $h\in\Fq$. We assume that either 
   \par
\emph{(1)} the sheaf $\sheaf{F}$ is self-dual (so that $K$ is real-valued) and $\uple{\gamma}$ is normal
\par
\emph{(2)} the $\sheaf{F}$ is of $\SL$-type of rank $k\geq 3$,
$q>r$, and $\uple{\gamma}$ is
$k$-normal or $k$-normal w.r.t. the special involution of
$\sheaf{F}$, if it exists.
\par
\emph{(3)} or $h\not=0$.

We have
$$\mcC(K,\uple{\gamma},h)=\frac{1}q\sum_{x}K(\gamma_1\cdot x)\ldots K(\gamma_l\cdot x)\ov{K(\gamma'_1\cdot x)\ldots K(\gamma'_l\cdot x)}\eq(xh)\ll_{l,\CF}\frac{1}\qde.$$

\end{theorem}

\proof We discuss the proof only in the self-dual case for simplicity. We group together identical $\gamma_i,\gamma_j'$ and the sum becomes
$$\frac{1}q\sum_{x}K(\gamma''_1\cdot x)^{m_1}\ldots K(\gamma''_{t}\cdot x)^{m_t}\eq(xh)$$
where $t\leq 2l$, the $\gamma''_i$ are distinct and by hypothesis one of the $m_i$ is odd. The above sum is associated to the trace function of the sheaf
$$\bigotimes_{i=1}^t \Std(\gamma_i''^*\mcF)^{\otimes{m_i}}\otimes \mcL_\psi$$
where $\psi(\cdot)=\eq(h\cdot)$ and $\Std$ is the tautological representation.
We decompose each representation into irreducible
$$\rho_{m,0}=\Std(G)^{\otimes m}=\sum_{r}m_r(\rho_{m,0})r$$
and are reduced to considering various sheaves of the shape
\begin{equation}\label{eqlesheaf}
\bigotimes_{i=1}^t r_i(\gamma_i''^*\mcF)\otimes \mcL_\psi
\end{equation}
where  $(r_i)_{i\leq t}$ is a tuple of irreducible representations of $G$; by our hypothesis, we know that either $\mcL_\psi$ is not trivial or at least one of the $r_i$ is not trivial (and necessarily of dimension $>1$).

It is then sufficient to show that, under these assumptions, the sheaves \eqref{eqlesheaf} are irreducible. For this we consider the direct sum sheaf
$$\bigoplus_i\gamma_i''^*\mcF$$
and let $\Ggeomd{\oplus}\subset \prod_i G$ be the Zariski closure of the image of $\Ggeom$ under the sum of representations.
The following very useful criterion is due to Katz

\begin{theorem}[Goursat-Kolchin-Ribet criterion] Let $(\mcF_i)_i$ be a tuple of geometrically irreducible sheaves lisse on $U\subset \Aa^1_{\Fq}$, pure of weight $0$, with geometric monodromy groups $G_i$. We assume that
\begin{itemize}
	\item For every $i$, $G_i=\Sp_{k_i}$ or $\SL_{k_i}$,
	\item for any rank $1$ sheaf $\mcL$ and any $i\not=j$ there is no geometric isomorphism between $\mcF_i\otimes\mcL$ and $\mcF_j$,
	\item for any rank $1$ sheaf $\mcL$ and any $i\not=j$ there is no geometric isomorphism between $\mcF_i\otimes\mcL$ and $D(\mcF_j)$.
\end{itemize}
Then the geometric monodromy group of the sheaf
$\bigoplus_i\mcF_i$ equals $\prod_i G_i$ .
\end{theorem}
Our assumptions (the projective automorphism group of $\mcF$ is trivial, $\uple{\gamma}$ is normal
and the geometric monodromy group is either $\SL$ or $\Sp$) imply that the above criterion holds and this implies that
$$\bigotimes_i r_i(\gamma_i''^*\mcF)\otimes \mcL_\psi$$
is always irreducible.
\qed

\subsection{Application to non-vanishing of Dirichlet $L$-functions}
We now discuss a beautiful application of bounds for multicorrelation sums due to R. Khan and H. Ngo \cite{KhNg}. It concerns the proportion of non-vanishing of Dirichlet $L$-functions at the central point $1/2$. The interest in this kind of problems from analytic number theory was renewed with the work of Iwaniec and Sarnak in their celebrated attempt to prove the non-existence of a Landau-Siegel zero \cite{IS1}. Their approach was based on the following general problem: {\em given  a  family of $L$-functions
$$\{L(f,s)=\sum_{n\geq 1}\frac{\lf(n)}{n^s},\ f\in\mcF\}$$indexed by a "reasonable" family of automorphic forms $\mcF$\footnote{A reasonable definition of the notion of "reasonable" can be found in \cite{K,SST}}, show that for many $f\in\mcF$, one has $$L(f,1/2)\not=0.$$ }

	In their work \cite{IS1}, Iwaniec and Sarnak  showed specifically that for $\mcF=\mcS_2(q)$ the set of holomorphic new-forms of weight $2$ and  prime level $q$ (with trivial nebentypus), if one could show that for  $q$ large enough at least $(25+2017^{-2017})\%$ of the central $L$-values $L(f,1/2)$ do not vanish (more precisely that at least $(25+2017^{-2017})\%$ of these central values are larger than $\log^{-2017}q$ ) then there would be no Landau-Siegel zero. They eventually proved 
	\begin{theorem}[\cite{IS1}] As $q\ra\infty$ along the primes one has
$$\frac{|\{f\in \mcS_2(q),\ L(f,1/2)\geq \log^{-2}q\}|}{|\mcS_2(q)|}\geq 1/4-o(1).$$
\end{theorem}
This is "just" at the limit.
	
The possibility of producing a positive proportion of non-vanishing is not limited to this specific family and one of the most powerful and general tools  to achieve this is via the {\em mollification method}. The principle of mollification method is as follows: given the family $\mcF$, one considers for some parameter $L\geq 1$ and some suitable vector $\bfx_L=(x_\ell)_{\ell\leq L}\in\Cc^{\ell}$ the linear 
form
\begin{equation}\label{linear}
\mcL(\mcF,\bfx_L):=\frac{1}{|\mcF|}\sum_{f\in\mcF}L(f,1/2)M(f,\bfx_L)	
\end{equation}
and the quadratic form
\begin{equation}\label{quadratic}
\mcQ(\mcF,\bfx_L):=\frac{1}{|\mcF|}\sum_{f\in\mcF}|L(f,1/2)M(f,\bfx_L)|^2
\end{equation}
where $M(f,\bfx_L)$ is the linear form (called "mollifier")
$$M(f,\bfx_L)=\sum_{\ell\leq L}\frac{\lf(\ell)}{\ell^{1/2}}x_\ell$$
and the $x_\ell$ are coefficients to be chosen in an optimal way with the idea of approximating the inverse $L(f,1/2)^{-1}$. Such coefficients are almost bounded,  i.e.~satisfy: 
$$x_\ell= |\mcF|^{o(1)}.$$
By Cauchy's inequality one has
$$\frac{|\{f\in\mcF,\ L(f,1/2)\not=0\}|}{{|\mcF|}}\geq \frac{|\mcL(\mcF,\bfx_L)|^2}{\mcQ(\mcF,\bfx_L)}.$$
For suitable families one can evaluate asymptotically $\mcL(\mcF,\bfx_L)$ and $\mcQ(\mcF,\bfx_L)$ (the hard case being $\mcQ$) when $L=|\mcF|^\lambda$ for $\lambda>0$ some fixed constant
and (upon minimizing $\mcQ(\mcF,\bfx_L)$ with respect to $\mcL(\mcF,\bfx_L)$) one usually shows that 
\begin{equation}\label{ratiomolli}
\frac{|\mcL(\mcF,\bfx_L)|^2}{\mcQ(\mcF,\bfx_L)}=F(\lambda)+o(1)
\end{equation}
for $F$ some increasing rational fraction with $F(0)=0$. In \cite{IS1}, Iwaniec and Sarnak have also implemented this strategy for the (simpler) family of Dirichlet $L$-functions of modulus $q$
$$\{L(\chi,s)=\sum_{n\geq 1}\frac{\chi(n)}{n^s},\ \chi\in\what{(\Zz/q\Zz)^\times}\}$$
and were able to evaluate \eqref{linear} and \eqref{quadratic} for any $\lambda<1/2$ and to prove \eqref{ratiomolli} with
$$F(\lambda)=\frac{\lambda}{\lambda+1}$$
hence: 
\begin{theorem}[\cite{IS2}] As $q\ra\infty$ along the primes one has
$$\frac{|\{\chi\mods q,\ L(\chi,1/2)\not=0\}|}{|\{\chi\mods q\}|}\geq 1/3-o(1).$$
\end{theorem}
 Thus  the proportion of non-vanishing  can be arbitrarily close to $33.33\dots\%$.
Shortly after, Michel and Vanderkam \cite{MvdK} obtained the same proportion by a slightly different method: taking into account the fact that for a complex character, the $L$-function $L(\chi,s)$ is not self-dual ($L(\chi,s)\not=L(\ov\chi,s)$) and has root number  
$$\eps_\chi=i^{\mfa}\frac{\tau(\chi)}{q^{1/2}},\ \mfa=\frac{\chi(-1)-1}2$$ were $\tau(\chi)$ is the  Gauss sum, they introduced a symmetrized mollifier of the shape
$$M^s(\chi,\bfx_L)=M(\chi,\bfx_L)+\ov{\eps_\chi}M(\ov \chi,\bfx_L)=\sum_{\ell\leq L}\frac{\chi(\ell)+\ov{\eps_\chi}.\ov\chi(\ell)}{\ell^{1/2}}x_\ell.$$
Because of the oscillation of the root number $\eps_\chi$, they could evaluate \eqref{quadratic} only in the shorter range $\lambda<1/4$. However this weaker range is offset by the fact that the symmetrized mollifier is more effective: indeed the rational fraction $F(\lambda)$ is then replaced by
$$F^s(\lambda)=\frac{2\lambda}{2\lambda+1}$$
which takes value $1/3$ at $\lambda=1/4$. 

Recently R. Khan and H. Ngo founds a better method to bound the exponential sums considered in \cite{MvdK} building on Theorem \ref{cor-concrete} and they increased the allowed range from $\lambda<1/4$ to $\lambda<3/10$:
\begin{theorem}[\cite{KhNg}]\label{KhNgthm} As $q\ra\infty$ along the primes one has
$$\frac{|\{\chi\mods q,\ L(\chi,1/2)\not=0\}|}{|\{\chi\mods q\}|}\geq 3/8-o(1).$$
\end{theorem}
The key step in their proof is the asymptotic evaluation of the second mollified moment
\begin{equation}\label{2ndchi}
\frac{1}{\vphi(q)}\sum_{\chi\mods q}|L(\chi,1/2)|^2|M^s(\chi,\bfx_L)|^2
\end{equation}
for $L=q^\lambda$, and any fixed $\lambda<3/10$. By (nowadays) standard methods\footnote{inappropriately called "approximate functional  equation"} the $L$-value $L(\chi,1/2)$ can be written as a sum of rapidly converging series (cf.~\cite[Theorem 5.3]{IwKo}): for $q$ prime and $\chi\not=1$
$$|L(\chi,1/2)|^2=2\sum_{n_1,n_2\geq 1}\frac{\chi(n_1)\ov\chi(n_2)}{(n_1n_2)^{1/2}}V(\frac{n_1n_2}{q})$$
where $V$ is a rapidly decreasing function which depends on $\chi$ only through its parity $\chi(-1)=\pm 1$. Plugging this expression in the second moment \eqref{2ndchi} and unfolding, one finds that  the key point is to obtain a bound of the following shape\footnote{for simplicity we ignore the dependency of $V$ in the parity of the $\chi$'s}
\begin{equation}\label{KhNggoal}
\sumsum_\stacksum{\ell_1,\ell_2\leq L,n_1,n_2}{(l_1l_2n_1n_2,q)=1}\frac{x_{l_1}\ov{x_{l_2}}}{(ql_1l_2n_1n_2)^{1/2}}V(\frac{n_1n_2}q)e\left(\frac{n_2\ov{l_1l_2n_1}}{q}\right)\ll q^{-\delta}	
\end{equation}
for some $\delta=\delta(\lambda)>0$ for any fixed $\lambda<3/10$. This sum can  be decomposed in various sub-sums in which the variables are localized to specific ranges. The problem becomes essentially that of bounding by $O(q^{-\delta})$ the family of bilinear sums
$$\Sigma(L_1,L_2,N_1,N_2)=\frac{1}{(qL_1L_2N_1N_2)^{1/2}}\sumsum_\stacksum{l_i\sim L_i,i=1,2}{n_1,n_2}x_{l_1}\ov{x_{l_2}}W(\frac{n_1}{N_1})W(\frac{n_2}{N_2})e\left(\frac{n_2\ov{l_1l_2n_1}}{q}\right)$$
where $W\in\mcC_c(]1/2,2[)$, $L_1,L_2\leq L$ and $N_1N_2\leq q$.

The $n_2$-sum is essentially a geometric series bounded by
$$\ll \min(N_2,{\|\ov{l_1l_2n_1}/q\|^{-1}})$$ 
where $\|\cdot\|$ is the distance to the nearest integer. Hence
\begin{align}
\nonumber 
\Sigma(L_1,L_2,N_1,N_2)&\ll \frac{q^\eps}{(qL_1L_2N_1N_2)^{1/2}}\sum_{m\approx L_1L_2N_1}\min(N_2,{\|\ov{m}/q\|^{-1}})\\
&\ll\frac{q^{2\eps}}{(qL_1L_2N_1N_2)^{1/2}}\max_{1\leq U\leq q/2}\min(N_2,\frac{q}{U})\sumsum_\stacksum{m\approx L_1L_2N_1,\ ,\ u\sim U}{um\equiv \pm 1\mods q}1\nonumber\\
\nonumber
&\ll\frac{q^{2\eps}}{(qL_1L_2N_1N_2)^{1/2}}\max_{1\leq U\leq q/2}\min(N_2,\frac{q}{U})(\frac{L_1L_2N_1U}{q}+1)\\
&\ll
q^{2\eps}\frac{L}{q^{1/2}}(\frac{N_1}{N_2})^{1/2}.\label{KNbound1}
\end{align}
(Observe that for $\frac{L_1L_2N_1U}{q}\ll 1$ the equation $um\equiv\pm 1\mods q$ has no solution unless $L_1L_2N_1U\ll 1$).

Alternatively, applying the Poisson summation formula to the $n_1$ variable we obtain a sum of the shape 
$$\Sigma(L_1,L_2,N_1,N_2)=\frac{1}{(qL_1L_2N_1N_2)^{1/2}}\frac{N_1}{q^{1/2}}\sumsum_\stacksum{l_i\sim L_i,i=1,2}{n_1,n_2}x_{l_1}\ov{x_{l_2}}\widetilde W(\frac{n_1}{q/N_1})W(\frac{n_2}{N_2})\Kl_2(\ov{l_1l_2}n_1n_2;q)$$
where $\widetilde W$ is bounded and rapidly decreasing. Bounding this sum trivially (using that $|\Kl_2(m;q)|\leq 2$) yields
\begin{equation}\label{trivialbound12}
\Sigma(L_1,L_2,N_1,N_2)\ll q^{\eps}L(\frac{N_2}{N_1})^{1/2}.	
\end{equation}
The expression $\min(\frac{L}{q^{1/2}}(\frac{N_1}{N_2})^{1/2},L(\frac{N_2}{N_1})^{1/2})$ is maximal for $\frac{N_1}{N_2}=q^{1/2}$ and equals $L/q^{1/4}$ which is $O(q^{-\delta})$ if $\lambda<1/4$.
 
 The bound \eqref{trivialbound12} did not exploit cancellation from the $n_1,n_2,l_1,l_2$ averaging and indeed this is not evident because in the limiting case $N_1=q^{3/4},\ N_2=q/N_1=q^{1/4}$, $L_1=L_2=L=q^{1/4}$, one has $$n_1\approx n_2\approx l_1\approx l_2\approx q^{1/4}$$ which is pretty short. Nevertheless Khan and Ngo where able to detect further cancellation from summing of these short variables. The idea, which we have met already, is to group some of these variables to form longer variables. One possibility could be to group together $n_1$, $n_2$ on the one hand and $l_1$, $l_2$ on the other hand with the idea of applying the methods of \S \ref{secbilinear}. However, the new variables would have size $q^{1/2}$, which is the P\'olya-Vinogradov range at which point the standard completion method just fails. Instead, one can group $n_1$, $n_2$ and $l_2$ together and leave $l_1$ alone. The variable $r=n_1n_2\ov{l_2}\mods q$ takes essentially $q^{3/4}$ distinct values but over all of $\Fqt$ and does not vary along an interval. To counter this defect, one uses the Holder inequality instead of Cauchy-Schwarz.
 
Proceeding as above, we write
 $$\Sigma(L_1,L_2,N_1,N_2)=\frac{1}{(qL_1L_2N_1N_2)^{1/2}}\frac{N_1}{q^{1/2}}\sumsum_{r\in\Fqt,l_1}x_{l_1}\nu(r)\Kl_2(\ov{l_1}r;q)$$
 where
 $$\nu(r)=\sumsum_\stacksum{l_2,n_1,n_2}{r=n_1n_2\ov{l_2}(q)}\ov{x_{l_2}}\widetilde W(\frac{n_1}{q/N_1})W(\frac{n_2}{N_2}).$$

Under the assumption
\begin{equation}\label{moment2hyp}
	L_2\frac{q}{N_1}N_2<q/100\Longrightarrow  L_2\frac{N_2}{N_1}<1/100
\end{equation}
we have
 $$\sum_r|\nu(r)|+\sum_r|\nu(r)|^2\ll q^\eps L_2\frac{q}{N_1}N_2.$$
 Indeed under \eqref{moment2hyp} one has
 $$\ov l_2n_1n_2\equiv \ov l'_2n'_1n'_2\mods q \Longleftrightarrow  l'_2n_1n_2\equiv l_2n'_1n'_2\mods q \Longleftrightarrow l'_2n_1n_2= l_2n'_1n'_2$$
 and the choice of $l'_2,n_1,n_2$ determines $l_2,n'_1,n'_2$ up to $O(q^\eps)$ possibilities.
Hence, applying Cauchy's inequality twice, we obtain 
 $$\Sigma(L_1,L_2,N_1,N_2)=\frac{q^\eps}{(qL_1L_2N_1N_2)^{1/2}}\frac{N_1}{q^{1/2}}(L_2\frac{q}{N_1}N_2)^{3/4}\left(\sum_{r\in\Fqt}|\sum_{l\sim L_1}x_l\Kl_2(\ov lr;q)|^4\right)^{1/4}.$$
 
Now (using that $\Kl_2(n;q)\in\Rr$)
$$\sum_{r\in\Fqt}|\sum_{l\sim L_1}x_l\Kl_2(\ov lr;q)|^4\ll q^\eps\sum_\bfl|\sum_{r\in\Fqt}\prod_{i=1}^4\Kl_2(\ov{l_i}r;q)|$$
where $\bfl=(l_1,l_2,l_3,l_4)\in[L_1,2L_1[^4$.

Theorem \ref{cor-concrete}, applied to the Kloosterman sheaf, gives
$$\sum_{r\in\Fqt}\prod_{i=1}^4\Kl_2(\ov{l_i}r;q)\ll q^{1/2}$$
unless there exists a partition $\{1,2,3,4\}=\{i,j\}\sqcup\{k,l\}$ such that 
$$l_i=l_j,\ l_k=l_l.$$ 
In this case, we use the trivial bound
 $$\sum_{r\in\Fqt}\prod_{i=1}^4\Kl_2(\ov{l_i}r;q)\ll q.$$
 Hence
$$\sum_\bfl|\sum_{r\in\Fqt}\prod_{i=1}^4\Kl_2(\ov{l_i}r;q)|\ll L_1^2q+L_1^4q^{1/2}$$
and
\begin{align}\nonumber
\Sigma(L_1,L_2,N_1,N_2)&\ll\frac{q^\eps}{(qL_1L_2N_1N_2)^{1/2}}\frac{N_1}{q^{1/2}}(L_2\frac{q}{N_1}N_2)^{3/4}(L_1^{1/2}q^{1/4}+L_1q^{1/8})\\
&\ll q^{\eps}L(\frac{N_2}{N_1})^{1/2}(Lq\frac{N_2}{N_1})^{-1/4}(L^{-1/2}q^{1/4}+q^{1/8}).	\label{KhNgbound3}
\end{align}
For $L\geq q^{1/4}$ (the range one would like to improve) one obtains under \eqref{moment2hyp} 
\begin{equation}\label{KhNgbound3final}
	\Sigma(L_1,L_2,N_1,N_2)
\ll q^{\eps}L(\frac{N_2}{N_1})^{1/2}(Lq^{1/2}\frac{N_2}{N_1})^{-1/4}.
\end{equation}

Suppose now we are in a limiting case for \eqref{trivialbound12}, namely
$L^2N_2/N_1=1$. Then \eqref{moment2hyp} holds as long as $L\gg 1$ and \eqref{KhNgbound3final} improves over \eqref{trivialbound12} by a factor $(q^{1/2}/L)^{1/4}$, which is $<1$ as long as $L< q^{1/2}$. 

A more detailed analysis combining \eqref{KNbound1}, \eqref{trivialbound12} and \eqref{KhNgbound3final} shows that \eqref{KhNggoal} holds for any fixed $\lambda<3/10$, and hence leads to Theorem \ref{KhNgthm}.

\section{Advanced completion methods: the $q$-van der Corput method}
In this section and the next ones, we discuss general methods to evaluate trace functions along intervals of length smaller than the P\'olya-Vinogradov range discussed in \S \ref{Secshort}.

\subsection{The $q$-van der Corput method} One of the most basic techniques encountered in analytic number to estimate sums of (analytic) exponentials is the {\em van der Corput method} (see \cite[Chap. 8]{IwKo}). The $q$-Van der Corput method is an arithmetic variant due to Heath-Brown which replace archimedean analysis with $q$-adic analysis. That method concerns $c$-periodic functions for $c$ a {\em composite number}. Suppose (to simplify the presentation) that $c=pq$ for two primes $p$ and $q$ and let
$$K_c=K_pK_q\colon \Zz/c\Zz\to \Cc$$ be some function modulo $c$ which is the product of two trace functions modulo $p$ and $q$ (of conductor bounded by some constant $C$). We consider the sum
$$S_V(K,N):=\sum_{n}K_c(n)V(\frac nN)=\sum_{n}K_p(n\mods p)K_q(n\mods q)V(\frac nN)$$
where $V\in\mcC^\infty(]1,2[)$ and $2N<c=pq$.
We will explain the proof of the following result
\begin{theorem}[$q$-van der Corput method] Let $c=pq$ a product of two primes and $K_c=K_p.K_q$ as above; assume that $K_q$ is the trace function associated with a geometrically irreducible sheaf $\mcF$, which is not geometrically isomorphic to a linear or quadratic phase (i.e.~not of the shape $[P]^*\mcL_\psi$ for $P$ a polynomial of degre $\leq 2$). Then for $2N<pq$, we have
$$S_V(K_c,N)\ll_C N^{1/2}(p+q^{1/2})^{1/2}.$$
\end{theorem}
\begin{remark}This bound is non trivial as long as $$N\geq \max(p,q^{1/2}),$$
which is a weaker condition than $N\geq (pq)^{1/2}$ as long as 
$$1<p<q.$$
We have therefore improved over the P\'olya-Vinogradov range; moreover the range of non triviality is maximal when $p\approx c^{1/3}$ and $q\approx c^{2/3}$. In that case,  one obtains
\begin{equation}\label{optimalvdC}
S_V(K,N)\ll_C N^{1/2}c^{1/6}	
\end{equation}
which is non-trivial as long as
$$N\geq c^{1/3}.$$
\end{remark}

\proof The proof makes use of the (semi-)invariance of $K$ under  translations:
$$K(n+ph)=K_p(n)K_q(n+ph).$$
For $H\leq N/100p$ we have
$$S_V(K,N)=\frac{1}{2H+1}\sum_{|h|\leq H}\sum_{n}K_p(n)K_q(n+ph)V(\frac {n+ph}N)
$$
$$= \frac{1}{2H+1}\sum_{|n|\leq 3N}K_p(n)\sum_{|h|\leq H}K_q(n+ph)V(\frac {n+ph}N)$$ 
$$\ll \frac{1}{2H+1}N^{1/2}\bigl(\sum_{|n|\leq 3N}\bigl|\sum_{|h|\leq H}K_q(n+ph)V(\frac {n+ph}N)\bigr|^2\bigr)^{1/2}$$
$$=\frac{N^{1/2}}{H}\bigl(\sumsum_{|h|,|h'|\leq H}\sum_{n}K_q(n+ph)\ov{K_q(n+ph')}W_{p,h,h'}(\frac{n}N)\bigr)^{1/2}$$
where
$$W_{p,h,h'}(\frac{n}N)=V(\frac {n+ph}N)\ov{V(\frac {n+ph'}N)}.$$
We split the $h,h'$-sum into its diagonal and non-diagonal contribution
$$\sumsum_{|h|,|h'|\leq H}\ldots=\sumsum_\stacksum{|h|,|h'|\leq H}{h=h'}\ldots+\sumsum_\stacksum{|h|,|h'|\leq H}{h\not=h'}\ldots\ .$$
The diagonal sum contributes by $O(NH)$ and it remains to consider the correlation sums
$$\mcC(K_q,h,h'):=\sum_{n}K_q(n+ph)\ov{K_q(n+ph')}W_{p,h,h'}(\frac{n}N)$$
for $h\not=h'$.

Observe that this is the sum of a trace function of modulus $q$ of length $\approx N$. By comparison with the initial sum, we had a trace function of modulus ${pq}$ of length $\approx N$ so the relative length of $n$ compared to the modulus has increased ! By the P\'olya-Vinogradov method, it is sufficient to determine whether the sheaf
$$[+ph]^*\mcF\otimes [+ph']^*D(\mcF)$$ 
has an Artin-Schreier sheaf in its irreducible components. This is equivalent to whether one has an isomorphism
$$[+p(h-h')]^*\mcF\simeq \mcF\otimes\mcL_\psi$$
for some Artin-Schreier sheaf. We will answer this question in a slighly more general form:
\begin{definition} For $d$ an integer satisfying  $1\leq d<q$, a polynomial phase sheaf of degree $d$ is a sheaf of the shape
$[P]^*\mcL_\psi$ for $P$ a polynomial of degree $d$ and $\psi$ a non-trivial additive character. It is lisse on $\Aa^1_{\Fq}$, ramified at infinity with Swan conductor equal to $d$ and its trace function equals $$x\mapsto \psi(P(x)).$$
	
\end{definition}
 We can now invoke the following
\begin{proposition}[\cite{Polymath8a}] Let $d$ be an integer satisfying  $1\leq d<q$. Suppose that $\mcF$ is geometrically irreducible, not isomorphic to a polynomial phase of degree $\leq d$  and that $\CF\leq q^{1/2}$. Then for any $h\in\Fq-\{0\}$ and any non-constant polynomial $P$ of degree $\leq d-1$,$$[+h]^*\mcF\hbox{ and } \mcF\otimes [P]^*\mcL_\psi$$
are not geometrically isomorphic.
	
\end{proposition}
\proof We will only give the easiest part of it and refer to \cite[Thm. 6.15]{Polymath8a} for the complete argument. Suppose that $\mcF$ is ramified at some point $x_0\in\Aa^1(\ov{\Fq})$, since polynomial phases are ramified only at $\infty$ the isomorphism $$[+h]^*\mcF\simeq \mcF\otimes [P]^*\mcL_\psi$$
restricted to the inertia group $I_x$
 implies that $\mcF$ is ramified at $x_0-h$ and iterating at $x_0-nh$ for any $n\in\Zz$, this would imply that $C(\mcF)\geq q$ which is excluded. It remains to deal with the case where $\mcF$ is ramified only at $\infty$.
\qed

Under our assumptions the above proposition implies that for $h\not=h'$
$$\mcC(K_q,h,h')=O(q^{1/2})$$ and that
$$S_V(K,N)\ll N^{1/2}(\frac{N}{H}+q^{1/2})^{1/2}$$
and we choose $H=N/100 p$ to conclude the proof.

\qed

\subsection{Iterating the method} Suppose more generally that $c$ is a squarefree number and that 
$$K_c=\prod_{q|c}K_q$$
is a product of trace functions associated to sheaves not containing any polynomial phases. One can repeat the above argument after factoring $c$ into a product of squarefree coprime moduli $r.s$ and decompose accordingly $$K_c=K_r.K_s.$$
Thus, we  have to bound sums of the shape
\begin{equation}\label{ssum}
	\sum_{n}K_s(n+rh)\ov{K_s(n+rh')}W_{r,h,h'}(\frac{n}{N})
\end{equation}
 This time we need to be a bit more careful and decompose the $h,h'$ sum according to the gcd $(h-h',s)$. After applying the Poisson summation formula (cf.~\eqref{eqpoisson}) we can factor the resulting Fourier transform modulo $s$ into sums over prime moduli $q|s$:
$$\what{K_s}(y)=\prod_{q|s}\what{K_q}(\ov{s_q}y\mods q),\ y\in\Zz/s\Zz,\ s_q=s/q.$$
If $q|h-h'$ we use the trivial bound $\what{K_q}(\ov{s_q}y\mods q)\ll q^{1/2}$ and if $q\!\!\not|h-h'$ we use the non-trivial bound $\what{K_q}(\ov{s_q}y\mods q)\ll 1$. We eventually obtain (see \cite{Polymath8a})
\begin{theorem} Let $C\geq 1$, let $c$ be squarefree and let $K_c\colon \Zz/c\Zz\ra\Cc$ be a product of trace functions $K_q$  such that for any prime $q|c$ the underlying sheaf $\mcF_q$ is of conductor $\leq C$ , is geometrically irreducible and is not geometrically isomorphic to any polynomial phase of degree $\leq 2$. Then 
$$S_V(K_c,N)\ll_{C,\eps} c^{\eps}N^{1/2}(r+s^{1/2})^{1/2}$$
	for any $\eps>0$.
\end{theorem}
If $s$ is not a prime, we could also iterate, factor $s$ into $s=r_2s_2$ and instead of applying the P\'olya-Vinogradov completion method to the sum \eqref{ssum}, we could also apply the $q$-van der Corput method
with the trace functions
$$n\mapsto K_q(n+rh)\ov{K_q(n+rh')},\ q|s_1.$$
This leads us to the quadruple correlation sum
$$\mcC(K_q,{\mathbf \gamma},\alpha)=\frac{1}q\sum_{x}K_q(\gamma_1\cdot x)K_q(\gamma_2\cdot x)\ov{K_q(\gamma'_1\cdot x)K_q(\gamma'_2\cdot x)}\eq(\alpha x)$$
where the $\gamma_i,\gamma'_j,\ i,j=i,2$ are unipotent matrices
$$\gamma_i=\begin{pmatrix}1&h_i\\0&1
\end{pmatrix},\ \gamma'_i=\begin{pmatrix}1&h'_j\\0&1
\end{pmatrix}
$$
In suitable situations, we can then apply Theorem \ref{cor-concrete} from the previous section.

An important example is when 
$$K_c(n)=\Kl_k(n;c)=\frac{1}{c^{(k-1)/2}}\sum_\stacksum{x_1,\ldots,x_k\in(\Zz/c\Zz)^\times}{x_1.\ldots.x_k=n}e\left(\frac{x_1+\ldots+x_k}{c}\right)$$ is a hyper-Kloosterman sum. For any $q|c$, one has $$K_q(y)=\Kl_k(\ov{c_q}^ky;q)\text{ with }\ c_q=c/q $$ and the underlying sheaf is the multiplicatively shifted Kloosterman sheaf $\mcF_q=[\times \ov{c_q}^k]^*\KL_k$. In that case Theorem \ref{cor-concrete} applies and we eventually obtain the bound
$$S_V(\Kl_k(\cdot;c),N)\ll_k c^{\eps}N^{1/2}\left(r+(N^{1/2}(s_1+s_2^{1/2}))^{1/2}\right)^{1/2}.$$
for any factorisation $c=rs_1s_2$. In particular, if there exists a  factorisation $c=rs_1s_2$ such that
$$r\approx c^{1/4},\ s_1\approx c^{1/4},\ s_2\approx c^{1/2}$$ we obtain 
$$S_V(\Kl_k(\cdot;c),N)\ll_k N^{1-\eta}$$
for some $\eta=\eta(\delta)>0$ as long as $$N\geq c^{1/4+\delta}.$$
Iterating once more we see that for any factorisation $c=rs_1s_2s_3$ one has
\begin{equation}\label{vdCiterate3}
S_V(\Kl_k(\cdot;c),N)\ll_{k,\eps} c^{\eps}N^{1/2}\left(r+(N^{1/2}(s_1+(N^{1/2}(s_2+s_3^{1/2}))^{1/2}))^{1/2}\right)^{1/2}	
\end{equation}
so if there exists a  factorisation $c=rs_1s_2s_3$ such that
$$r\approx c^{1/5},\ s_1\approx c^{1/5},\ s_2\approx c^{1/5},\ s_3\approx c^{2/5}$$ then
$$S_V(\Kl_k(\cdot;c),N)\ll_{k,\eps} N^{1-\eta}$$
for some $\eta=\eta(\delta)>0$   as long as $$N\geq c^{1/5+\delta}.$$
We can continue this way as long as enough factorisation for $c$ are available. Such availability is garanteed by the notion of friability:
\begin{definition} An integer $c\not=0$ is $\Delta$-friable if $$q|c\ (q\hbox{ prime })\Rightarrow q\leq \Delta.$$
\end{definition}
Using the reasoning above, Irving \cite{IrvingIMRN} proved  the following  result for $k=2$ (in a quantitative form):
\begin{theorem} For any $L\geq 2$ there exists $l=l(L)\geq 1$ and $\eta=\eta(L)>0$ such that for $c$ a squarefree integer which is $c^{1/l}$-friable and any $k\geq 2$, one has,
$$S_V(\Kl_k(\cdot;c),N)\ll_{k,V} N^{1-\eta}$$
whenever $N\geq c^{1/L}$.
\end{theorem}
Therefore one can obtain non-trivial bounds for extremely short sums of hyper-Kloosterman sums as long as their modulus is firable enough. In particular for $k=2$ we have seen in Remark \ref{remselberg} that improving on Selberg's $2/3$-exponent for the distribution of the divisor function in arithmetic progressions to large moduli (Theorem \ref{thmd2}) was essentially equivalent to bounding non-trivially sums of the shape
$$\sumsum_{n_1,n_2}\Kl_2(an_1n_2;c)V(\frac{n_1}{N^*_1})V(\frac{n_2}{N^*_2})$$
for $$N^*_1N^*_2\approx c^{1/2}.$$
If $N^*_1N^*_2\approx c^{1/2}$ then $\max(N^*_1,N^*_2)\gg c^{1/4}$ and we can use the \eqref{vdCiterate3} to bound non-trivially the above sum granted that $c$ is friable enough. This leads to the following theorem (compare with Theorem \ref{thmd2} for $c$ a prime):
\begin{theorem}{\cite{IrvingIMRN}} There exists $L\geq 4$ and $\eta>0$ such that for any $c\geq 1$ which is squarefree and $c^{1/L}$-friable and any $a$ coprime with $c$, one has for $c\leq X^{2/3+\eta}$ and any $A\geq 0$
$$E(d_2;c,a)\ll_A \frac{X}{c}(\log X)^{-A}.$$ 
\end{theorem}
See \cite{IrvingIMRN2} and \cite{WuPing} for further applications of these ideas.

\section{Around Zhang's theorem on bounded gaps between primes }

Some of the arguments of the previous chapter can be found in Yitang Zhang's spectacular proof of the existence of bounded gaps between the primes:
\begin{theorem}[\cite{YZhang}]\label{thmZhang} Let $(p_n)_{n\geq 1}$ be the sequence of primes in increasing order ($p_1=2,p_2=3,p_3=5,\ldots$). There exists an absolute constant $C$ such that
$$p_{n+1}-p_n\leq C$$
for infinitely many $n$.
\end{theorem}
Besides Zhang's original paper, we refer to \cite{Gran,KowBBK1} for a detailed description of Zhang's proof and the methods involved and historical background. Let us however mention a few important facts:
\begin{itemize}
\item The question of the existence of small gaps between primes has occupied analytic number theorists for a very long time and has been the motivations for the invention of many techniques, in particular the {\em sieve method} to detect primes with additional constraints. A conceptual breakthrough occurred with the work of Goldston, Pintz and Y\i ld\i r\i m \cite{GYP} who proved the weaker result
$$\liminf_{n}\frac{p_{n+1}-p_n}{\log p_n}=0$$
and who on this occasion invented a technique which is also key to Zhang's approach (see Soundararajan's account of their works \cite{SoundGYP}).
\item Zhang's theorem can be seen as an approximation to the twin prime conjecture:
$$\hbox{\em{There exist infinitely many primes $p$ such that $p+2$ is prime}}.$$
Indeed, Zhang's theorem with $C=2$ is equivalent to the twin prime conjecture.
\item A value for the constant $C$ can be given explicitly : Zhang himself gave 
$$C=70.10^6$$
and mentioned that this could certainly be improved. Improving the value of this constant was the objective of the Polymath8 project: following and optimizing Zhang's method in several aspects (some to be explained below), the value was reduced to
$$C=4680.$$
However Maynard \cite{Maynard} made independently another conceptual breakthrough, simplifying the whole proof and making it possible to obtain stronger results and improving the constant to 	
$$C=600.$$
Eventually the Polymath8 project joined with Maynard ; optimizing his argument, the value
$$C=246$$
 was reached (cf.~\cite{Polymath8b}).
 
A side-effect of Maynard's approach is that what we are going to describe now plays no role anymore in this specific application. Nevertheless, it adresses another important question in analytic number theory.
\end{itemize}

\subsection{The Bombieri-Vinogradov theorem and beyond}
The breakthrough of Goldston, Pintz and Y\i ld\i r\i m that is at the origin of Zhang's work builds on the use of sieve methods to detect the existence of infinitely many pairs of primes at distance $\leq C$ from one another. The fuel to be put in this sieve machine are results concerning the distribution of primes in arithmetic progressions to moduli large with respect to the size of the primes which are sought after. In this respect the Bombieri-Vinogradov theorem already discussed in \S \ref{Secternary} is a powerful substitute to GRH: 
\begin{theorem}[Bombieri-Vinogradov] For any $A>0$ there is $B=B(A)>0$ such that for $x\geq 2$
$$\sum_{q\leq x^{1/2}/\log^B x}\max_{(a,q)=1}\left|\psi(x;q,a)-\frac{\psi(x;q)}{\vphi(q)}\right|\ll \frac{x}{\log^Ax}.$$
\end{theorem}

For the question of the existence of bounded gaps between primes, the exponent $1/2$ appearing in the constraint $q\leq x^{1/2}/\log^B x$ turns out to be crucial. In their seminal work \cite{GPY}, Goldston-Pintz-Y\i ld\i r\i m had  pointed out that the Bombieri-Vinogradov theorem with the exponent $1/2$ replaced by any strictly larger constant would be sufficient to imply Theorem \ref{thmZhang}. 

The possibility of going beyond Bombieri-Vinogradov is not unexpected: the Elliott-Halberstam conjecture predicts that any fixed exponent $<1$ could replace $1/2$. That this conjecture is not wishful thinking comes from the work of Fouvry, Iwaniec and Bombieri-Friedlander-Iwaniec from the 80's \cite{FIActaAr, Fou,BFI} who proved versions of the Bombieri-Vinogradov theorem  with exponents $>1/2$ but for "fixed" congruences classes (for instance with the sum involving the difference $|\psi(x;q,1)-\frac{\psi(x;q)}{\vphi(q)}|$ instead of $\max_{(a,q)=1}|\psi(x;q,a)-\frac{\psi(x;q)}{\vphi(q)}|$). Zhang's groundbreaking insight has been to nail down a beyond-Bombieri-Vinogradov type theorem that could be established unconditionally and would be sufficient to establish the existence of bounded gaps between primes. The following theorem is a variant of Zhang's  theorem (\cite[Thm 1.1]{Polymath8a}). Let us recall that an integer $q\geq 1$ is $\Delta$-friable if any prime $p$ dividing $q$ is $\leq \Delta$.
\begin{theorem}\label{thmzhang} Let $\bfa=(a_p)_{p\in\mcP}$ be a sequence of integers indexed by the primes such that  $a_p$ is coprime with $p$ for all $p$. For any squarefree integer $q$, let $a_q\mods q$ be the unique congruence class modulo $q$ such that $$\forall p|q,\ a_q\equiv a_p\mods p;$$
 in particular $a_q\in(\Zz/ q\Zz)^\times$.
There exist absolute constants $\theta>1/2$ and $\delta>0$, independent of $\bfa$, such that for any $A>0$, $x>2$ one has
$$\sum_\stacksum{q\leq x^{\theta},\ {sqfree}}{q\ x^\delta-{friable}}|\psi(x;q,a_q)-\frac{\psi(x;q)}{\vphi(q)}|\ll \frac{x}{\log^Ax}.$$
Here the implicit constant depends only on $A$, but not on  $\bfa$.
\end{theorem}
\begin{remark} Zhang essentially proved this theorem for $\theta=1/2+1/585$ and in an effort to improve Zhang's constant, the Polymath8 project improved $1/585$ to $7/301$.
\end{remark}
We will now describe some of the principles of the proof of this theorem and especially at the points where algebraic exponential sums occur. We refer to the introduction of \cite{Polymath8a} and to E. Kowalski's account in the Bourbaki seminar \cite{KowBBK1}. 

Let us write $c(q)$ for $\mu^2(q)$ times the sign of the difference $\psi(x;q,a_q)-\frac{\psi(x;q)}{\vphi(q)}$. The above sum equals
$$\sum_\stacksum{q\leq x^{\theta}}{q\ x^\delta-\hbox{friable}}c(q)\sum_{n\leq x}\Lambda(n)\Delta_{\bfa}(n;q).$$
where
$$\Delta_{\bfa}(n):=\delta_{n\equiv a_q\mods q}-\frac{\delta_{(n,q)=1}}{\vphi(q)}$$
As is usual when counting primes numbers, the next step is to decompose the von Mangoldt function $\Lambda(n)$ into a sum of convolution of arithmetic functions (for instance by using Heath-Brown's identity Lemma \ref{lemHB} as in \S \ref{Secprimes}): we essentially arrive at the problem of bounding $(\log x)^{O_J(1)}$ of the following model sums (for $j\leq J$ and $J$ is a fixed and large integer)
$$\Sigma(\bfM;\bfa,Q):=\sum_\stacksum{q\sim Q}{q\ x^\delta-\hbox{friable}}c(q)\sumsum_{m_1,\ldots, m_{2j}}\mu(m_1)\ldots\mu(m_j) V_{1}	\Bigl(\frac{m_{1}}{M_{1}}\Bigr)\ldots V_{2j}	\Bigl(\frac{m_{2j}}{M_{2j}}\Bigr)\Delta_{a_q}(m_1\ldots m_{2j})	$$
where $V_i,\ i=1,\ldots, 2j$ are smooth  functions compactly supported in $]1,2[$ and $\bfM=(M_1,\ldots,M_{2j})$ is a tuple satisfying
$$Q\leq x^{\theta},\ M_i=:x^{\mu_i},\ \forall i\leq j,\ \mu_i\leq 1/J,\ \sum_{i\leq 2j}{\mu_i}=1+o(1).$$

Our target is the bound
\begin{equation}\label{zhangtarget}
\Sigma(\bfM;\bfa,Q)
\overset{?}{\ll} \frac{x}{\log^A x}.
\end{equation}
The most important case is when $$Q= x^\theta=x^{1/2+\varpi}$$ for some fixed sufficiently small $\varpi>0$.

The variables with index $j+1\leq \leq 2j$ are called {\em smooth} because they are weighted by smooth functions and this makes it possible to use the Poisson summation formula on them to analyze the congruence condition mod $q$. This is going to be efficient if the range $M_i$ is sufficiently big relatively to $q\sim Q$.  The variables with indices $1\leq i\leq j$ are weighted by the M\"obius function but (at least as long as some strong form of the Generalized Riemann Hypothesis is not available) we cannot exploit this information and we will consider the M\"obius functions like arbitrary bounded functions. The tradeoff to non-smoothness is that the range of these variables is pretty short $M_i\leq x^{1/J}$, especially if $J$ is choosen large.

 As we did before we will aggregate  some of the variables $m_i,\ i=1,\ldots,2j$ so as to form two new variables whose ranges are located adequately (similarly to what we did in \S \ref{Secprimes}) and will use different methods to bound the sums depending on the size and the type of these new variables.

More precisely, we define 
$$\alpha_i(m)=\begin{cases}\mu(m)V_{i}	\Bigl(\frac{m}{M_{i}}\Bigr)&\ 1\leq i\leq j\\
V_{i}\Bigl(\frac{m}{M_{i}}\Bigr)&\ j+1\leq i\leq 2j	.
\end{cases}
$$
Given some partition of the set of $m$-indices $$\{1,\ldots,2j\}=\bfI\sqcup\bfJ$$
let
$$M=\prod_{i\in\bfI}M_i,\ N=\prod_{j\in\bfJ}M_j$$
and
$$\mu_\bfI:=\sum_{i\in\bfI}\mu_i,\ \mu_\bfJ:=\sum_{i\in\bfJ}\mu_i.$$
We have $$ \mu_\bfI+\mu_\bfJ=1+o(1),\ M=x^{\mu_{\bfI}},\ N=x^{\mu_{\bfJ}}.$$
In the sequel we will always make the convention that $N\leq M$ or equivalently  $\mu_\bfI\geq\mu_\bfJ$.

Finally we define the Dirichlet convolution functions
$$\alpha(m):=\star_{i\in\bfI}\alpha_i(m),\ \beta(n):=\star_{i\in\bfJ}\alpha_i(n).$$
We are reduced to bound sums of the shape
\begin{equation}\label{partitionsum}
\sum_\stacksum{q\sim Q}{x^\delta-\text {friable}}c(q)\sumsum_\stacksum{m\sim M}{n\sim N}\alpha(m)\beta(n)\Delta_{a_{q}}(mn)
\overset{?}{\ll} \frac{x}{\log^A x}.	
\end{equation}

Observe that the functions $\alpha,\beta$ are essentially bounded
$$\forall\eps>0, \alpha(m),\beta(n)\ll x^\eps$$
so we need only to improve slightly over the trivial bound.

\subsection{Splitting into types}
The sums \eqref{partitionsum} will be subdivided into three different types and their treatment will depend on which type the sum belong.

This subdivision follows from the following simple combinatorial Lemma (cf.~\cite[Lem. 3.1]{Polymath8a}):
\begin{lemma} Let $1/10<\sigma<1/2$ and let $\mu_i,\ i=1,\ldots 2j$ be some non-negative real numbers such that
$$\sum_{i=1}^{2j}\mu_i=1.$$
One of the following holds
\begin{itemize}
\item Type 0: there exists $i$ such that $\mu_i\geq 1/2+\sigma$.
\item Type II: there exists a partition $$\{1,\ldots,2j\}=\bfI\sqcup\bfJ$$ such that
$$1/2-\sigma\leq \sum_{i\in\bfJ}\mu_i\leq \sum_{i\in\bfI}\mu_i<1/2+\sigma.$$
\item Type III: there exist distincts $i_1,i_2,i_3$	such that
$$2\sigma\leq \mu_{i_1}\leq \mu_{i_2}\leq \mu_{i_3}\leq 1/2-\sigma\hbox{ and }\mu_{i_1}+\mu_{i_2}\geq 1/2+\sigma.$$
\end{itemize}
\end{lemma}

\begin{remark}\label{remTypeIII} If $\sigma>1/6$ the Type III
situation never occurs since $2\sigma>1/2-\sigma$.	
\end{remark}

Given $\sigma$ such that  $$1/10<\sigma<1/2$$ 
we assume that $J$ is choosen large enough so that
\begin{equation}\label{Jcond}
1/J\leq \min(1/2-\sigma,\sigma).	
\end{equation}
We say that a sum \eqref{partitionsum} is of
\begin{itemize}
	
\item Type 0, if there exists some $i_0$ such that $\mu_{i_0}\geq 1/2+\sigma$. We choose
$$\bfI=\{i_0\}\hbox{  and $\bfJ$ the complement.}$$ Since for any $i\leq j$, one has $\mu_i\leq 1/J<1/2+\sigma$, necessarily $i_0\geq j+1$ corresponds to a smooth variable; the corresponding sum therefore equals
\begin{equation}
\label{type0}\sum_\stacksum{q\sim Q}{x^\delta-\hbox{friable}}c(q)\sumsum_{m\geq 1, n\sim N}V(\frac{m}{M_{i_0}})\beta(n)\Delta_{a_{q}}(mn).	
\end{equation}

\item Type I/II if one can partition the set of indices
$$\{1,\ldots,2j\}=\bfI\sqcup\bfJ$$ in a way that the corresponding ranges
$$M=\prod_{i\in\bfI}M_i=x^{\mu_\bfI}\geq N=\prod_{i\in\bfJ}M_{i}=x^{\mu_\bfJ}$$
satisfy
\begin{equation}\label{NcondtypeI}
{1/2-\sigma}\leq \mu_\bfJ=\sum_{i\in\bfJ}\mu_{i}\leq {1/2}
\end{equation}

\item Type III if we are neither in the Type 0 or Type I/II situation: there exist distinct indices $i_1,i_2,i_3$	such that
$$2\sigma\leq \mu_{i_1}\leq \mu_{i_2}\leq \mu_{i_3}\leq 1/2-\sigma\hbox{ and }\mu_{i_1}+\mu_{i_2}\geq 1/2+\sigma.$$
We choose
$$\bfI=\{i_1,i_2,i_3\}\hbox{  and $\bfJ$ to be the complement.}$$
Again, since $1/J<2\sigma$ by \eqref{Jcond}, the indices $i_1,i_2,i_3$ are associated to smooth variables and the Type III sums are of the shape
$$\sum_\stacksum{q\sim Q}{x^\delta-\hbox{friable}}c(q)\sumsum_\stacksum{m_1,m_2,m_3}{n\sim N}V(\frac{m_1}{M_{i_1}})
V(\frac{m_2}{M_{i_2}})V(\frac{m_3}{M_{i_3}})\beta(n)\Delta_{a_{q}}(m_1m_2m_3n).$$
\end{itemize}

\begin{remark}In the paper \cite{Polymath8a} the "Type II" sums introduced here were split into two further types that were called "Type I" and "Type II". These are the sums for which the $N$ variable satisfies
\begin{gather*}
\hbox{ Type I: }x^{1/2-\sigma}\leq N< x^{1/2-\varpi-c} \\
\hbox{ Type II: }x^{1/2-\varpi-c}\leq N\leq x^{1/2}	
\end{gather*}
for some extra parameter $c$ satisfying
$$1/2-\sigma<1/2-\varpi-c<1/2.$$
This distinction was necessary for optimisation purposes and especially to achieve the exponent $1/2+7/301$ in Theorem \ref{thmzhang}. 
	
\end{remark}

Zhang's Theorem now essentially follows from
\begin{theorem} There exist $\varpi,\sigma>0$ with $1/10<\sigma<1/2$ such that the bound \eqref{partitionsum} holds for the Type 0, II and III sums.
\end{theorem}
For the rest of this section we will succinctly describe how each type of sum is handled.

The case of Type 0 sums \eqref{type0} is immediate: one applies the Poisson summation formula to the $m$ variable to decompose the congruence $mn\equiv a_q\mods q$. The zero frequency contribution is cancelled up to an error term  by the second term of $\Delta_{a_{q}}(mn)$ while the non-zero frequencies contribute a negligible error term as long as the range of the $m$ variable is larger than the modulus, i.e.~
$$1/2+\sigma>1/2+\varpi$$ which can be assumed. 
 
\subsection{Treatment of type II sums}

\subsubsection{The art of applying Cauchy-Schwarz}
The  Type II sums are more complicated to deal with because we have essentially no control on the shape of the coefficients $\alpha(m),\beta(n)$ (except that they are being essentially bounded). The basic principle is to consider the largest variable $m\sim M$, to make it smooth using the Cauchy-Schwarz inequality and then resolve the congruence 
$$m\equiv \ov n a_q\mods q$$
using the Poisson summation formula. This is the essence of the {\em dispersion method} of Linnik.

When implementing this strategy one has to decide which variables to put "inside" the Cauchy-Schwarz inequality and which to leave "outside". To be more specific, suppose we need to bound a general trilinear sum 
$$\sumsum_{m\sim M,n\sim N}\sum_{q\sim Q}\alpha_m\beta_n\gamma_q K(m,n,q)$$
and wish to smooth the $m$ variable using Cauchy-Schwarz. There are two possibilities, either
$$
\sumsum_{m\sim M,n\sim N}\sum_{q\sim Q}\alpha_m\beta_n\gamma_q K(m,n,q)\ll \|\alpha\|_2\|\gamma\|_2\biggl(\sumsum_{m\sim M,q\sim Q}|\sum_{n\sim N}\beta_n K(m,n,q)|^2\biggr)^{1/2}$$
or
$$
\sumsum_{m\sim M,n\sim N}\sum_{q\sim Q}\alpha_m\beta_n\gamma_q K(m,n,q)\ll \|\alpha\|_2\biggl(\sum_{m\sim M}|\sumsum_{n\sim N,q\sim Q}\beta_n\gamma_q K(m,n,q)|^2\biggr)^{1/2}	$$

In the first case the inner sum of the second factor equals
$$\sumsum_{n_1,n_2\sim N}\beta_{n_1}\ov{\beta_{n_2}}\sumsum_{m\sim M,q\sim Q} K(m,n_1,q)\ov{K(m,n_2,q)}
$$
and in the second case
$$\sumsum_{n_1,n_2\sim N}\sumsum_{q_1,q_2\sim Q}\beta_{n_1}\gamma_{q_1}\ov{\beta_{n_2}\gamma_{q_2}}\sum_{m\sim M} K(m,n_1,q_1)\ov{K(m,n_2,q_2)}.$$

In either case, one expects to be able to detect cancellation from the $m$-sum, at least when the other variables $(n_1,n_2)$ or $(n_1,n_2,q_1,q_2)$ 
are not located on the diagonal (i.e.~$n_1=n_2$ or $n_1=n_2,\ q_1=q_2$). If the other variables are on the diagonal, no cancellation is possible but the diagonal is small compared to the space of variables.

We are faced with the following trade-off:
\begin{itemize}
	
\item For the first possibility, the $m$-sum is simpler (it involves three parameters $n_1,n_2,q$) but the ratio "size of the diagonal"$/$"size of the set of parameters" is  $N/N^2=N^{-1}$. 	
\item For the second possibility, the $m$-sum is more complicated as it involves more auxiliary parameters $n_1,n_2,q_1,q_2$ but the ratio "size of the diagonal"$/$"size of the set of parameters" $NQ/N^2Q^2=1/NQ$ is smaller (hence more saving can be obtained from the diagonal part). 
\end{itemize}

\subsubsection{The Type II sums}

We illustrate this discussion in the case of Type II sums. If we apply Cauchy with the $q$ variable outside the diagonal $n_1=n_2$ would not provide enough saving. If, on the other hand, we apply Cauchy with $q$ inside, then the diagonal is large but we have to analyze the congruence
$$mn_1\equiv a\mods{q_1},\ mn_2\equiv a\mods{q_2}$$
which is a congruence modulo $[q_1,q_2]$. Assuming we are in the generic case of $q_1,q_2$ coprime, the resulting modulus is $q_1q_2\sim Q^2=x^{1+2\varpi}$ while $m\sim M\leq x^{1/2}$, which is too small for the Poisson formula to be efficient. 

There is fortunately a middle-ground: we can use the extra flexibility (due to Zhang's wonderful insight) that our problem involves {\em friable} moduli: by the greedy algorithm, one can factor  $q\sim Q$ into a product $q=rs$ where $r$ and $s\sim Q/r$ vary over ranges that we can essentially choose as we wish (up to a small  indeterminacy of $x^\delta$ for $\delta$ small). In other words, we are reduced to bounding sums of the shape
$$\Sigma(M,N;\bfa,R,S)=\sumsum_\stacksum{r\sim R,\ s\sim S}{rs\ x^\delta-\hbox{friable}}c(rs)\sumsum_\stacksum{m\sim M}{n\sim N}\alpha(m)\beta(n)\Delta_{a_{rs}}(mn)$$
for any factorisation $RS=Q$ that fits with our needs. Now, when applying Cauchy-Schwarz, we have the extra flexibility of having the $r$ variable "out" and the $s$ variable "in".

 We do this and get
\begin{multline*}
\sumsum_{r\sim R,s\sim S}c(rs)\sumsum_\stacksum{m\sim M}{n\sim N}\alpha(m)\beta(n)\Delta_{a_{rs}}(mn)=\sum_{r\sim R}\sum_{m\sim M}\alpha(m)\sum_sc(rs)\sum_{n\sim N}\beta(n)\Delta_{a_{rs}}(mn)\\
\ll_\eps R^{1/2}M^{1/2+\eps}\biggl(\sum_r\sumsum_{s_1,s_2,n_1,n_2}c(rs_1)\ov{c(rs_2)}\beta(n_1)\ov{\beta(n_2)}\sum_{m}V(\frac{m}M)\Delta_{a_{rs_1}}(mn_1)\Delta_{a_{rs_2}}(mn_2)\biggr)^{1/2}	
\end{multline*}
for $V$ a smooth function compactly supported in $[M/4,4M]$. We choose $R$ of the shape $$R=Nx^{-\eps}\leq Mx^{-\eps}$$
for $\eps>0$ but small.

Expanding the square, we obtain a sum involving four terms. The most important one comes from the product
\begin{equation}\label{Deltaproduct}
\Delta_{a_{rs_1}}(mn_1)\Delta_{a_{rs_2}}(mn_2)=(\delta_{mn_1\equiv a_{rs_1}\mods{rs_1}}-\frac{\delta_{(n,rs_1)=1}}{\vphi(rs_1)})(\delta_{mn_2\equiv a_{rs_2}\mods {rs_2}}-\frac{\delta_{(n,rs_2)=1}}{\vphi(rs_2)})	.
\end{equation}
We will concentrate on the contribution of this term from now on. 

The generic and main case is when $(s_1,s_2)=1$, so that $m$ satisfies a congruence modulo $rs_1s_2\sim RS^2=Mx^{2\varpi+\eps}$ which is not much larger than $M$ if $\varpi$ is small. Observe that 
$$mn_i\equiv a_{rs_i}\mods{rs_i},\ i=1,2\Longrightarrow \ n_1\equiv n_2\mods r.$$
We can therefore write $n_1=n,\ n_2=n+rl$ with $|l|\ll N/R=x^{\eps}$.
By the Poisson summation formula, we have 
$$\sum_{m}V(\frac{m}M)\delta_{m\equiv b\mods{rs_1s_2}}=\frac{M}{rs_1s_2}\what V(0)+\frac{M}{rs_1s_2}\sum_{h\not=0}\what V(\frac{h}{rs_1s_2/M})e\left(\frac{hb}{rs_1s_2}\right)$$
where $b=b(n,l)\mods {rs_1s_2}$ is such that
$$b\equiv a_{rs_1s_2}\ov n\mods r,\ b\equiv a_{rs_1s_2}\ov n\mods{s_1}, b\equiv a_{rs_1s_2}\ov{n+lr}\mods{s_2}.$$
The $h=0$ contribution provides a main term which is cancelled up to an admissible error term by the main contributions coming from the other summands of \eqref{Deltaproduct}. The contribution of the frequencies $h\not=0$ will turn out to be error terms. We have to show that
$$\sum_r\sumsum_{s_1,s_2,n,l}c(rs_1)\ov{c(rs_2)}\beta(n)\ov{\beta(n+rl)}\frac{M}{rs_1s_2}\sum_{h\not=0}\what V(\frac{h}{rs_1s_2/M})e\left(\frac{hb}{rs_1s_2}\right){\ll} \frac{MN^2}Rx^{-\eta}=x^{1-\eta+\eps}$$
for some fixed $\eta>0.$
 The length of the $h$ sum is essentially
$$H=RS^2/M=Q^2N/(xR)=x^{2\varpi+\eps}$$
which is small (if $\varpi$ and $\eps$ are). We therefore essentially need to prove that
\begin{multline}\label{TypeIIinter}\frac{1}{H}\sum_{r\sim R}\sum_{l\ll N/R}\sum_{n}\beta(n)\ov{\beta(n+lr)}\sum_{0\not=h\ll H}\left|\sum_{s_1,s_2}c(rs_1)\ov{c(rs_2)}e\left(h\frac{a_{rs_1s_2}\ov n}{rs_1}+h\frac{a_{rs_1s_2}\ov{n+lr}}{rs_2}\right)\right|\\
{\ll}x^{1-\eta+\eps}.	
\end{multline}

We can now exhibit cancellation in the $n$-sum by smoothing out the $n$ variable using the Cauchy-Schwarz inequality for any fixed $r,l$: letting the $h$ variable "in" we obtain exponential sums of the shape

$$\sum_{n\sim N}e\left(h\frac{a_{rs_1s_2}\ov n}{rs_1}-h'\frac{a_{rs'_1s'_2}\ov n}{rs'_1}+h\frac{a_{rs_1s_2}\ov{n+lr}}{rs_2}-
h'\frac{a_{rs'_1s'_2}\ov{n+lr}}{rs'_2}\right).$$

The generic case is when $h-h',s_1,s_2,s'_1,s'_2$ are all coprime. In that case
the above exponential sum has length $$N\in[x^{1/2-\sigma},x^{1/2}]$$ and the moduli involved are of size $$RS^4=Q^4/R^3=x^{O(\eps)}Q^4/N^3=[x^{1/2+4\varpi+O(\eps)},x^{1/2++4\varpi+3\sigma+O(\eps)}].$$
Therefore if $\sigma,\varpi,\eps$ are small, the length $N$ is not much smaller than the modulus so we could apply the completion method to improve over the trivial bound $O(N)$ for the $n$-sum. If we apply the P\'olya-Vinogradov method, the trivial bound is replaced by $O((RS^4)^{1/2+o(1)})$  and we find that  the left-hand side of \eqref{TypeIIinter} is bounded by 
$$\frac{1}HR.\frac{N}RN^{1/2}(H^2S^4(RS^4)^{1/2+o(1)})^{1/2}=x^{O(\eps)+o(1)}N^{3/2}S^3R^{1/4}=x^{\frac{7}8+3\varpi+\frac{5}4\sigma+O(\eps)+o(1)}$$
which is $\ll x^{1-\eta}$ for some $\eta>0$ whenever $\sigma<1/10$ and $\varpi$ and $\eps$ are small enough.

Instead of using the P\'olya-Vinogradov bound, we could take advantage of the fact that the modulus $rs_1s'_1s_2s'_2$ is $x^\delta$-friable (again we can take $\delta>0$ as small as we need) and apply the $q$-van der Corput method from the previous section. Factoring $rs_1s'_1s_2s'_2$ into a product $r's'$ such that $r'\sim (rs_1s'_1s_2s'_2)^{1/3+O(\delta)}$, $s'\sim (rs_1s'_1s_2s'_2)^{2/3+O(\delta)}$, a suitable variant of
\eqref{optimalvdC} bounds the $n$-sum by $O(N^{1/2}(RS^4)^{1/6+O(\delta)+o(1)})$ and the 
left-hand side of \eqref{TypeIIinter} is bounded by 
$$\frac{R}H\frac{N}RN^{\frac{1}2}(H^2S^4N^{1/2}(RS^4)^{1/6})^{\frac{1}2+o(1)+O(\delta)}=x^{O(\eps+\delta)+o(1)}N^{7/4}S^{7/3}R^{1/12}=x^{\frac{11}{12}+\frac73\varpi+\frac{1}2\sigma+O(\eps+\delta)+o(1)}$$
which is $\ll x^{1-\eta}$ for some $\eta>0$ whenever $\sigma<1/6$ and $\varpi$ and $\eps$ are small enough.

\subsection{Treatment of type III sums}
Our objective for the Type III sums is the following bound: for some $\eta>0$, we have
\begin{equation}\label{typeIIIgoal}
\sum_\stacksum{q\sim Q}{x^\delta-\hbox{friable}}c(q)\sum_{n\sim N}\beta(n)\sum_{m}\tau_{3,\bfM}(m)\Delta_{a_{q}}(m_1m_2m_3n){\ll}x^{1-\eta},	
\end{equation}
where $\bfM=(M_{i_1},M_{i_2},M_{i_3})$ and
$$\tau_{3,\bfM}(m):=\sum_{m_1m_2m_3=m}V(\frac{m_1}{M_{i_1}})
V(\frac{m_2}{M_{i_2}})V(\frac{m_3}{M_{i_3}})$$
and $M_{i_1},M_{i_2},M_{i_3}$ satisfy 
$$M=M_{i_1}M_{i_2}M_{i_3}\geq x^{1/2+3\sigma}.$$

The function 
$$m\mapsto \tau_{3,\bfM}(m)$$
is basically a smoothed version of the ternary divisor function $m\mapsto \tau_3(m)$ that we have discussed in \S \ref{Secternary}.

In fact, while describing the proof of Theorem \ref{thmd3}, we have shown that for $M=x$, and for $q$ a prime satisfying 
$$q\sim x^{1/2+\varpi},\ \varpi=1/47$$ one has
$$ \sum_{m}\tau_{3,\bfM}(m)\Delta_{a_{q}}(m_1m_2m_3n)\ll \frac{x^{1-\eta}}q$$ for some $\eta>0$. We have therefore the required bound but for individual moduli instead of having it on average.

As we have observed when discussing Type II sums, the parameter $\sigma$ can be taken as close to $1/6$ as we wish and in particular
 $M\in [x^{1+3(\sigma-\frac16)},x]$ can be made as close as we wish from $x$ and $N\in [1,x^{3(\frac16-\sigma)}]$ as we wish from $x$ (in the logarithmic scale). In particular, this establishes  \eqref{typeIIIgoal} for prime moduli $q\sim Q$ for some value of $\sigma$ (close enough to $1/6$), and some value of $\varpi$ (close enough to $0$) and some $\eta>0$.
 
 The case of $x^\delta$-friable moduli uses similar methods and (besides some elementary technical issues) is maybe simpler than in the prime modulus case because of the extra flexibility provided by the friable moduli.  

\begin{remark} By a more elaborate treatment, involving different uses of the Cauchy-Schwarz inequality and iterations of the $q$-van der Corput method, it is possible to bounds successfully all the Type II sums associated to some explicit parameter $\sigma>1/6$. As pointed out in Remark \ref{remTypeIII}, this makes the section devoted to Type III sums (and in particular the theory of hyper-Kloosterman sums $\Kl_3(x;q)$) unnecessary. The interest of this remark comes from the fact that the trace functions occurring in the treatment of the sums of Type II are exclusively algebraic exponentials:
$$x\mapsto \eq({f(x)}),\ \hbox{for }f(X)\in\Fq(X).$$
For such trace functions, Corollary \ref{delignecor} "only" uses Weil's resolution of the Riemann Hypothesis for curves over finite fields \cite{Weil0} and not the full proof of the Weil conjectures by Deligne \cite{WeilII}.
\end{remark}

\section{Advanced completions methods: the  $+ab$ shift}
In this last section, we describe another method allowing to break the P\'olya-Vinogradov barrier for prime moduli. This method has its origins in the celebrated work of Burgess on short sums of Dirichlet characters \cite{Bur}.

\subsection{Burgess's bound}\label{Bursec}

Let $q$ be a prime and le $\chi\colon \Fqt\ra\Ct$ be a non trivial multiplicative character. Consider the sum
$$S_V(\chi,N):=\sum_{n}\chi(n)V(\frac nN)$$
where $V\in\mcC^\infty(]1,2[)$. 

\begin{theorem}[Burgess]\label{Burgessbound} For any $N\geq 1$ and $l\geq 1$ such that 
\begin{equation}\label{Burgesscond}
q^{1/2l}\leq N< \frac12q^{1/2+1/4l}	
\end{equation}
we have
$$S_V(\chi,N)\ll_{V,l}q^{o(1)}N(N/q^{1/4+1/4l})^{-1/l}.$$	
\end{theorem}
\begin{remark} Observe that this bound is non-trivial (sharper than $S_V(\chi,N)\ll N$) whenever
$$q^{1/4+1/4l+o(1)}\leq N< \frac12q^{1/2+1/4l}.$$
Moreover, for $N\geq  \frac12q^{1/2+1/4l}$, the P\'olya-Vinogradov bound $S_V(\chi,N)\ll q^{1/2}$ is non trivial, therefore, we see that by taking $l$ large enough, that \eqref{Burgessbound} yields a non-trivial bound for $S_V(\chi,N)$ as long as $$N\geq q^{1/4+\delta}$$  for some fixed $\delta>0$.
\end{remark}

\proof
Burgess's argument exploits two features in a critical way: the first one is that an interval is "essentially" invariant under sufficiently small additive translations and the second is the multiplicativity of the Dirichlet character.

Let $A,B\geq 1$ be parameters such that 
$AB\leq N/2$; we will also assume that $2B<q$. 

We have 
$$S_V(\chi,N)=\frac{1}{AB}\sum_{|n|\leq 2N}\sumsum_{a\sim A,b\sim B}\chi(n+ab)V(\frac{n+ab}N).$$
The next step is to invoke the Fourier inversion formula to separate the variables $n$ and $ab$: one has
$$V(\frac{n+ab}N)=\int_\Rr\what V(t)e(\frac{tn}N)e(\frac{tab}N)dt.$$
Plugging this formula in our sum, we obtain 
\begin{align*}
S_V(\chi,N)
&=\frac{1}{AB}\int_\Rr\sum_{|n|\leq 2N}e(\frac{tn}N)\sumsum_{a\sim A,b\sim B}\chi(n+ab)e(\frac{tab}N)\what V(t)dt	\\
&\leq  \frac{1}{AB}\int_\Rr\sum_{|n|\leq 2N}\sum_{a\sim A}\bigl|\frac{\chi(a)}{a}\what V(\frac{t}a)\bigr|\bigl|\sum_{b\sim B}\chi(\ov an+b)e(\frac{tb}N)\bigr|
dt\\
&\leq \frac{1}{AB}\int_\Rr \sum_{|n|\leq 2N}\sum_{a\sim A}\bigl|\sum_{b\sim B}\chi(\ov an+b)e(\frac{tAb}N)\bigr||W(t)|dt
\end{align*}
for $W$ some bounded rapidly decaying function. 
\begin{remark}\label{remchi}
	Observe that the factor $\chi(a)$ coming from the identity
\begin{equation}\label{chimult}
\chi(n+ab)=\chi(a(\ov an+b))=\chi(a)\chi(\ov an+b)	
\end{equation}
has been absorbed in the absolute value of the first inequality above.
\end{remark}

The innermost sum can be rewritten
$$\sum_{|n|\leq 2N}\sum_{a\sim A}\bigl|\sum_{b\sim B}\chi(\ov an+b)e(\frac{tAb}N)\bigr|=\sum_{r\in\Fqt}\nu(x)|\sum_{b\sim B}\eta_b\chi(r+b)\bigr|$$
where $\eta_b=e(\frac{tAb}N)$ and
$$\nu(r):=|\{(a,n)\in[A,2A[\times[-2N,2N],\ \ov an=r\mods q\}|.$$
Consider the map
$$(a,n)\in[A,2A[\times[-2N,2N]\mapsto \ov an\mods q=r\in \Fq.$$
The function $\nu(r)$ is the size of the fiber of that map above $r$. We will show that this map is "essentially injective" (has small fibers on average). Suppose that $A$ is chosen such that
$4AN<q;$ then one has
$$\sum_r\nu(r)\ll AN,\ \sum_r\nu^2(r)\ll (AN)^{1+o(1)}$$
where the first bound is obvious while for the second we observe that
$$\sum_r\nu^2(r)=|\{(a,a',n,n'), \ a,a'\in[A,2A[,\ |n|,|n'|\ll N,\ an'\equiv an\mods q\}|,$$
then use the fact that $AN<q$ and that the integer $an'$ has at most $(an')^{o(1)}$ decomposition of the shape $an'=a'n$.

This map however is not surjective nor even close to being so in general, so that the change of variable $\ov a.n\leftrightarrow x$ is not very effective. A way to moderate ineffectiveness is to use H\"older's inequality.

Let $l\geq 1$ be some integer parameter. Applying H\"older's inequality with $1/p=1-1/2l,\ 1/q=1/2l$ and the above estimate one obtains
\begin{align*}
\sum_{x\in\Fqt}\nu(x)|\sum_{b\sim B}\eta_b\chi(x+b)\bigr|&\leq  	(\sum_x\nu(x)^{\frac{2l}{2l-1}})^{1-1/2l}(\sum_x
|\sum_{b\sim B}\eta_b\chi(x+b)\bigr|^{2l})^{1/2l}\\
&\ll (AN)^{1-1/2l+o(1)}(\sum_x
|\sum_{b\sim B}\eta_b\chi(x+b)\bigr|^{2l})^{1/2l}.
\end{align*}
The $x$-sum in the rightmost factor equals
$$\sum_{\bfb}\eta_\bfb\sum_{r\in\Fq}\chi(\frac{\prod_{i=1}^l(r+b_i)}{\prod_{i=i}^l(r+b_{k+i})})$$
where
$\bfb=(b_1,\ldots,b_{2l})\in[B,2B[^{2l}$ and $\eta_\bfb=\prod_{i=1}^{2l}\eta_{b_i}$. Consider the fraction
$$F_\bfb(X):=\frac{\prod_{i=1}^l(X+b_i)}{\prod_{i=i}^l(X+b_{k+i})}\in \Qq(X)$$
and the function on $\Fq$
$$r\in\Fq\mapsto \chi(F_\bfb(r))$$
(extended by $0$  for $r=-b_i\mods q,\ i=1,\ldots,2l$). This function is the trace function of the rank one sheaf $[F_\bfb]^*\mcL_\chi$ whose conductor is bounded in terms of $l$ only and (because it is of rank $1$) which is geometrically irreducible if not-geometrically constant. If not geometrically constant one has\footnote{It is not necessary to invoke Deligne's main theorem here: this follows from A. Weil's proof of the Riemann hypothesis for curves \cite{Weil0}.}
$$\sum_{r\in\Fq}\chi(F_\bfb(r))\ll_lq^{1/2}.$$
If $q>\max(l,2B)$ this occurs precisely when $F_\bfb(X)$ is not constant nor a $k$-th power, where $k$ is the order of $\chi$. Hence this holds for $\bfb$ outside an explicit set $\mcB^{bad}\subset [B,2B[^{2l}$ of size bounded by $O(B^{l})$. If $\bfb\in \mcB^{bad}$, we use the triv,ial bound
$$|\sum_{r\in\Fq}\chi(F_\bfb(r))|\leq q.$$
All in all, we eventually obtain
$$\sum_{\bfb}\eta_\bfb\sum_{x}\chi\left(\frac{\prod_{i=1}^l(x+b_i)}{\prod_{i=i}^l(x+b_{k+i})}\right)\ll |\mcB^{bad}|q+|\mcB-\mcB^{bad}|q^{1/2}\ll B^lq+B^{2l}q^{1/2}.$$
Choosing $B=q^{1/2l}$ (so as to equal the two terms in the bound above) and $A\approx Nq^{-1/2l}$  with the condition $4AN<q$, which is equivalent to \eqref{Burgesscond}, we obtain that
$$S_V(\chi,N)\ll_l \frac{q^{o(1)}}{AB}(AN)^{1-1/2l}(q^{3/2})^{1/2l}\ll q^{o(1)}N^{1-1/{l}}q^{3/4l-(1-1/2l)/2l}=q^{o(1)}N(N/q^{1/4+1/4l})^{-1/l}.$$
\qed

\subsection{The $+ab$-shift for type I sums}
It is natural to try to extend this method to other trace functions; unfortunately the above argument breaks down because the identity \eqref{chimult} is not valid in general. It is however possible to mitigate this problem by introducing an extra average. 

This technique goes back to Karatsuba and Vinogradov (for the function $x\mapsto \chi(x+1)$). It has been also used by Friedlander-Iwaniec \cite{FrIw} (for the function $x\mapsto e\left(\frac{\ov x}q\right)$), Fouvry-Michel \cite{FoMi} and Kowalski-Michel-Sawin \cite{KMS,KMS2}. 

 Instead of a single sum $S_V(K,N)$, one considers the following average of multiplicative shifts
 $$B_V(K,\uple{\alpha},N):=\sum_{m\sim M}\alpha_m\sum_n V(\frac{n}{N})K(mn)$$
 where $1\leq M<q$ and $(\alpha_m)_{m\sim M}$ is a sequence of complex numbers of modulus $\leq 1$ (this includes the averaged sum
 $\sum_{m\sim M}\bigl|\sum_nK(mn)V(\frac{n}{N})\bigr|=\sum_m|S_V([\times m]^*K,N)|$). The 
  objective here is to improve over the trivial bound
 $$B_V(K,\uple{\alpha},N)\ll \|K\|_\infty MN.$$
 
 Proceeding as above we have
\begin{align*}
B_V(K,\uple{\alpha},N)&=\frac{1}{AB}\sum_{m}\alpha_m\sum_{n}\sumsum_{a\sim A,b\sim B}K(m(n+ab))V(\frac{n+ab}N)\\
&\leq \frac{1}{AB}\int_\Rr \sum_{m\sim M}\alpha_m\sum_{|n|\leq 2N}\sum_{a\sim A}\bigl|\sum_{b\sim B}K(am(\ov an+b))e(\frac{tAb}N)\bigr||W(t)|dt.
\end{align*}
We have
$$\sum_{m\sim M}\alpha_m\sum_{|n|\leq 2N}\sum_{a\sim A}\bigl|\sum_{b\sim B}K(am(\ov an+b))e(\frac{tAb}N)\bigr|=\sumsum_{r,s\in\Fq}\nu(r,s)\bigl|\sum_{b\sim B}\eta_b K(s(r+b))\bigr|$$
with
$$\nu(r,s)=\sum_{m\sim M}\sum_{|n|\leq 2N}\sum_{a\sim A}\alpha_m\delta_{\ov an=r,am=s\mods q}.$$
Assuming that $4AN<q$ and evaluating the number of solutions to the equations
$$am=a'm',\ a\ov n\equiv a'\ov n'\mods q,\ (a,m,n)\in[A,2A[\times[M,2M[\times[N,2N[$$
one finds that
$$\sumsum_{r,s\in\Fq}|\nu(r,s)|\ll AMN,\ \sumsum_{r,s\in\Fq}|\nu(r,s)|^2\ll q^{o(1)}AMN$$
which we interpret as saying that the map
$$(a,m,n)\in[A,2A[\times[M,2M[\times[N,2N[\ra (r,s)=(\ov a.n,am)\in \Fq\times [AM,4AM[$$ is essentially injective (i.e.~has small fibers on average). As before, this map is far from being surjective but one can dampen this with H\"older's inequality:
$$\sumsum_\stacksum{r\in\Fq}{1\leq s\leq 4AM}\nu(r,s)\bigl|\sum_{b\sim B}\eta_b K(s(r+b))\bigr|\ll \big(\sumsum_{r,s}|\nu(r,s)|^{\frac{2l}{2l-1}}\big)^{1-1/2l}\big(\sumsum_{r,s}\bigl|\sum_{b\sim B}\eta_b K(s(r+b))\bigr|^{2l}\big)^{1/2l}$$
$$\ll q^{o(1)}(AMN)^{1-1/2l}\bigl(\sum_{\bfb}\eta_{\bfb}\sum_{r,s}\prod_{i=1}^lK(s(r+b_i))\ov{K(s(r+b_{i+l}))}\bigr)^{1/2l}.$$
We are now reduced to the problem of bounding the two variable sum 
\begin{equation}\label{generalKVsum}
\sum_{r,s}\prod_{i=1}^lK(s(r+b_i))\ov{K(s(r+b_{i+l}))}=\sum_r\sum_{s}\bfK(sr,s\bfb)=\sum_{r}\bfR(r,\bfb)	
\end{equation}
(say) where
\begin{equation}\label{KRdef}
\bfK(r,\bfb):=\prod_{i=1}^lK(r+b_i)\ov{K(r+b_{i+l})},\ \bfR(r,\bfb)=\sum_{s}\bfK(sr,s\bfb).	
\end{equation}

The bound will depend on the vector $\bfb\in[B,2B[^{2l}$. 
To get a feeling of what is going on, let us consider one of cases treated in \cite{FoMi}: let
$$K(x)=\eq(\ov x+x).$$
We have
$$\bfR(sr,s\bfb)=\sum_{s\in\Fqt}\eq(\ov s\sum_{i=1}^l(\ov{r+b_i}-\ov{r+b_{i+l}})+s\sum_{i=1}^l({b_i}-{b_{i+l}})).$$
This sum is either
\begin{enumerate}
\item Equal to $q-1$, if and only if the vector $(b_1,\ldots, b_l)$ equals the vector $(b_{l+1},\ldots, b_{2l})$ up to permutation of the entries.
\item Equal to $-1$ if $\bfb$ is not as in (1) but is in the hyperplane with equation $\sum_{i=1}^l({b_i}-{b_{i+l}})=0$.
\item The Kloosterman sum $$\bfR(r,\bfb)=q^{1/2}\Kl_2\left(\frac{\sum_{i=1}^l(\ov{r+b_i}-\ov{r+b_{i+l}})}{\sum_{i=1}^l({b_i}-{b_{i+l}})};q\right)$$
otherwise.	
\end{enumerate}
The last case is the most interesting. Given $\bfb$ as in the last situation, we have to evaluate
$$q^{1/2}\sum_r\Kl_2(G_\bfb(r);q)$$
where
\begin{equation}\label{Gbdef}
G_\bfb(X)=\frac{\sum_{i=1}^l(\ov{X+b_i}-\ov{X+b_{i+l}})}{\sum_{i=1}^l({b_i}-{b_{i+l}})}.
\end{equation}
\begin{lemma}For $\bfb=(b_1,\ldots,b_{2l})\in\Fq^{2l}$ such that 
\begin{equation}\label{bcond} (b_1,\ldots, b_l)\hbox{ is not equal to }(b_{l+1},\ldots, b_{2l}) \hbox{ up to permutation and }\sum_{i=1}^l({b_i}-{b_{i+l}})\not=0,	
\end{equation}
 one has
	$$\sum_r\Kl_2(G_\bfb(r);q)\ll_l q^{1/2}.$$
\end{lemma}
\proof
The function
$$r\mapsto \Kl_2(G_\bfb(r);q)$$
is the trace function of the rank $2$ sheaf $[G_\bfb]^*\KL_2$ obtained by pull-back from the Kloosterman sheaf $\KL_2$ of morphism
$$x\mapsto G_\bfb(x)$$
which is non-constant by assumption.

Moreover, one can show that he conductor of $[G_\bfb]^*\KL_2$ is bounded in terms of $l$ only, and  moreover the geometric monodromy group of $[G_\bfb]^*\KL_2$ is obtained as the (closure of the) image of the representation $\rho_{\KL_2}$ restricted to a finite index subgroup of $\Gal(\Ksep/\ov{\Fq}.K)$. Since the geometric monodromy group of $\KL_2$ is $\SL_2$ which has no finite index subgroup, the geometric monodromy group of $[G_\bfb]^*\KL_2$ is $\SL_2$ as well. It follows that the sheaf $[G_\bfb]^*\KL_2$ is geometrically irreducible (and not geometrically trivial because of rank $2$) and the estimate follows by Deligne's theorem.
\qed

It follows from this analysis that
$$\sumsum_{r,s}\bigl|\sum_{b\sim B}\eta_b K(s(r+b))\bigr|^{2l}\ll B^{l}q^{2}+B^{2l}q,$$
hence choosing $B=q^{1/l}$, $AB\approx N$ and $A\approx Nq^{-1/l}$ we obtain
$$B_V(K,\uple{\alpha},N)\ll\frac{q^{o(1)}}{AB}(AMN)^{1-1/2l}q^{3/2l}=q^{o(1)}MN(\frac{N^2M}{q^{1+1/l}})^{-1/2l}.$$
To resume we have therefore proven the
\begin{theorem} Let $K(x)=\eq(\ov x+x)$ and $M,N,l\geq 1$ and $(\alpha_m)_{m\sim M}$ be a sequence of complex numbers of modulus bounded by $1$. Assuming that
$$q^{1/l}\leq N<\frac12 q^{1/2+1/2l}$$ we have
$$\sum_{m\sim M}\alpha_m\sum_{n}V(\frac{n}N)K(mn)\ll q^{o(1)}MN(\frac{N^2M}{q^{1+1/l}})^{-1/2l}.$$ 
\end{theorem}

This bound is non trivial (sharper than $\ll MN$) as long as\footnote{If $N\geq \frac12 q^{1/2+1/2l}$ the P\'olya-Vinogradov inequality is non trivial already.}
$$N^2M\geq q^{1+1/l}.$$
For instance, if $M=q^\delta$ for some $\delta>0$,  the above bound is nontrivial for $l$ large enough and $N\geq q^{1/2+\delta/3}$. Alternatively if $M=N$,  this bound is non trivial as long as  
$$N=M\geq q^{1/3+\delta}$$ if $l$ is taken large enough. Therefore this method improves the range of non-triviality in Theorem \ref{thmbilinear}.

\subsection{The $+ab$-shift for type II sums}
With this method, it is also possible to deal with the more general (type II) bilinear sums 
$$B(K,\uple{\alpha},\uple{\beta})=\sumsum_{m\sim M,n\sim N}\alpha_m\beta_n K(mn)$$
where $(\alpha_m)_{m\sim M}$, $(\beta_n)_{n\sim N}$ are sequences of complex numbers of modulus bounded by $1$.

We leave it to the interested reader to fill in the details (or to look at \cite{FoMi,KMS} or \cite{KMS2}). The first step is to apply Cauchy-Schwarz to smooth out the $n$ variable: for a suitable smooth function $V$, compactly supported in $[1/2,5/2]$ and bounded by $1$, one has
$$\bigl|\sumsum_{m\sim M,n\sim N}\alpha_m\beta_n K(mn)\bigr|\leq N^{1/2}\bigl(\sumsum_{m_1,m_2\sim M}\alpha_{m_1}\ov{\alpha_{m_2}}\sum_{n}V(\frac{n}N)K(m_1n)\ov{K(m_2n)}\bigr)^{1/2}.$$
The next step is to perform the $+ab$-shift on the $n$ variable and to make the change of variables
$$(a,m_1,m_2,n)\in[A,2A[\times[M,2M[^2\times[N,2N[\longleftrightarrow (\ov an,am_1,am_2)\mods q=(r,s_1,s_2)\in\Ff_q^3.$$
Considering the fiber counting function for that map, namely
$$\nu(r,s_1,s_2):=\sumsum_\stacksum{(a,n,m_1,m_2)}{a\sim A, |n|\leq 2N,\ m_i\simeq M}\alpha_{m_1}\ov{\alpha_{m_2}}\delta_{\ov an=r,\ am_i=s_i\mods q}$$
one shows that for $AN<q/2$ one has
$$\sumsum_{(r,s_1,s_2)\in\Fq^3}|\nu(r,s_1,s_2)|\ll AM^2N,\ \sumsum_{(r,s_1,s_2)\in\Fq^3}|\nu(r,s_1,s_2)|^2\leq q^{o(1)}AM^2N.$$
Applying H\"older's inequality leads us to the problem of bounding the following complete sum indexed by the parameter $\bfb$ 
\begin{equation}\label{typeIIcomplete}
\sum_{r\in\Fq}|\bfR(r,\bfb)|^2-q\sum_{r\in\Fq}|\bfK(r,\bfb)|^2.	
\end{equation}
We will explain what is expected in general in a short moment but let us see what happens for our previous case $K(x)=\eq(\ov x+x)$: for $\bfb=(b_1,\ldots,b_{2l})\in\Fq^{2l}$ satisfying \eqref{bcond} the sum \eqref{typeIIcomplete} equals
$$q\sum_\stacksum{r\in\Fq}{r\not=-b_i}|\Kl_2(G_\bfb(r);q)|^2-q\sum_\stacksum{r\in\Fq}{r\not=-b_i}1=q\sum_\stacksum{r\in\Fq}{r\not=-b_i}(|\Kl_2(G_\bfb(r);q)|^2-1)+O_l(q)$$
where $G_\bfb(X)$ is defined in \eqref{Gbdef} 
\begin{lemma}For $\bfb=(b_1,\ldots,b_{2l})\in\Fq^{2l}$ satisfying \eqref{bcond}, one has
	$$\sum_r(|\Kl_2(G_\bfb(r);q)|^2-1)\ll_l q^{1/2}.$$
\end{lemma}
\proof This follows from the fact that $[G_\bfb]^*\KL_2$ is geometrically irreducible with geometric monodromy group equal to $\SL_2$: since the tensor product of the standard representation of $\SL_2$ with itself equals the trivial representation plus the symmetric square of the standard representation which is non-trivial and irreducible, 
$$x\mapsto |\Kl_2(G_\bfb(r);q)|^2-1$$
is the trace function of a geometrically irreducible sheaf.\qed

Using this bound and trivial estimates for $\bfb$ not satisfying \eqref{bcond}, one eventually obtains
\begin{theorem}\label{kloosfracprop} Let $K(x)=\eq(\ov x+x)$, $1\leq M,N<q$ and $l\geq 1$ some integer. Assuming that
$$N<\frac12 q^{1/2+1/2l},$$ one has
$$B(K,\uple{\alpha},\uple{\beta})=\sumsum_{m\sim M,n\sim N}\alpha_m\beta_n K(mn)\ll q^{o(1)}MN(\frac1M+(\frac{MN}{q^{3/4+3/4l}})^{-1/4l})^{1/2}.$$ 
\end{theorem}
\begin{remark}
For $l$ large enough, this bound is non-trivial as long as $M\geq q^\delta$ and $MN\geq q^{3/4+\delta}$, again improving on Theorem
\ref{thmbilinear} in this specific case.	
\end{remark}

\subsection{The $+ab$-shift for more general trace functions} 

For applications to analytic number theory, it is highly desirable to extend the method of the previous section to trace functions as general as possible. This method may be axiomatized in the following way.
Let $q$ be a prime, $K\colon \Fq\ra \Cc$ a complex valued function bounded by $1$ in absolute value, $1\leq M,N<q$ some parameters and 
$\uple{\alpha}=(\alpha_m)_{m\sim M}$, $\uple{\beta}=(\beta_n)_{n\sim N}$  sequences of complex number bounded by $1$. We define the type I sum
$$B(K,\uple{\alpha},1_N)=\sumsum_{m\sim M,n\sim N}\alpha_m K(mn)$$
and the type II sum
$$B(K,\uple{\alpha},\uple{\beta})=\sumsum_{m\sim M,n\sim N}\alpha_m\beta_n K(mn).$$

For $l\geq 1$ an integer, let $\bfK(r,\bfb)$ and $\bfR(r,\bfb)$ be the functions of the variables $(r,\bfb)\in\Fq\times\Fq^{2l}$ given by \eqref{KRdef}. For $B\geq 1$ we set 
$$\mcB=\Zz^{2l}\cap[B,2B[^{2l}.$$

An axiomatic treatment of the type I sums $B(K,\uple{\alpha},1_N)$ is provided by the following:

\begin{theorem}\label{thmtypeI} Notations as above, let $B,C\geq 1$ and $\gamma\in [0,2]$ be some real numbers. 
\begin{itemize}
	\item  Let $\mcB^{\Delta}\subset \mcB$ be the set of $\bfb\in\mcB$ for which
	\begin{equation}\label{Rcontrol}
	\hbox{ there exists }r\in\Fq\hbox{ satisfying } |\bfR(r,\bfb)|> Cq^{1/2}.
	\end{equation}

	\item  Let $\mcB_I^{bad}\subset \mcB$ be the union of $\mcB^{\Delta}$ and the set of $\bfb\in\mcB$ such that
		\begin{equation}\label{RScontrol}\bigl|\sum_{r\in\Fq}\bfR(r,\bfb)\bigr|> Cq.
		\end{equation}
\end{itemize}
Suppose that for any $1\leq B<q/2$ one has
\begin{equation}\label{badsetbound1}|\mcB^{\Delta}|\leq CB^l,\ |\mcB_I^{bad}|\leq B^{(2-\gamma)l}	.
\end{equation}

Then, if $N$ satisfies $$q^{1/l}\leq N\leq \frac12q^{1/2+1/2l},$$ one has for any $\eps>0$
\begin{equation}\label{typeIKMS}B(K,\uple{\alpha},1_N)\ll_{C,l,\eps}q^\eps MN(\frac{q^{1+1/l}}{MN^2}+\frac{q^{3/2-\gamma+1/l}}{MN^2})^{1/2l}.	
\end{equation}
\end{theorem}

An axiomatic treatment of the type II sums $B(K,\uple{\alpha},\uple{\beta})$ is provided by the following
\begin{theorem}\label{thmtypeII} Notations as above, let $B,C\geq 1$ and $\gamma\in [0,2]$ be some real numbers,  
\begin{itemize}
	\item  Let $\mcB^{\Delta}\subset \mcB$ be the set of $\bfb\in\mcB$ for which
	$$\hbox{ there exists }r\in\Fq\hbox{ satisfying } |\bfR(r,\bfb)|> Cq^{1/2}.$$
	\item  Let $\mcB_{II}^{bad}\subset \mcB$ be the union of $\mcB^{\Delta}$ and the set of $\bfb\in\mcB$ such that	\begin{equation}\label{RScontrol2}\bigl|\sum_{r\in\Fq}|\bfR(r,\bfb)|^2-q\sum_{r\in\Fq}|\bfK(r,\bfb)|^2\bigr|> C q^{3/2}	.
	\end{equation}

\end{itemize}
Assume that for any  $B\in[1,q/2[$ one has
\begin{equation}\label{badsetbound2}|\mcB^{\Delta}|\leq CB^l,\ |\mcB_{II}^{bad}|\leq CB^{(2-\gamma)l}.	
\end{equation}
Then, if $N$ satisfies $$q^{3/2l}\leq N\leq \frac12q^{1/2+3/4l},$$  one has  for any $\eps>0$,
\begin{equation}\label{typeIIKMS}
B(K,\uple{\alpha},\uple{\beta})\ll_{C,l,\eps}q^\eps MN\bigl(\frac{1}{M}+(\frac{q^{1-\frac{3}4\gamma+\frac{3}{4l}}}{MN}+\frac{q^{\frac{3}{4}+\frac{3}{4l}}}{MN})^{\frac1{l}}\bigr)^{1/2}.	
\end{equation}

\end{theorem}

We conclude these lectures with a few remarks concerning these two theorems:

\begin{enumerate}
\item 	In the case $K(x)=\eq(\ov x+x)$, we have just verified that the conditions \eqref{badsetbound1} and \eqref{badsetbound2} hold with $\gamma=1$. In \cite{FoMi}, this was shown to hold more generally for the trace functions
$$K(x)=\eq(x^{-k}+ax),\ a\in\Fq,\ k\geq 1.$$

\item For more general trace functions, the first condition in \eqref{badsetbound1} and \eqref{badsetbound2} can be verified using some variant of the "sums of products" Theorem \ref{cor-concrete} and does not constitute a big obstacle. One should also notice that Theorem \ref{cor-concrete} implies that for any $\bfb=(b_1,\ldots,b_{2l})$ on the "first" diagonal (i.e.~$b_1=b_{l+1},\ldots,b_{l}=b_{2l}$) one has
$$\bfR(r,\bfb)=\sum_{s}\prod_{i=1}^l|K(s(r+b_i))|^{2}=|K(0)|^{2l}+\sum_{s\not=0}\prod_{i=1}^l|K(s(r+b_i))|^{2}\gg_l q$$
and therefore 
$$|\mcB^\Delta|\geq B^{l}.$$
It follows that the first bound in \eqref{badsetbound1} and \eqref{badsetbound2} is sharp and for the second condition one cannot expect  $\gamma$ to be greater than $1$.
\item In order to reach the best available bound by the above method, it is not necessary to aim for $\gamma=1$: it is sufficient to establish \eqref{badsetbound1} with $\gamma\geq 1/2$ and \eqref{badsetbound2} with $\gamma\geq 1/3$. In such a case, the bounds of Theorem \ref{thmtypeI} and Theorem \ref{thmtypeII} are non trivial as long as
$$MN^2\geq q^{1+1/l},\ MN\geq q^{3/4+3/4l},$$
respectively.

\item Checking the second bound in \eqref{badsetbound1} and \eqref{badsetbound2} for general trace functions is much more difficult. In \cite{KMS}, with specific applications in mind, these bounds have been established for $l=2$ and $\gamma=1/2$  for the  hyper-Kloosterman sums
$$K(x)=\Kl_k(x;q),\ k\geq 2.$$
 Because $l=2$ is too small, this alone is not sufficient to improve over the P\'olya-Vinogradov type bound of Theorem \ref{thmbilinear} (one would have needed $l\geq 4$). A more refined treatment is necessary: instead of letting (somewhat wastefully) the variables $s=am\mods q$ or $s_1=am_1,s_2=am_2 \mods q$ vary freely over the whole interval $[0,q-1]\simeq \Fq$, one uses the fact that $s,s_1,s_2$ belong to the shorter interval $[AM,4AM[$. Applying the P\'olya-Vinogradov completion method to detect this inclusion with additive characters, this leads to bounds for complete sums analogous to \eqref{RScontrol} and \eqref{RScontrol2} but for the additively twisted variant of $\bfR(r,\bfb)$ defined by
 $$\bfR(r,\lambda,\bfb)=\sum_{s}\bfK(sr,s\bfb)e\left(\frac{\lambda s}q\right),\ \hbox{ for }\lambda\in\Fq.$$
Specifically, the bounds are: for all $ \bfb\in\mcB-\mcB^{\Delta},$ we have
 $$\forall\lambda\in\Fq,\ |\bfR(r,\lambda,\bfb)|\leq Cq^{1/2},$$
 and for all $\bfb\in\mcB-\mcB_I^{bad}$, we have
  $$\forall\lambda\in\Fq,\ |\sum_r\bfR(r,\lambda,\bfb)|\leq Cq,$$
  and for all $\bfb\in\mcB-\mcB_{II}^{bad}$, we have
$$\forall\lambda,\lambda'\in\Fq,\ \Bigl|\sum_r\bfR(r,\lambda,\bfb)\ov{\bfR(r,\lambda',\bfb)}-q\delta_{\lambda=\lambda'}\sum_{s}\prod_{i=1}^l|K(s(r+b_i))|^{2}\Bigr|\leq Cq^{3/2}.	
$$
In \cite{KMS}, these bounds were established for $l=2$ and $\bfb$ outside the sets $\mcB^\Delta$, $\mcB_{I}^{bad}$ and $\mcB_{II}^{bad}$
satisfying
$$|\mcB^\Delta|\leq B^2,\ |\mcB_{I,II}^{bad}|\leq CB^{3}.$$

\item In the  paper \cite{KMS2}, the bounds \eqref{badsetbound1} and \eqref{badsetbound2} are established for the  hyper-Kloosterman sums and generalized Kloosterman sums for every $l\geq 2$ and $\gamma=1/2$. 
\end{enumerate}

\subsection{Some applications of the $+ab$-shift bounds}

The problem of estimating bilinear sums of trace functions below the critical P\'olya-Vinogradov range already had several applications in analytic number theory. We list some of them below with references for the interested remaining reader(s).
\begin{itemize}
\item This method was used by Karatsuba and Vinogradov, for the function 
 $$K(n)=\chi(n+a)$$
where $(a,q)=1$ and $\chi\mods q$ is a non-trivial Dirichlet character, to bound non-trivially its sum along the primes over short intervals (now a special case of Theorem \ref{thmprimesumthm}). In particular, Karatsuba \cite{Kar} proved for any $\eps>0$, the bound
$$\sum_\stacksum{p\leq x}{p\ \mathrm{ prime}}\chi(p+a)\ll x^{1-\eps^2/1024}$$
whenever $x\geq q^{1/2+\eps}$. This bound is therefore non-trivial in a range which is wider than that established in Theorem \ref{thmprimesumthm} for general trace functions.
\item The method was used by Friedlander-Iwaniec for the function
$$K(n)=e\left(\frac{\ov n}q\right),\ n.\ov n\equiv 1\mods q$$
to show that the ternary divisor function $d_3(n)$ is well distributed in arithmetic progressions to modulus $q\leq x^{1/2+1/230}$, passing for the first time the Bombieri-Vinogradov barrier (see Theorem \ref{thmd3}).
\item In the case of the Kloosterman sums
 $$K(n)=\Kl_2(n;q),$$
the bound established in \cite{KMS} together with \cite{BlMi,BFKMM} leads to an asymptotic formula for the second moment of character twists of modular $L$-functions: for $f$ a fixed primitive cusp form, one has
$$\frac{1}{q-1}\sum_{\chi\mods q}|L(f\otimes\chi,1/2)|^2= MT_f(\log q)+O_f(q^{-1/145})$$
for $q$ prime, where $MT_f(\log q)$ is a polynomial in $\log q$ (of degree $\leq 1$) depending on $f$. This completes the work of Young for $f$ an Eisenstein series \cite{Young} and of Blomer-Milicevic for $f$ cuspidal and $q$ suitably composite \cite{BlMi}.
\item Using  this method, Nunes \cite{nunes} obtained non-trivial bounds, below the P\'olya-Vinogradov range, for the  (smooth) bilinear sum
$$\sumsum_\stacksum{m\leq M}{n\leq N}K(mn^2)$$
where $K$ is the Kloosterman-like trace function
$$K(n;q):=\frac{1}{q^{1/2}}\sum_{x\in\Fqt}e_q(a\ov x^{2}+bx)$$
(where $a,b$ are some integral parameters such that $(ab,q)=1$). He deduced from this bound that the characteristic function of squarefree integers
is well-distributed in arithmetic progression to prime modulus 
$$q\leq x^{2/3+1/57}.$$
\end{itemize}
The previous best result, due to Prachar \cite{Prachar}, was $q\leq x^{2/3-\eps}$ (similar to Selberg's Theorem \ref{thmd2} for the divisor function $d_2(n)$) dated to 1958 !

\end{document}